\documentclass[11pt,a4paper]{amsart}

\usepackage{microtype}
\usepackage{mathtools}     
\usepackage{amssymb}       
\usepackage{amsthm}        
\usepackage{tensor}
\usepackage{bbm}           

\usepackage{graphicx}
\graphicspath{{./images/}}
\usepackage{tikz}
\usetikzlibrary{calc, patterns} 

\usetikzlibrary{matrix}
\usepackage{tikz-cd}

\usepackage[shortlabels]{enumitem}

\usepackage[hidelinks,hypertexnames=false]{hyperref}
\usepackage{autonum}

\usepackage{float}

\theoremstyle{plain}
\newtheorem{theorem}{Theorem}[section]
\newtheorem{lemma}[theorem]{Lemma}
\newtheorem{prop}[theorem]{Proposition}
\newtheorem{corollary}[theorem]{Corollary}

\theoremstyle{definition}
\newtheorem{definition}[theorem]{Definition}

\newtheorem{ntn}[theorem]{Notation}

\theoremstyle{remark}
\newtheorem{remark}[theorem]{Remark}
\newtheorem{example}[theorem]{Example}

\newcommand{\bb}[1]{\mathbb{#1}}
\newcommand{\cc}[1]{\mathcal{#1}}
\newcommand{\vp}{\varphi}

\numberwithin{equation}{section}

\newcommand{\vlg}{\cc{V}}
\newcommand{\tns}[2]{\tensor{#1}{#2}}
\newcommand{\Cd}{\check C}
\newcommand{\coll}{\kappa}
\newcommand{\lcol}{\nu}
\newcommand{\cob}{\cc{C}ob}

\title{Composing Lagrangian cobordisms without slowing down}
\author[A. Kapota]{Asaf Kapota}
\address{Institute of Mathematics, The Hebrew University of Jerusalem}
\email{asaf.kapota@mail.huji.ac.il}
\author[J. Solomon]{Jake P. Solomon}
\address{Institute of Mathematics, The Hebrew University of Jerusalem}
\email{jake@math.huji.ac.il}
\keywords{Lagrangian cobordism, Lagrangian correspondence, $\infty$-groupoid, Kan complex, Lagrangian boundary conditions, manifolds with corners}
\subjclass[2020]{53D12 (Primary) 55U10, 18N60, 57R90 (Secondary)}
\date{November 2026}

\begin{document}

\begin{abstract}

We establish a framework for Lagrangian cobordism that expresses the condition to be a cobordism between two Lagrangian submanifolds by imposing Lagrangian boundary conditions on the cobordism.
This framework allows Lagrangian cobordisms that are geodesics in the space of positive Lagrangians. Such geodesic cobordisms cannot be trivialized at their ends as required in other approaches to Lagrangian cobordism.  The framework includes the definition of a possibly immersed $k$-cobordism over an arbitrary $k$-manifold with corners.
A Lagrangian $1$-cobordism is the composition of two others if there exists a $2$-cobordism over the $2$-simplex that restricts over the edges to the three given $1$-cobordisms. We show that compositions exist and are unique up to a $2$-cobordism, which is unique up to a $3$-cobordism and so on. More precisely, we show that Lagrangian cobordisms form an $\infty$-groupoid.

\end{abstract}
\maketitle

\tableofcontents

\section{Introduction}
\label{sec:Introduction}

\subsection{Overview}
\label{subsec:Overview}

Let $(X,\omega_X)$ be a symplectic manifold, and let $L_0,L_1 \subset X$ be Lagrangian submanifolds. The notion of a Lagrangian cobordism $W \subset X \times T^\ast [0,1]$ between $L_0$ and $L_1$ has been well studied, tracing back to Arnold~\cite{Arnold1, Arnold2}. More recent works~\cite{BiranCornea1, BiranCornea2,NadlerTanaka} relate the Fukaya category of $X$ to a cobordism category whose objects are Lagrangian submanifolds, morphisms are Lagrangian cobordisms, and composition is concatenation.
A fundamental assumption in these works is that a Lagrangian cobordism $W$ is trivial near its boundary: Identify $ [0,1] $ with the zero section in $ T^\ast [0,1] $. Then,  for some $\varepsilon>0$,
\begin{equation}
	\label{eq:trivialCollars}
	W|_{X \times T^\ast [0,\varepsilon)} = L_0 \times [0,\varepsilon),
		\qquad
		W|_{X \times T^\ast (1-\varepsilon,1]} = L_1 \times (1-\varepsilon,1].
\end{equation}
While this assumption simplifies the composition of Lagrangian cobordisms, reducing it to concatenation, it is somewhat restrictive. In particular, let $\phi_t: X \to X$ be a Hamiltonian isotopy, and suppose a cobordism $W: L \to \phi_1(L)$ arises as the suspension of the path of Lagrangians $t \mapsto \phi_t(L)$. The condition that the cotangent fiber coordinate remains constant on each end forces the generating Hamiltonian $H_t: X \to \bb{R}$ to be constant on $\phi_t(L)$ for $t \in [0,\varepsilon) \cup (1-\varepsilon,1]$.
Our motivation for this work stems from the study of geodesics in the space of positive Lagrangians~\cite{Solomon13}, which are paths of Lagrangians satisfying a geometric PDE that forces the path to be constant if the Hamiltonian $H_t$ is constant for any one $t.$ The geodesic equation can be expressed in terms of the suspension cobordism corresponding to the path~\cite{RubinsteinSolomon17}, and thus one can define a more general notion of geodesic cobordism. Geodesic cobordisms arise naturally in the Lagrangian web construction of~\cite{KS26a}. Thus, we need a more flexible framework for Lagrangian cobordisms.

Dropping the trivial-end assumption~\eqref{eq:trivialCollars} leads to two problems. First, we need to redefine what it means for $W$ to be a cobordism from $L_0$ to $L_1.$ We do this by imposing a Lagrangian boundary condition on the Lagrangian with boundary $ W. $
Specifically, think of a Lagrangian submanifold as a Lagrangian over the point, and a Lagrangian cobordism as a Lagrangian over the interval.
We define a notion of \emph{pushforward} and \emph{pullback} of Lagrangian submanifolds over base manifolds by maps between the base manifolds.
Then $W$ is a cobordism $L_0 \to L_1$ if, for each $i \in \{0,1\}$, the boundary component of $W$ over $\{i\} \subset \partial [0,1] $ lies in the pushforward of $ L_i $ along the inclusion map $\{i\} \hookrightarrow [0,1] ,$ which is the Lagrangian
\begin{equation}
	L_i \times T_i^*[0,1] \subset  X \times T^\ast [0,1].
\end{equation}
Equivalently, the pullback of $W$ along the inclusion $\{i\} \to [0,1]$ is $L_i.$
Moreover, we allow $ L_0 $ and $ L_1 $ themselves to have boundary lying in a fixed embedded Lagrangian $ \vlg \subset X $. Accordingly, $ W $ acquires a boundary component lying in a Lagrangian cobordism from $ \vlg $ to $ \vlg $ given by the pullback of $ \vlg $, which is $\vlg \times [0,1]$. Thus, $ W $ is allowed to be a Lagrangian with corners.
The idea to consider Lagrangian boundary conditions for Lagrangians comes from~\cite{SolomonYuval20, SolomonYuval21}.

The second problem is that without cylindrical ends there is no canonical composition. Given Lagrangian cobordisms $ W_{01}: L_0 \to L_1 $ and $W_{12}: L_1 \to L_2$, one can not in general concatenate them: their horizontal boundaries in the pushforward of $L_1$ need not match in the cotangent fiber direction. Instead, one can interpolate the Hamiltonian isotopies locally generating $W_{01}$ and $W_{12}$ near $L_1$ to construct a composition cobordism $ W_{02}: L_0 \to L_2$ that agrees with $W_{01} \cup W_{12}$ away from $L_1$. We prove that this composition fits into a \emph{second-order} Lagrangian 2-cobordism $W_{012} \subset X \times T^\ast \Delta^2$ over the 2-simplex, that simultaneously relates $W_{01}$, $W_{12}$ and $W_{02}$.

While the cobordism $W_{02}$ is not unique, we show that it is unique up to a 2-isomorphism, in the following sense: If $W_{02}'$ is another composition, there exists a Lagrangian 2-cobordism that exhibits $W_{02}$ as the composition of $W_{02}'$ with the identity cobordism on $ L_2$. Algebraically, our Lagrangian cobordisms naturally assemble into a \textit{Kan complex}: a simplicial set whose $k$-simplices are Lagrangian $k$-cobordisms over the $k$-simplex. The Kan condition, which ensures that every horn can be filled, corresponds to a composition property that we establish for higher cobordisms. This structure provides a model for an $\infty$-groupoid, which can be viewed as an enhanced space of Lagrangian submanifolds relative to $ \vlg $, and leads us to our main motivation.

The space $\cc{L}(X,L)$ of Lagrangian submanifolds diffeomorphic to $L$ is an infinite dimensional manifold. In~\cite{Solomon13}, the second author introduced a functional $\cc{C}$ on homotopy classes of exact Lagrangian paths in $\cc{L}(X,L)$, generalizing the classical Calabi functional. In the Calabi-Yau setting, $\cc{C}$ interacts naturally with the geometry of positive Lagrangians and provides a framework for studying the existence and uniqueness of special Lagrangians. In an accompanying paper, we show that $\cc{C}$ extends naturally to Lagrangian cobordisms. This extension is invariant under higher cobordisms, and additive under the composition of cobordisms defined here. Furthermore, we demonstrate that the notion of geodesics of positive Lagrangians admits a natural extension to Lagrangian cobordisms.
Thus, the $\infty$-groupoid structure on the space of Lagrangians provides a flexible framework for studying the variational properties of $\cc{C}$.

\subsection{Statements of results}

\begin{definition}
	[{\cite[Definition~2.1]{SolomonYuval20}}]
	Let $X$ and $L$ be smooth manifolds with corners. An \emph{immersed submanifold} $W \subset X$ of type $L$ is an equivalence class of immersions
	\begin{equation}
		W = [\Psi: L \to X],
	\end{equation}
	where the immersions $\Psi$, $\Psi'$ are \emph{equivalent} if there exists a diffeomorphism of manifolds with corners $\vp: L \to L$ such that $\Psi' = \Psi \circ \vp$ .
	We call $\Psi$ \emph{free} if the only diffeomorphism $\vp:L\to L$
	satisfying $\Psi\circ\vp=\Psi$ is the identity. We call $W$ \emph{free}
	if some, and thus any representative is free.

	A \emph{point} $p \in W$ is an equivalence class of pairs $p = [(\Psi,q)]$, where $\Psi: L \to X$ represents $W$ and $q \in L$. The pairs $(\Psi,q), (\Psi', q')$ are \emph{equivalent} if there exists a diffeomorphism $\vp:L \to L$ such that
	\begin{equation}
		\Psi' = \Psi \circ \vp,
		\quad \quad
		\vp(q')=q.
	\end{equation}
	For $p = [(\Psi,q)]$  denote $p_0 = \Psi(q)$, and define the tangent space
	\begin{equation}
		T_{p}W := d\Psi_q(T_qL) \subset T_{p_0}X.
	\end{equation}
\end{definition}

\begin{definition}
	Let $(X,\omega_X)$ be a $2n$-dimensional symplectic manifold, and let $B$ be a $k$-dimensional manifold with corners. Define the symplectic product
	\begin{equation}
		(X_B,\omega_{X_B}) := \bigl(T^\ast B \times X,\; \omega_{T^\ast B} \oplus \omega_X \bigr).
	\end{equation}
	An \emph{immersed Lagrangian} $W \subset X_{B}$ is an immersed submanifold of type $L$ with $\dim L = n+k$ such that, for some and thus for any representative $\Psi: L \to X_B$, one has
	\begin{equation}
		\Psi^*\omega_{X_B} = 0.
	\end{equation}
	We call $W$ \emph{compact} if $L$ is compact.
	We denote by $\Psi_X$, $\Psi_{T^\ast B}$, and $\Psi_B$ the compositions of $\Psi$ with the canonical projections. We refer to immersed submanifolds of $X_B$ as immersed submanifolds \emph{over $B$}.
\end{definition}

Consider the symplectic product
\begin{equation}
	(X_{\overline B},\omega_{X_{\overline B}}) := \bigl(T^\ast B \times X,\; (-\omega_{T^\ast B}) \oplus \omega_X \bigr).
\end{equation}
The canonical symplectomorphism
\begin{equation}
	X_B \to X_{\overline B}, \quad ((b,\xi),x) \mapsto ((b,-\xi),x),
\end{equation}
identifies immersed Lagrangians $W \subset X_B$ with immersed Lagrangians $\overline{W} \subset X_{\overline B},$ which are also called \emph{Lagrangian correspondences} $\overline{W}: T^\ast B \dasharrow X$ (Definition~\ref{def:LagrangianCor}).
We will use this identification to define base-change operations by geometric composition (Definition~\ref{def:geometricComposition}  and Lemma~\ref{lem:geometricComposition}). Recall~\cite{Weinstein} that a smooth map $f: A \to B$ determines a Lagrangian correspondence
\begin{equation}
	T^\ast f: T^\ast B \dasharrow T^\ast A,
\end{equation}
known as the \emph{cotangent lift} (Definition \ref{def:cotangentLift}).
Let $(T^\ast f)^t: T^\ast A \dasharrow T^\ast B$ be the transpose correspondence (Definition~\ref{def:LagrangianCor}).
\begin{definition}
	For any immersed Lagrangian $W \subset X_B$ such that $f$ and $\Psi_B$ are transverse, we define the \emph{pullback of $W$ along $f$} to be the immersed Lagrangian $f^*W \subset X_A$ given by
	\begin{equation}
		\overline{f^\ast W} := \overline{W} \circ (T^\ast f)^t: T^\ast A \dasharrow X.
	\end{equation}
	If $f$ is an immersion, for any immersed Lagrangian $V \subset X_A$ we define the \emph{pushforward of $V$ along $f$} to be the immersed Lagrangian $f_*V \subset X_B$ given by
	\begin{equation}
		\overline{f_\ast V} := \overline{V} \circ T^\ast f : T^\ast B \dasharrow X.
	\end{equation}
\end{definition}

Let $W \subset X_B$ be an immersed submanifold represented by $\Psi: L \to X_B$.
For $d \geq 0$, let $C_d(L)$ denote the codimension-$d$ corner stratum of $L$, which is a manifold with corners with canonical immersion $\iota_d: C_d(L) \to L$. For definitions and notation regarding manifolds with corners, see Appendix \ref{app: ManifoldsWithCorners}.
Suppose that $\Psi_B: L \to B$ is a \emph{b-submersion} (Definition \ref{def:b_submersions}), which means that $d\Psi_B$ maps the boundary-normal directions of $L$ onto the boundary-normal directions of $B$. Then $\Psi_B$ induces a decomposition of $\partial L$ into horizontal (over $\partial B)$ and vertical (over $\mathrm{int}(B)$) boundary components, giving a decomposition (Definition~\ref{def:submersion_mwc_tensors})
\begin{equation}
	C_d(L) = \bigcup_{h+v=d} \tensor{C}{_h^v}(L),
\end{equation}
where $\tensor{C}{_h^v}(L)$ denotes corners meeting $h$ horizontal and $v$ vertical boundary components.
We denote by $\tensor{C}{_h^v}(W) \subset X_B$ the immersed submanifold represented by the composition
\begin{equation}
	\tensor{C}{_h^v}(L) \xrightarrow{\;\tensor{\iota}{_h^v}\;} L \xrightarrow{\;\Psi\;} X_B,
\end{equation}
where $\tensor{\iota}{_h^v}$ is the canonical immersion.

\begin{definition}
	\label{def:LagrangiankCobordism}
	Let $\vlg \subset X$ be a smoothly embedded Lagrangian submanifold. Denote by $\vlg_B \subset X_B$ the pullback of $\vlg$ along the constant map $B \to \mathrm{pt}$, which is the product of $\vlg$ with the zero section in $T^\ast B$.
	An \emph{immersed Lagrangian $k$-cobordism with corners of type $L$ over $B$ relative to $\vlg$} is an immersed Lagrangian $W \subset X_B$ of type $L$ such that:
	\begin{enumerate}
		\item \label{itm:bsubmersion}
		      For some and thus any Lagrangian immersion $\Psi: L \to X_B$ representing $W,$ the projection $\Psi_B: L \to B$ is a b-submersion.
		\item \label{itm:verticalBoundary}
		      Vertical boundary condition: $\tensor{C}{_0^1}(W) \subset X_B$ is an immersed submanifold of $\vlg_B$.
		\item \label{itm:verticalTransversality}
		      Vertical transversality: for all $p \in \tensor{C}{_0^1}(W)$ we have $T_p W \neq T_{p_0}\vlg_B$.
		\item \label{itm:verticalCorners}
		      No vertical corners: $\tensor{C}{_0^v}(W) = \emptyset$ for $v \geq 2$.
	\end{enumerate}
	From now on we refer to $W$ simply as a \emph{Lagrangian cobordism over $B$ relative to $\vlg$}.
\end{definition}

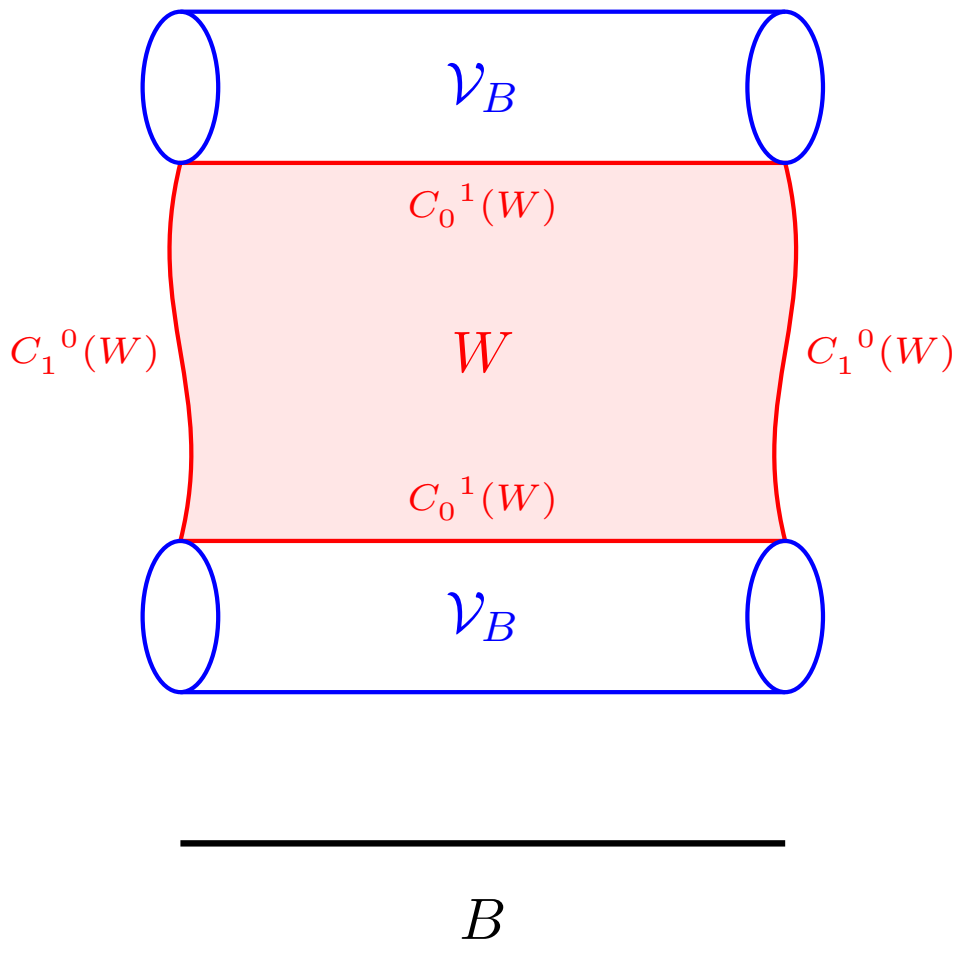
\begin{figure}[htbp]
	\centering
	\begin{tikzpicture}
		\coordinate (A) at (0,0);   
		\coordinate (A') at (0,-1);
		\coordinate (A'') at (0,-2);
		\coordinate (B) at (4,0);   
		\coordinate (B') at (4,-1);
		\coordinate (B'') at (4,-2);
		\coordinate (C) at (4,2.5); 
		\coordinate (C') at (4,3.5);
		\coordinate (D) at (0,2.5); 
		\coordinate (D') at (0,3.5);

		\filldraw[fill=red!10, draw=red, thick]
		(A) -- (B)
		.. controls (3.75,1) and (4.25,1.5) .. (C)
		-- (D)
		.. controls (-0.25,1.5) and (0.25,1) .. (A)
		-- cycle;

		\draw[thick, blue] (A') -- (B');
		\draw[thick, blue] (C') -- (D');

		\draw[very thick] (A'') -- (B'');

		\draw[thick, blue] ($(A)!0.5!(A')$) ellipse [x radius=0.25, y radius=0.5];
		\draw[thick, blue] ($(B)!0.5!(B')$) ellipse [x radius=0.25, y radius=0.5];
		\draw[thick, blue] ($(C)!0.5!(C')$) ellipse [x radius=0.25, y radius=0.5];
		\draw[thick, blue] ($(D)!0.5!(D')$) ellipse [x radius=0.25, y radius=0.5];

		\node[left, font=\tiny, red] at ($(A)!0.5!(D)$) {$\tensor{C}{_1^0}(W)$};
		\node[right, font=\tiny, red] at ($(B)!0.5!(C)$) {$\tensor{C}{_1^0}(W)$};
		\node[above, font=\tiny, red] at ($(A)!0.5!(B)$) {$\tensor{C}{_0^1}(W)$};
		\node[below, font=\tiny, red] at ($(C)!0.5!(D)$) {$\tensor{C}{_0^1}(W)$};

		\node[blue] at (2,3) {$\vlg_B$};
		\node[red] at (2,1.25) {$W$};
		\node[blue] at (2,-0.5) {$\vlg_B$};
		\node at (2,-2.5) {$B$};

	\end{tikzpicture}
	\caption{A Lagrangian 1-cobordism over $B=[0,1]$ relative to $\vlg$.}
\end{figure}

Taking $B = \mathrm{pt}$ in Definition \ref{def:LagrangiankCobordism}, a Lagrangian $0$-cobordism relative to $\vlg$ is an immersed Lagrangian with boundary $W \subset X$, whose boundary is an immersed submanifold of $\vlg$.
The pullback of a Lagrangian cobordism $W \subset X_B$ along a map $f:A \to B$ transverse to $\Psi_B$ is a Lagrangian cobordism over $A$. In particular, for the inclusion of a component $A \subset C_h(B)$, transversality to the b-submersion $\Psi_B$ is automatic, so the \emph{restriction} $ W|_A := f^\ast W $  is a Lagrangian $(k-h)$-cobordism.

\begin{definition}
	Let $W$ be a Lagrangian $k$-cobordism over $B$ relative to $\vlg.$ Let $B_1, \dots,B_l$ be the boundary components of $B$, and let $\iota_i: B_i \to B$ be the canonical inclusion.
	Let $W_i$ be the Lagrangian $(k-1)$-cobordism over $B_i$ given by restriction,
	\[
		W_i = W|_{B_i} = \iota_i^\ast W.
	\]
	Then, we say that $W$ is a cobordism \emph{between} $W_1,\ldots,W_l.$
\end{definition}

The pushforward of a Lagrangian cobordism need not be a Lagrangian cobordism since the projection to the base need not be a b-submersion. Nevertheless, pushforwards appear  as Lagrangian boundary conditions.
Let $W$ be a Lagrangian $k$-cobordism over $B$ between Lagrangian $(k-1)$-cobordisms $W_1,\dots,W_l$ over boundary components $B_1,\dots,B_l$. Then the boundary of $W$ is canonically immersed in the pushforwards $(\iota_i)_\ast W_i$ as shown in Section~\ref{sec:bndry}.

One natural way to construct a Lagrangian cobordism over $B$ relative to $\vlg$ is by suspension of a Hamiltonian isotopy. Suppose $W_0 \subset X$ is an immersed Lagrangian of type $L$ with boundary in a smoothly embedded Lagrangian $\vlg$. Let $\phi: B \to \mathrm{Ham}(X,\omega_X)$ be a smooth family of Hamiltonian symplectomorphisms such that $\phi_b(\vlg) = \vlg$ for all $b \in B$. Assume that the family $b\mapsto\phi_b(W_0)$ is exact relative
to $\vlg$ in the sense of Definition~\ref{def:exactFamily}.
The \emph{suspension} of the family $(\phi_b(W_0))_{b \in B}$ is then a Lagrangian cobordism over $B$ relative to $\vlg$ of type $L \times B$, with $\Psi_B: L \times B \to B$ the projection.
In general, however, we only require $\Psi_B$ to be a b-submersion, not a submersion. Consequently, we obtain cobordisms between non-isotopic Lagrangians in $X$ that may have different topological types.

\begin{definition}
	Let $\mathbf{\Delta}$ be the category of linearly ordered sets $[n] = \{0,\dots,n\}$ and order-preserving maps. A \emph{simplicial set} is a functor $S: \mathbf{\Delta}^\mathrm{op} \to \text{Set}$. We write $S_n := S([n])$; elements of $S_n$ are called \emph{$n$-simplices}.
	For $0\leq i \leq n$, let
	\begin{equation}
		\delta_{n,i}: [n-1] \to [n]
	\end{equation}
	be the injective map that skips $i$, and
	\begin{equation}
		\sigma_{n,i}: [n+1] \to [n]
	\end{equation}
	the surjective map that identifies $i$ and $i+1$. The maps $\delta_{n,i}$ and $\sigma_{n,i}$ generate all morphisms of the category $\mathbf{\Delta}$. The associated \emph{face} and \emph{degeneracy} maps are
	\begin{equation}
		d_{n,i} = S(\delta_{n,i}): S_{n} \to S_{n-1},
		\quad \quad
		s_{n,i} = S(\sigma_{n,i}): S_{n} \to S_{n+1}.
	\end{equation}
	A \emph{map of simplicial sets} is a natural transformation between such functors.
\end{definition}

\begin{definition}
	Let $\Delta^n \subset \bb{R}^{n+1}$ be the \emph{geometric $n$-simplex} given by the convex hull of the standard basis vectors (Definition \ref{def:n-simplex}).
	Let $\Delta[n]$ denote the \emph{standard simplicial $n$-simplex} $\hom_{\mathbf\Delta}(-,[n])$. A $k$-simplex of $\Delta[n]$ is an order-preserving map $\alpha: [k] \to [n]$. For a simplicial set $S$, the Yoneda Lemma gives a bijection between $S_n$ and maps $\Delta[n] \to S$ of simplicial sets.
\end{definition}

\begin{definition}
	\label{def:n-i-horn}
	A simplicial set $T$ is a \emph{subset} of a simplicial set $S$ if $T_n \subset S_n$ and $S(\alpha)|_{T_n} = T(\alpha)$ for all $n$ and $\alpha: [k] \to [n]$. For $0\leq i \leq n$, the \emph{$(n,i)$-horn} $\Lambda_i[n] \subset \Delta[n]$ is the subset of $\Delta[n]$ obtained by removing the $n$-simplex and the facet opposite the $i$-th vertex. Equivalently, the $k$-simplices are given by
	\begin{equation}
		\Lambda_i[n]_k =
		\left\{
		\alpha: [k] \to [n] \text{ order preserving} \;|\;
		\mathrm{im}(\alpha) \subset [n]\setminus \{j\}
		\text{ for some }
		j \neq i
		\right\}.
	\end{equation}
\end{definition}

The following lemma is proved in Section~\ref{subsec:kan-complex}.
\begin{lemma}
	\label{lem:simplicial-set-cob}
	Let $(X,\omega_X)$ be a symplectic manifold and $\vlg \subset X$ an embedded Lagrangian. The following data determines a simplicial set  $\cob(X,\vlg): \mathbf\Delta^\mathrm{op} \to \mathrm{Set}$:
	\begin{enumerate}
		\item For $n\geq 0$, let $\cob(X,\vlg)_n$ be the set of compact Lagrangian $n$-cobordisms $W \subset X_{\Delta^n}$ over the geometric $n$-simplex, relative to $\vlg$, such that $W|_F$ is free for every nonempty face $F\subseteq\Delta^n$, including $\Delta^n$ itself.
		\item For $ 0 \leq i \leq n $, the $i$-th face map $d_{n,i}: \cob(X,\vlg)_n \to \cob(X,\vlg)_{n-1}$ is pullback along the face inclusion
		      \begin{equation}
			      \bar\delta_{n,i}
			      : \Delta^{n-1} \to \Delta^{n},
			      \quad \quad
			      (s_0,\dots,s_{n-1})
			      \mapsto (s_0,\dots, s_{i-1}, 0, s_i, \dots, s_{n-1}),
		      \end{equation}
		      i.e. $d_{n,i}(W)= (\bar\delta_{n,i})^\ast W$.
		\item For $ 0\leq i \leq n $, the $i$-th degeneracy map $s_{n,i}: \cob(X,\vlg)_n \to \cob(X,\vlg)_{n+1}$ is given by pullback along the map
		      \begin{equation}
			      \bar\sigma_{{n},i}
			      : \Delta^{n+1} \to \Delta^{n},
			      \quad
			      (s_0,\dots,s_{n+1})
			      \mapsto (s_0,\dots,s_{i-1},s_i+s_{i+1},s_{i+2},\dots,s_{n+1}).
		      \end{equation}
		      i.e. $s_{n,i}(W)= (\bar\sigma_{n,i})^\ast W$.
	\end{enumerate}
\end{lemma}

\begin{definition}
	\label{def:Lag-over-sSet}
	Let $S$ be a simplicial set.
	A \emph{Lagrangian $W$ over $S$} is a map of simplicial sets
	$W: S \to \cob(X,\vlg)$.
\end{definition}

The following theorem is the main result of this work. It shows that $\cob(X,\vlg)$ is a \emph{Kan complex}: a simplicial set in which every horn can be filled. The proof is given in Section~\ref{subsec:kan-complex}.
\begin{theorem}
	\label{thm:hornFilling}
	Fix $n\geq 1$ and $0 \leq i \leq n$. Let $W: {\Lambda_i[n]} \to \cob(X,\vlg)$ be a Lagrangian over the $(n,i)$-horn.
	Then there exists a Lagrangian $\widetilde{W}: \Delta[n] \to \cob(X,\vlg)$ extending $W,$ that is, the diagram
	\begin{equation}
		\label{eq:hornFilling}
		\begin{tikzcd}
			\Lambda_i[n]
			\ar[d, hook']
			\ar[dr,"W"]
			\\
			\Delta[n]
			\ar[r, "\widetilde{W}"', dashed]
			&
			\cob(X,\vlg)
		\end{tikzcd}
	\end{equation}
	commutes. In particular, $\cob(X,\vlg)$ is a Kan complex.
\end{theorem}

Recall that a \emph{groupoid} is a category in which every morphism is invertible. Likewise, an \emph{$\infty$-groupoid} is an $\infty$-category in which all $n$-morphisms are invertible for every $n\geq 1$. A theorem of Joyal~\cite{Joyal} shows that Kan complexes model $\infty$-groupoids.
\begin{corollary}
	\label{cor:groupoid}
	There exists an $\infty$-groupoid $\cob(X,\vlg)$ whose objects are compact freely immersed Lagrangians in $X$ with boundary in $\vlg$, and whose $n$-morphisms are Lagrangian $n$-cobordisms over $\Delta^n$ relative to $\vlg$. In particular, the homotopy 1-category $h\cob(X,\vlg)$ is a groupoid with the same objects, and morphisms given by homotopy classes (equivalently, 2-cobordism classes) of Lagrangian 1-cobordisms.
\end{corollary}

\subsection{Outline}
In Section~\ref{sec:LagrangianCobordisms}  we define the pullback and pushforward
of immersed Lagrangians over a base, describe the boundary conditions
for Lagrangian cobordisms in terms f pushforwards, and construct suspensions of exact families. In~\ref{sec:verticalSystems} we outline a generalization of Definition \ref{def:LagrangiankCobordism} that replaces $\cc{V}_B$ by a non-trivial Lagrangian cobordism and permits higher vertical corners.

In Section~\ref{sec:collars} we construct compatible exact collar models
near the corners of a Lagrangian $n$-cobordism over $\Delta^n$
(Theorem~\ref{thm:compatible_exact_models}). In each collar the cobordism
is represented by the graph of an exact one-form in a Weinstein chart.
The proof relies on a compatible construction of the collar maps, the Weinstein charts,
and the generating functions in
Sections~\ref{sec:extend-collars}--\ref{sec:gen-fn}.

In Section~\ref{sec:composition} we prove the horn-filling theorem.
In Sections~\ref{sec:extend-gen-empty} and~\ref{sec:ext_Lag_horns}
we extend the generating functions to the facet collars of $\Lambda_i^n$
and glue the resulting cobordisms over a neighborhood of the horn.
In Section~\ref{subsec:kan-complex} we use this construction to fill the horn, proving
Theorem~\ref{thm:hornFilling}. Finally, we describe composition in the homotopy groupoid $h\cob(X,\vlg)$.

Appendix~\ref{app: ManifoldsWithCorners} reviews manifolds with corners
and fixes notation for horizontal and vertical corners.

\subsection{AI declaration} ChatGPT was used to proofread the manuscript and it helped fix minor errors.

\subsection{Acknowledgements}
The authors would like to thank K.~Daniel and E.~Kosloff for helpful conversations. The authors were partially supported by BSF grants 2020329 and 2024293. J.~S. was partially supported by the Miriam and Julius Vinik Chair in Mathematics.

\section{Immersed Lagrangian cobordisms with corners}
\label{sec:LagrangianCobordisms}

\subsection{Lagrangian correspondences}
\label{sec:LagCorr}
For a symplectic manifold $(M,\omega)$, we denote by $\overline M$ the symplectic manifold $(M,-\omega)$.

\begin{definition}
	\label{def:LagrangianCor}
	Let $(M_0,\omega_0)$, $(M_1,\omega_1)$ be symplectic manifolds with corners, and let $L$ be a manifold with corners.
	A \emph{Lagrangian correspondence} $W: M_0 \dasharrow M_1$ of type $L$ is an immersed Lagrangian $W \subset \overline M_0 \times M_1$ of type $L$.
	The \emph{transpose} correspondence $W^t: M_1 \dasharrow M_0$ is the image of $W$ under the swap isomorphism $\overline M_0 \times M_1 \simeq \overline M_1 \times M_0$.
\end{definition}

\begin{definition}
	\label{def:geometricComposition}
	Let $(M_i,\omega_i)$, $i=0,1,2$, be symplectic manifolds.
	Let
	\begin{equation}
		\Psi_{01}: L_{01} \to \overline M_0 \times M_1,
		\qquad
		\Psi_{12}: L_{12} \to \overline M_1 \times M_2,
	\end{equation}
	represent Lagrangian correspondences $W_{01}: M_0 \dasharrow M_1$ and $W_{12}: M_1 \dasharrow M_2$.
	For a map $\Psi_{ij}$ with target $\overline M_i \times M_j$, we write
	\begin{equation}
		\Psi_{(i)j}: L_{ij} \to M_j
		\qquad\text{and}\qquad
		\Psi_{i(j)}: L_{ij} \to M_i
	\end{equation}
	for its components. Assume the maps
	\begin{equation}
		\Psi_{(0)1}: L_{01} \to M_1,
		\qquad
		\Psi_{1(2)}: L_{12} \to M_1
	\end{equation}
	are transverse in the sense of Definition \ref{def:transversality}, so that $L_{01} \times_{M_1} L_{12}$ is a manifold with corners by Theorem \ref{thm:fiber_products_corners}.
	Define
	\begin{equation}
		\label{eq:geometric_composition}
		\Psi_{02}: L_{01} \times_{M_1} L_{12} \to \overline M_0 \times M_2,
		\qquad
		\Psi_{02}(p,q)=\bigl(\Psi_{0(1)}(p),\, \Psi_{(1)2}(q)\bigr).
	\end{equation}
\end{definition}

\begin{lemma}
	\label{lem:geometricComposition}
	The map \eqref{eq:geometric_composition} is a Lagrangian immersion.
\end{lemma}

\begin{proof}
	It follows from the diagram
	\begin{equation}
		\begin{tikzcd}
			&
			L_{01} \times_{M_1} L_{12}
			\ar[r, "\mathrm{pr}_1"]
			\ar[d, "\mathrm{pr}_2"]
			&
			L_{01}
			\ar[r, "\Psi_{0(1)}"]
			\ar[d, "\Psi_{(0)1}"]
			&
			M_0
			\\
			M_2
			&
			L_{12}
			\ar[l, "\Psi_{(1)2}"']
			\ar[r, "\Psi_{1(2)}"]
			&
			M_1
		\end{tikzcd}
	\end{equation}
	that $\mathrm{pr}_1^\ast \Psi_{(0)1}^\ast \omega_1 = \mathrm{pr}_2^\ast \Psi_{1(2)}^\ast \omega_1 $. Since $\Psi_{01}$, $\Psi_{12}$ are Lagrangian, we have $\Psi_{(0)1}^\ast \omega_1 = \Psi_{0(1)}^\ast \omega_0$ and $\Psi_{(1)2}^\ast \omega_2 = \Psi_{1(2)}^\ast \omega_1$. Putting this together, we get
	\begin{equation}
		\Psi_{02}^\ast (-\omega_0 \oplus \omega_2)
		= -\mathrm{pr}_1^\ast\Psi_{0(1)}^\ast \omega_0 + \mathrm{pr}_2^\ast\Psi_{(1)2}^\ast \omega_2
		= -\mathrm{pr}_1^\ast\Psi_{(0)1}^\ast \omega_1 + \mathrm{pr}_2^\ast\Psi_{1(2)}^\ast \omega_1 =0,
	\end{equation}
	so $\Psi_{02}$ is Lagrangian. We now prove that $\Psi_{02}$ is an immersion.
	Let $ 0 \neq v \in T(L_{01} \times_{M_1} L_{12}) $ and assume without loss of generality that $d\mathrm{pr}_1(v) \neq 0$. Suppose that
	\begin{equation}
		\label{eq:vanishingProjection01}
		d(\Psi_{0(1)} \circ \mathrm{pr}_1)(v)=0
	\end{equation}
	and assume by contradiction that
	\begin{equation}
		\label{eq:vanishingProjection02}
		d(\Psi_{(1)2} \circ \mathrm{pr}_2)(v) = 0.
	\end{equation}
	Since $\Psi_{01}$ is immersion, by \eqref{eq:vanishingProjection01} we must have $d(\Psi_{1(2)}\circ \mathrm{pr}_2)(v)=d(\Psi_{(0)1}\circ \mathrm{pr}_1)(v)\neq 0$. Since $\omega_1$ is non-degenerate and the maps $\Psi_{(0)1}$, $\Psi_{1(2)}$ are transverse, there exists $u \in TL_{01}$ and $w \in TL_{12}$ such that
	\begin{equation}
		\label{eq:nondegenerate01}
		\omega_1(d(\Psi_{1(2)}\circ \mathrm{pr}_2)(v), d\Psi_{(0)1}(u) + d\Psi_{1(2)}(w)) \neq 0.
	\end{equation}
	Since $\Psi_{01}$ is Lagrangian, we have
	\begin{align}
    \MoveEqLeft[1] \omega_1 (d(\Psi_{1(2)}\circ \mathrm{pr}_2)(v), d\Psi_{(0)1}(u) ) \notag\\
		 & =
		\omega_1 (d(\Psi_{(0)1}\circ \mathrm{pr}_1)(v), d\Psi_{(0)1}(u) )
		\\
		 & =
		\Psi_{(0)1}^\ast \omega_1(d\mathrm{pr}_1(v), u)
		\\
		 & = \Psi_{(0)1}^\ast \omega_1(d\mathrm{pr}_1(v), u) - \Psi_{0(1)}^\ast \omega_0(d\mathrm{pr}_1(v),u)
		 &                                                                                                    &
		\text{by } \eqref{eq:vanishingProjection01}
		\\
		 & = \Psi_{01}^\ast (-\omega_0 \oplus \omega_1)(d\mathrm{pr}_1(v),u) =0.
		\label{eq:nondegenerate02}
	\end{align}
	Since $\Psi_{12}$ is Lagrangian, we have
	\begin{align}
        \MoveEqLeft[1] \omega_1(d(\Psi_{1(2)}\circ \mathrm{pr}_2)(v), d\Psi_{1(2)}(w)) \notag\\
		 & =
		\Psi_{1(2)}^\ast \omega_1(d\mathrm{pr}_2(v),w)
		\\
		 & =
		\Psi_{1(2)}^\ast \omega_1(d\mathrm{pr}_2(v),w) - \Psi_{(1)2}^\ast \omega_2(d\mathrm{pr}_2(v),w)
		 &   &
		\text{by } \eqref{eq:vanishingProjection02}
		\\
		 & =
		-\Psi_{12}^\ast (-\omega_1 \oplus \omega_2)(d\mathrm{pr}_2(v),w) =0.
		\label{eq:nondegenerate03}
	\end{align}
	From \eqref{eq:nondegenerate02} and \eqref{eq:nondegenerate03} we get a contradiction to \eqref{eq:nondegenerate01}. Thus \eqref{eq:vanishingProjection01} and \eqref{eq:vanishingProjection02} cannot hold simultaneously and $\Psi_{02}$ is an immersion.
\end{proof}

\begin{definition}
	Continue with the assumptions of Definition \ref{def:geometricComposition}. By Lemma \ref{lem:geometricComposition}
	\begin{equation}
		W_{12}\circ W_{01} := [\Psi_{02}: L_{01} \times_{M_1} L_{12} \to \overline{M}_0 \times M_2]
	\end{equation}
	is a Lagrangian correspondence $W_{12}\circ W_{01}: M_0 \dasharrow M_2$, the \emph{geometric composition} of $W_{01}$ and $W_{12}$. By the universal property of fiber products, this immersed Lagrangian is well defined independently of the representatives $\Psi_{01}$ and $\Psi_{12}$.
\end{definition}

\subsection{The cotangent lift}
\label{sec:cotLift}

\begin{definition}
	\label{def:cotangentLift}
	Let $f:A \to B$ be a smooth map of manifolds with corners. Define
	\begin{equation}
		\label{eq:cotangentLift}
		\Psi_f:\ f^\ast T^\ast B \longrightarrow \overline{T^\ast B}\times T^\ast A,
		\qquad
		\xi \longmapsto \bigl(\,\tilde f(\xi),\, df^\ast(\xi)\,\bigr),
	\end{equation}
	where $\tilde f: f^\ast T^\ast B \to T^\ast B$ is the natural projection and
	$df^\ast: f^\ast T^\ast B \to T^\ast A$ is the canonical pullback map.
\end{definition}

\begin{lemma}
	The map \eqref{eq:cotangentLift} is a Lagrangian immersion, representing a Lagrangian correspondence
	\begin{equation}
		T^\ast f:\ T^\ast B \dasharrow T^\ast A,
	\end{equation}
	called the \emph{cotangent lift} of $f$.
\end{lemma}

\begin{proof}
	Let $\theta_A$, $\theta_B$ be the canonical 1-forms on $T^\ast A$, $T^\ast B$. We first claim that
	\begin{equation}
		\label{eq:canonical1Forms}
		(df^\ast)^\ast \theta_A = (\tilde f)^\ast \theta_B.
	\end{equation}
	Taking exterior derivative of both sides gives $\Psi_f^*(-\omega_B \oplus \omega_A)=0$ so $\Psi_f$ is Lagrangian. Indeed, consider the commutative diagram
	\begin{equation}
		\begin{tikzcd}
			T^\ast A \arrow[rd, "\pi_A"'] & f^\ast T^\ast B \arrow[l, "df^\ast"'] \arrow[r, "\tilde f"] \arrow[d, "\text{pr}_A"] & T^\ast B \arrow[d, "\pi_B"] \\
			& A \arrow[r, "f"]                                                     & B
		\end{tikzcd}
	\end{equation}
	For a manifold $M$, denote by $\left<,\right>_M$ the natural pairing of $T^\ast M$ and $TM$.  Let $(a,\xi) \in f^\ast T^\ast B$ and $v \in T_{(a,\xi)}(f^\ast T^\ast B)$. Evaluating both sides of \eqref{eq:canonical1Forms} on $v$ and using the canonical property of $\theta_A$, $\theta_B$ we get
	\begin{align}
		\left<\theta_A, d(df^\ast)_{(a,\xi)}(v) \right>_{T^\ast A}
		 & =
		\left<df^\ast (\xi), d\pi_A(d(df^\ast)_{(a,\xi)}(v)) \right>_{A}
		\\
		 & =
		\left<\xi, d (f \circ \pi_A \circ df^\ast)_{(a,\xi)}(v) \right>_B
		\\
		 & =
		\left<\xi, d ( \pi_B \circ \tilde f)_{(a,\xi)}(v) \right>_B
		\\
		 & =
		\left<\theta_B, d\tilde f_{(a,\xi)}(v) \right>_{T^\ast B}
	\end{align}
	which proves \eqref{eq:canonical1Forms}. Note that if $v \neq 0$ and $d (df^\ast)_{(a,\xi)}(v) =0$ then 
    \[
        d(\text{pr}_A)_{(a,\xi)}(v) =0,
    \] 
    so $d\tilde f_{(a,\xi)}(v) \neq 0$, which implies that $\Psi_f$ is an immersion.
\end{proof}

\begin{remark}
	In fact, $\Psi_f$ is an embedding. One can equivalently define $T^\ast f$ as the image, under the canonical isomorphism $T^\ast (A\times B) \simeq \overline{T^\ast B} \times T^\ast A$, of the conormal bundle $N^\ast \Gamma_f$ of the graph $\Gamma_f \subset A \times B$.
\end{remark}

\begin{lemma}
	The cotangent lift is functorial: Let $f:A \to B$ and $g: B \to C$. Then
	\begin{equation}
		T^\ast (g\circ f) = T^\ast f \circ T^\ast g.
	\end{equation}
\end{lemma}

\begin{proof}
	Let $\Psi_f: f^\ast T^\ast B \to \overline{T^\ast B} \times T^\ast A$ and $\Psi_g: g^\ast T^\ast C \to \overline{T^\ast C} \times T^\ast B$. Let $\vp$ be the canonical isomorphism
    \begin{equation}
      \begin{aligned}
        (g\circ f)^\ast T^\ast C
          &\xrightarrow{\sim} f^\ast g^\ast T^\ast C
           \xrightarrow{\sim} g^\ast T^\ast C \times_B A\\
          &\xrightarrow{\sim}
           g^\ast T^\ast C \times_{T^\ast B}(T^\ast B\times_B A)
           \xrightarrow{\sim}
           g^\ast T^\ast C \times_{T^\ast B}f^\ast T^\ast B.
      \end{aligned}
    \end{equation}
	It is straightforward to verify that the diagram
	\begin{equation}
		\begin{tikzcd}
			(g\circ f)^\ast T^\ast C \arrow[rrd, "\Psi_{g\circ f}"] \arrow[d, "\vp"]               &  &                                    \\
			g^\ast T^\ast C \times_{T^\ast B} f^\ast T^\ast B \arrow[rr, "\Psi_g \circ \Psi_f"'] &  & \overline{T^\ast C}\times T^\ast A
		\end{tikzcd}
	\end{equation}
	commutes, so $[\Psi_{g\circ f}] = [\Psi_g] \circ [\Psi_f]$.
\end{proof}

\begin{corollary}
	\label{cor:functorialPullback}
	Pullback and pushforward of immersed Lagrangians are functorial. Let  $f: A \to B$ and $g: B \to C$ be smooth maps of manifolds with corners.
	\begin{enumerate}
		\item Let $W\subset X_C$ be an immersed Lagrangian. Then $(g\circ f)^\ast W = f^\ast g^\ast W$ whenever these pullbacks are defined.
		\item Let $W \subset X_A$ be an immersed Lagrangian and suppose that $f$ and $g$ are immersions. Then $(g\circ f)_\ast W = g_\ast f_\ast W$.
	\end{enumerate}
\end{corollary}

\subsection{Pullback and pushforward}
\label{sec:PullPush}
In this subsection, we will construct the pullback and pushforward of Lagrangian cobordisms, discussed in the introduction.
\begin{prop}
	\label{prop:pushforward}
	Let $ V \subset X_A $ be an immersed Lagrangian of type $ L $ represented by $ \Psi: L \to X_A $. Let $ f:A \to B $ be a smooth immersion of manifolds with corners. Then $ f_\ast V \subset X_B $ is an immersed Lagrangian of type
	\begin{equation}
		f_\ast L:= L\times_{T^\ast A} f^\ast T^\ast B
	\end{equation}
	over $ B $, represented by
	\begin{equation}
		f_\ast \Psi: f_\ast L \to X_B,
		\quad \quad
		(q,\xi) \mapsto (\tilde f(\xi), \Psi_X(q)),
	\end{equation}
	as illustrated in the commutative diagram
	\begin{equation}
		\begin{tikzcd}
			&
			f_\ast L
			\arrow[d, "\mathrm{pr}_2"']
			\arrow[r, "\mathrm{pr}_1"]
			&
			L
			\arrow[d, "\Psi_{T^\ast A}"]
			\arrow[r, "\Psi_X  "]
			&
			X
			\\
			T^\ast B
			&
			f^\ast T^\ast B
			\arrow[r, "df^\ast"]
			\arrow[l, "\tilde f"']
			&
			T^\ast A
			&
		\end{tikzcd}
	\end{equation}
\end{prop}

\begin{proof}
	Since $f$ is an immersion, $df^\ast$ is a submersion and $df^\ast \pitchfork \Psi_{T^\ast A}$. The Proposition follows from Lemma \ref{lem:geometricComposition} by taking $\Psi_{01} = \Psi_f$ and $\Psi_{12} = \Psi$.
\end{proof}

\begin{prop}
	\label{prop:pullback}
	Let $ W \subset X_B $ be a Lagrangian cobordism of type $ L $ relative to $ \vlg $, represented by $ \Psi: L \to X_B $. Let $ f:A \to B $ be a smooth map of manifolds with corners which is transverse to $ W $,  i.e.  $f\pitchfork \Psi_B$ in the sense of Definition~\ref{def:transversality}. Then $ f^\ast W \subset X_A $ is a Lagrangian cobordism of type
	\begin{equation}
		f^\ast L:= L\times_{T^\ast B} f^\ast T^\ast B
	\end{equation}
	over $ A $ relative to $ \vlg $, represented by
	\begin{equation}
		\label{eq:pullback}
		f^\ast \Psi: f^\ast L \to X_A,
		\quad \quad
		(q,\xi) \mapsto (df^\ast(\xi), \Psi_X (q)),
	\end{equation}
	as illustrated in the commutative diagram
	\begin{equation}
		\begin{tikzcd}
			&
			f^\ast L
			\arrow[d, "\mathrm{pr}_2"']
			\arrow[r, "\mathrm{pr}_1"]
			&
			L
			\arrow[d, "\Psi_{T^\ast B}"]
			\arrow[r, "\Psi_X  "]
			&
			X
			\\
			T^\ast A
			&
			f^\ast T^\ast B
			\arrow[r, "\tilde f "]
			\arrow[l, "df^\ast"']
			&
			T^\ast B
			&
		\end{tikzcd}
	\end{equation}
\end{prop}

\begin{proof}
	We first claim that \eqref{eq:pullback} is a Lagrangian immersion. By taking $\Psi_{01} = \Psi_f^t$ and $\Psi_{12} = \Psi$ in Lemma \ref{lem:geometricComposition}, it suffices to show that
	\begin{equation}
		f \pitchfork \Psi_B \iff \tilde f \pitchfork \Psi_{T^\ast B}.
	\end{equation}
	Indeed, since $\text{pr}_A$ is a surjective submersion
	\begin{equation}
		f \pitchfork \Psi_B \iff f \circ \text{pr}_A \pitchfork \Psi_B \iff \pi_B \circ \tilde f \pitchfork \pi_B \circ \Psi_{T^\ast B}.
	\end{equation}
	Since $\tilde f$ is fiberwise surjective, the base transversality $\pi_B \circ \tilde f \pitchfork \pi_B \circ \Psi_{T^\ast B}  $ is equivalent to $\tilde f \pitchfork \Psi_{T^\ast B} $.

	We now verify conditions \eqref{itm:bsubmersion}-\eqref{itm:verticalCorners} in Definition \ref{def:LagrangiankCobordism}. Consider the fiber product diagram
	\begin{equation}
		\begin{tikzcd}
			f^\ast L
			\ar[r, "\mathrm{pr}_2"]
			\ar[d, "\mathrm{pr}_1"']
			&
			f^* T^* B
			\ar[r, "\mathrm{pr}_A"]
			\ar[d, "\tilde f"]
			&
			A
			\ar[d, "f"]
			\\
			L
			\ar[r, "\Psi_{T^\ast B}"]
			&
			T^*B
			\ar[r,"\pi_B"]
			&
			B
		\end{tikzcd}
	\end{equation}
	The lower horizontal arrows compose to $ \Psi_B $, which is a b-submersion since $ \Psi $ satisfies Definition \ref{def:LagrangiankCobordism}\eqref{itm:bsubmersion}. Thus the upper horizontal arrows compose to a b-submersion  $ (f^\ast\Psi)_A $, so $ f^\ast \Psi $ satisfies Definition \ref{def:LagrangiankCobordism}\eqref{itm:bsubmersion}. By Lemma \ref{lem:verticalPullback}, for $ v \geq 0 $  we have a canonical identification
	\begin{equation}
		\tns{C}{_0^v}(f^\ast L)
		\simeq
		f^\ast \tns{C}{_0^v}(L)
		\simeq
		\tns{C}{_0^v}(L) \times_{T^*B} f^\ast T^\ast B,
	\end{equation}
	which immediately implies Definition \ref{def:LagrangiankCobordism}\eqref{itm:verticalCorners} for $ f^*\Psi $.
	Fix
	\begin{equation}
		(q,\xi) \in
		\tns{C}{_0^1}(L) \times_{T^*B} f^\ast T^\ast B
		\simeq
		\tns{C}{_0^1}(f^*L).
	\end{equation}
	Since $ \Psi $ satisfies Definition \ref{def:LagrangiankCobordism}\eqref{itm:verticalBoundary},  $ \Psi_X(q) \in \vlg $ and $ \tilde{f}(\xi)=\Psi_{T^*B}(q)$ lies on the zero section of $T^*B$. In particular, $df^*(\xi)=0$, so
	\begin{equation}
		f^\ast \Psi (q,\xi) = (\Psi_X(q), df^*(\xi)) \in \vlg_A
	\end{equation}
	and $ f^*\Psi $ satisfies Definition \ref{def:LagrangiankCobordism}\eqref{itm:verticalBoundary}.

	We show that $f^*\Psi$ satisfies Definition \ref{def:LagrangiankCobordism}\eqref{itm:verticalTransversality}.
	Set $a := (f^\ast\Psi)_A(q,\xi)$, $b:=\Psi_B(q)$. We have
	\begin{equation}
		T_{(q,\xi)}(f^\ast L)\;\cong\; \Bigl\{(v,e)\in  T_qL \oplus T_\xi(f^*T^*B) \ \Big|\ d(\Psi_{T^*B})_q(v) =d\tilde{f}_\xi(e)\Bigr\}.
	\end{equation}
	Since $\xi$ is in the zero section, we have canonical splittings
	\begin{equation}
		T_{\xi}(f^*T^*B) \simeq T_aA \oplus T_b^*B,
		\qquad
		T_{\tilde{f}(\xi)}(T^*B) \simeq T_bB \oplus T_b^*B,
	\end{equation}
	giving a canonical identification
    \begin{equation}
        \label{eq:tangent-pullback}
        T_{(q,\xi)}(f^\ast L)
        \cong \bigl\{(v,(u,\eta))\bigm|
        d(\Psi_{T^*B})_q(v)=(df_a(u),\eta)\bigr\},
    \end{equation}
	where $v\in T_qL$, $u\in T_aA$, and $\eta\in T_b^*B$.
	We first claim that
	\begin{equation}
		\label{eq:almost-vertical-transversality}
		\exists (v,(u,\eta)) \in T_{(q,\xi)}(f^\ast L),
		\qquad
		d\Psi_q(v) \notin T_{\Psi(q)}\vlg_B.
	\end{equation}
	Since $\Psi$ satisfies Definition~\ref{def:LagrangiankCobordism}\eqref{itm:verticalTransversality} there exists $w\in T_qL$ such that $d\Psi_{q }(w)\notin T_p\vlg_B$. Let $s:=d(\Psi_B)_q(w)\in T_{\Psi_B(q)}B$.
	Recall Definition \ref{def:transversality}, and let $j\geq 0$ such that $a\in S^j(A)$. Since $\Psi_B,f$ are transverse, there exist $u \in T_aS^j(A)$, $w' \in T_q\tns{C}{_0^1}(L)$ with
	\begin{equation}
		d(\Psi_B)_q(w')+df_a(u)=s.
	\end{equation}
	Set $v:=w-w'$. Then
	\begin{equation}
		d(\Psi_B)_q(v)=d(\Psi_B)_q(w)-d(\Psi_B)_q(w')= s - (s - df_a(u))=df_a(u).
	\end{equation}
	Taking $\eta \in T_b^*B$ to be the vertical part of $d(\Psi_{T^*B})_q(v)$, we get that $(v,(u,\eta)) \in T_{(q,\xi)}(f^*L)$ by \eqref{eq:tangent-pullback}.
	By Definition~\ref{def:LagrangiankCobordism}\eqref{itm:verticalBoundary} we have $d\Psi_q(w')\in T_{\Psi(q)}\vlg_B$, so
	\begin{equation}
		d\Psi_q(v)=d\Psi_q(w)-d\Psi_q(w')\notin T_{\Psi(q)}\vlg_B
	\end{equation}
	since $T_{\Psi(q)}\vlg_B$ is a vector subspace. This proves \eqref{eq:almost-vertical-transversality}.
	We now show that
	\begin{equation}
		\label{eq:actual-vertical-transversality}
		d(f^\ast\Psi)_{(q,\xi)}(v,(u,\eta))\notin T_{f^\ast\Psi(q,\xi)}\vlg_A.
	\end{equation}
	We have two cases. The first case is when $d(\Psi_X)_q(v)\notin T_{\Psi_X(q)}\vlg$. Since $(f^\ast\Psi)_X=\Psi_X$ we have
	\begin{equation}
		d\bigl((f^\ast \Psi)_X\bigr)_{(q,\xi)}(v,(u,\eta))
		= d(\Psi_X)_q(v) \notin T_{\Psi_X(q)}\vlg
		= T_{(f^\ast \Psi)_X(q,\xi)} \vlg.
	\end{equation}
	In particular, we get \eqref{eq:actual-vertical-transversality}.

	The second case is when $d(\Psi_X)_q(v)\in T_{\Psi_X(q)}\vlg$. Since $d\Psi_q(v) \notin T_{\Psi(q)}\vlg_B$, we must have
	\begin{equation}
		\label{eq:nonvanishing-vertical}
		\eta \neq 0.
	\end{equation}
	For every \(t\in T_q\tns{C}{_0^1}(L)\), we have \(d\Psi_q(t)\in T_{\Psi(q)}\vlg_B\) by Definition~\ref{def:LagrangiankCobordism}\eqref{itm:verticalBoundary}. Hence, the Lagrangian condition on \(\Psi\) gives
	\begin{align}
		0
		 & =\omega_{X_B}\bigl(d\Psi_q(v),d\Psi_q(t)\bigr)
		\\
		 & =\omega_{T^\ast B}\bigl((df_a(u),\eta),(d(\Psi_B)_q(t),0)\bigr)
		+\omega_X\bigl(d(\Psi_X)_q(v),d(\Psi_X)_q(t)\bigr).
	\end{align}
	The second term vanishes because both \(d(\Psi_X)_q(v)\) and \(d(\Psi_X)_q(t)\) lie in \(T_{\Psi_X(q)}\vlg\), and \(\vlg\) is Lagrangian. Thus,
	\begin{equation}
		\eta(d(\Psi_B)_q(t)) = \omega_{T^\ast B}\bigl((df_a(u),\eta),(d(\Psi_B)_q(t),0)\bigr)=0.
	\end{equation}
	Since this holds for every $t\in T_q \tns{C}{_0^1}(L)$, we get
	\begin{equation}
		\label{eq:vanishing-vertical}
		\eta\mid_{d(\Psi_B)_q(T_q\tns{C}{_0^1}(L))}\equiv 0.
	\end{equation}
	Assume by contradiction that \eqref{eq:actual-vertical-transversality} fails.
	Then the \(T^\ast A\) component of the lefthand side of \eqref{eq:actual-vertical-transversality} is tangent to the zero section, hence its vertical part vanishes. We have $d(df^*)_\xi(u,\eta)=(u,df^*_a\eta)$, hence
	\begin{equation}
		\label{eq:vanishing-vertical2}
		df_a^\ast(\eta)=0
		\qquad \Longrightarrow \qquad
		\eta|_{df_a(T_aS^j(A))} \equiv 0.
	\end{equation}
	The transversality of \(f\) and \(\Psi_B\) together with \eqref{eq:vanishing-vertical} and \eqref{eq:vanishing-vertical2} gives \(\eta \equiv 0\), a contradiction to \eqref{eq:nonvanishing-vertical}. Therefore \eqref{eq:actual-vertical-transversality} holds, proving Definition \ref{def:LagrangiankCobordism}\eqref{itm:verticalTransversality}.
\end{proof}

\begin{lemma}
	\label{lem:weakly-smooth-pullback}
	Let $W=[\Psi:L\to X_B]$ be a Lagrangian cobordism over $B$
	relative to $\vlg$, and let $f:A\to B$ be weakly smooth as in
	Definition~\ref{def:smoothMaps}. Suppose that $f$ and $\Psi_B$
	satisfy~\eqref{eq:transversality-tangent-spaces1}
	and~\eqref{eq:transversality-tangent-spaces2} whenever
	$f(a)=\Psi_B(q)$. Then formula~\eqref{eq:pullback} defines a
	Lagrangian cobordism $f^\ast W$ over $A$ relative to $\vlg$.
\end{lemma}

\begin{proof}
	Since $\Psi_B$ is a b-submersion,
	Lemma~\ref{lem:weakly-smooth-fiber-product} shows that
	$f^\ast L=L\times_B A$ is a manifold with corners and its
	projection to $A$ is a b-submersion. Its local vertical boundary
	components at $(q,a)$ correspond to the local boundary components
	of $L$ at $q$ that are vertical for $\Psi_B$.
	The map
	\[
		f^\ast\Psi:f^\ast L\to X_A,\qquad
		(q,a)\mapsto
		\bigl((df_a)^\ast\Psi_{T^\ast B}(q),\,\Psi_X(q)\bigr)
	\]
	is weakly smooth. Since the boundary
	coordinates of $X_A$ come from $A$ and the projection of
	$f^\ast\Psi$ to $A$ is a b-submersion, $f^\ast\Psi$ is
	smooth in the sense of Definition~\ref{def:smoothMaps}.

	With the tangent-space description~\eqref{eq:fiber-product-tangent-space}
	and the identification of the local vertical boundary components
	above, the proof of Proposition~\ref{prop:pullback} shows that
	$f^\ast\Psi$ is a Lagrangian immersion satisfying the required
	vertical boundary conditions.
\end{proof}

\begin{lemma}
	\label{lem:restriction_to_open_subset}
	Let  $W \subset X_B$ be a Lagrangian cobordism over $B$ relative to $\vlg$, represented by $\Psi: L \to X_B$, and let $f: A \to B$ be an open embedding. Then the canonical map $\mathrm{pr}_1: f^\ast L \to L$ is an open embedding and satisfies
	\[
		\widehat f \circ f^\ast \Psi \;=\; \Psi \circ \mathrm{pr}_1 .
	\]
	That is, if $\hat f: T^\ast A \to T^\ast B$ is the induced open embedding of cotangent bundles, and $\widehat f: X_A \to X_B$ is the open embedding defined by $\widehat f = \hat f \times \mathrm{id}_X$, then $f^\ast \Psi$ is the unique map making the diagram commute:
	\begin{equation}
		\begin{tikzcd}
			{f^\ast{L}}
			\arrow[d,"\mathrm{pr}_1" left]
			\arrow[r, "f^\ast \Psi", dashed ]
			&
			X_{A}
			\arrow[d, "\widehat f"]
			\\
			L
			\arrow[r, "\Psi" below]
			&
			X_{B}
		\end{tikzcd}
	\end{equation}
\end{lemma}

\begin{proof}
	We have a commutative diagram
	\begin{equation}
		\begin{tikzcd}
			f^\ast L
			\arrow[d, "\mathrm{pr}_2"']
			\arrow[r, "\mathrm{pr}_1"]
			&
			L
			\arrow[d, "\Psi_{T^\ast B}"]
			\arrow[r, "\Psi_X  "]
			&
			X
			\\
			f^\ast T^\ast B
			\arrow[r, "\tilde f "]
			\arrow[d, "df^\ast"']
			&
			T^\ast B
			&
			\\
			T^* A
			\arrow[ur, "\hat f "']
		\end{tikzcd}
	\end{equation}
	which yields
	\begin{align}
		\Psi \circ \mathrm{pr}_1
		 & = (\Psi_{T^\ast B} \circ \mathrm{pr}_1, \Psi_X \circ \mathrm{pr}_1)
		\\
		 & = (\hat f \circ df^\ast \circ \mathrm{pr}_2, \Psi_X \circ \mathrm{pr}_1)
		\\
		 & = \widehat{f} \circ (df^\ast \circ \mathrm{pr}_2, \Psi_X \circ \mathrm{pr}_1)
		= \widehat{f} \circ f^\ast \Psi. \qedhere
	\end{align}
\end{proof}

\begin{lemma}
	\label{lem:cobordismOverProduct}
	Let $ W \subset X_{A \times B} $ be a Lagrangian cobordism over $ A \times B $ relative to $ \vlg $, represented by $ \Psi: L \to X_{A \times B} $.
	Fix $ b_0 \in B $ and define $ f:A \to A \times B $ by $f(a)=(a,b_0)$. Suppose that $f$ is transverse to $W$.
	Then $f^\ast \Psi$ is the unique map making the following diagram commute:
	\begin{equation}
		\begin{tikzcd}
			f^* L
			\arrow[d,"\mathrm{pr}_1" left]
			\arrow[drr, "f^\ast \Psi", dashed ]
			\\
			L
			\arrow[r, "\Psi" below]
			&
			X_{A \times B}
			\arrow[r, "\mathrm{pr}_{X_A}" below]
			&
			X_A
		\end{tikzcd}
	\end{equation}
\end{lemma}

\begin{proof}
	Since $f$ is constant in the $B$-coordinate, we have a commutative diagram
	\begin{equation}
		\begin{tikzcd}
			f^\ast L
			\arrow[d, "\mathrm{pr}_2"']
			\arrow[r, "\mathrm{pr}_1"]
			&
			L
			\arrow[d, "\Psi_{T^\ast (A\times B)}"]
			\arrow[r, "\Psi_X  "]
			&
			X
			\\
			f^\ast T^\ast (A \times B)
			\arrow[r, "\tilde f "]
			\arrow[d, "df^\ast"']
			&
			T^*(A \times B)
			\arrow[dl, "\mathrm{pr}_{T^*A}"]
			&
			\\
			T^* A
		\end{tikzcd}
	\end{equation}
	which gives
	\begin{align}
		\mathrm{pr}_{X_A}\circ \Psi \circ \mathrm{pr}_1
		 & = (\mathrm{pr}_{T^* A}\circ \Psi_{T^\ast (A\times B)} \circ \mathrm{pr}_1, \Psi_X \circ \mathrm{pr}_1)
		\\
		 & = (df^\ast \circ \mathrm{pr}_2, \Psi_X \circ \mathrm{pr}_1) = f^* \Psi.
		\qedhere
	\end{align}
\end{proof}

\subsection{Lagrangian boundary conditions}
\label{sec:bndry}

\begin{prop}
	\label{prop:canonicalBoundaryImmersion}
	Let $W \subset X_B$ be a Lagrangian cobordism over $B$ relative to $\vlg$ represented by $\Psi: L \to X_B$. Let $0 \le h \le \dim(B)$ and let $\iota=\iota_h: C_h(B)\to B$ be the canonical immersion. The corners of $W$ satisfy the following Lagrangian boundary conditions in $X_B$:
	\begin{enumerate}
		\item \label{itm:horizontalCorners}
		      The corners $\tensor{C}{_h^0}(W) \subset X_B$ are canonically immersed in $\iota_\ast(\iota^\ast W)\subset X_B$, i.e.\ there exists an immersion $\tensor{\beta}{_h^0}$ such that the following diagram commutes:
		      \begin{equation}
			      \begin{tikzcd}
				      \tensor{C}{_h^0}(L)
				      \arrow[d, "\tensor{\iota}{_h^0}"']
				      \arrow[r, "\tensor{\beta}{_h^0}"]
				      &
				      \iota_\ast(\iota^\ast L)
				      \arrow[d, "\iota_\ast(\iota^\ast \Psi)"]
				      \\
				      L
				      \arrow[r, "\Psi"]
				      &
				      X_B
			      \end{tikzcd}
		      \end{equation}

		\item \label{itm:mixedCorners}
		      The corners $\tensor{C}{_h^1}(W) \subset X_B$ immerse in the Lagrangian $\iota_\ast\bigl(\vlg_{C_h(B)}\bigr)\subset X_B$.
	\end{enumerate}
\end{prop}

\begin{remark}
	\label{rem:conormalAlongIota}
	Note that
	\begin{equation}
		\iota_\ast\bigl(\vlg_{C_h(B)}\bigr)
		\;=\;
		\iota^\ast T^\ast B \times_{T^\ast C_h(B)} \vlg_{C_h(B)}
	\end{equation}
	is canonically isomorphic to the product of the conormal bundle along $\iota$ with $\vlg$,
	\begin{equation}
		N^\ast\!\iota \times \vlg
		\;=\;
		\tilde\iota\bigl(\ker(d\iota^\ast)\bigr) \times \vlg
		\;\subset\; X_B.
	\end{equation}
\end{remark}

\begin{proof}
	We prove \eqref{itm:horizontalCorners}.   Recall that $\iota^\ast L$ and $\iota_\ast \iota^\ast L$ are defined as the fiber products
	\begin{equation}
		\begin{tikzcd}
			\iota^\ast L \arrow[r, "\text{pr}_1"] \arrow[d, "\text{pr}_2"] & L \arrow[d, "\Psi_{T^\ast B}"] \\
			\iota^\ast T^\ast B \arrow[r, "\tilde \iota"]  & T^\ast B
		\end{tikzcd}
		\quad \quad \quad \quad
		\begin{tikzcd}
			\iota_\ast\iota^\ast L \arrow[r, "\text{pr}_3"] \arrow[d, "\text{pr}_4"] & \iota^\ast L \arrow[d, "(\iota^\ast\Psi)_{T^\ast C_h(B)}"] \\
			\iota^\ast T^\ast B \arrow[r, "d\iota^\ast"]                             & T^\ast C_h(B)
		\end{tikzcd}
	\end{equation}
	The map
	\begin{equation}
		\hat \beta: \iota^\ast L \to \iota_\ast \iota^\ast L,
		\quad \quad
		\hat \beta(q, \xi) = ((q,\xi),\xi),
	\end{equation}
	is well defined since $\Psi_{T^\ast B}(q) = \tilde\iota(\xi)$ implies $(\iota^\ast \Psi)_{T^\ast C_h(B)}(q,\xi)=d\iota^\ast (\xi)$ by definition. We have
	\begin{equation}
		(\iota_\ast \iota^\ast \Psi)_X(\hat\beta(q,\xi)) = (\iota^\ast \Psi)_X(q,\xi) = \Psi_X(q) ,
	\end{equation}
	and
	\begin{equation}
		(\iota_\ast \iota^\ast \Psi)_{T^\ast B} (\hat\beta(q,\xi))
		= \tilde\iota (\xi)
		= \Psi_{T^ \ast B}(q),
	\end{equation}
	thus the following diagram commutes:
	\begin{equation}
		\begin{tikzcd}
			\iota^\ast L
			\arrow[r, "\hat \beta"] \arrow[d, "\text{pr}_1"] & \iota_\ast\iota^\ast L \arrow[d, "\iota_\ast\iota^\ast\Psi"] \\
			L \arrow[r, "\Psi"]                 & X_B
		\end{tikzcd}
	\end{equation}
	Finally, note that we have a canonical isomorphism
	\begin{equation}
		\iota^\ast L = L \times_{T^\ast B} \iota^\ast T^\ast B \simeq L \times_{T^\ast B}(T^\ast B \times_B C_h(B)) \simeq L \times_B C_h(B)
	\end{equation}
	so by Lemma \ref{lem:horizontalPullback} we have a canonical isomorphism $ \vp: \tns{C}{_h^0}(L) \to \iota^\ast L$ which satisfies $\text{pr}_1 \circ \vp = \tns{\iota}{_h^0}$. Part \eqref{itm:horizontalCorners} follows by taking $\beta = \hat\beta\circ \vp$.
	We have a commutative diagram
	\begin{equation}
		\begin{tikzcd}
			\tns{C}{_h^1}(L)
			\ar[dr, "\tns{\iota}{_h^1}"]
			\ar[d, "\beta"]
			\ar[r, "\alpha"]
			&
			\tns{C}{_h^0}(L)
			\ar[d, "\tns{\iota}{_h^0}"]
			\\
			\tns{C}{_0^1}(L)
			\ar[r, "\tns{\iota}{_0^1}"']
			&
			L
		\end{tikzcd}
	\end{equation}
	where $\alpha$ and $\beta$ are the immersions forgetting the vertical and horizontal boundary components, respectively. Thus, for $q \in \tns{C}{_h^1}(L)$ we have $\Psi(\tns{\iota}{_h^1}(q)) = \Psi(\tns{\iota}{_0^1}(\beta(q))) \in \vlg_B$, and
	\begin{equation}
		\Psi_{T^*B}(\tns{\iota}{_h^1}(q))
		=
		\Psi_{T^*B}(\tns{\iota}{_h^0}(\alpha(q)))
		=
		\Psi_{T^*B}(\text{pr}_1 ( \vp(\alpha(q))))
		=
		\tilde \iota (\text{pr}_2 ( \vp(\alpha(q)))),
	\end{equation}
	so $\Psi(\tns{\iota}{_h^1}(q)) \in N^* \iota \times \vlg $ (see Remark \ref{rem:conormalAlongIota}), which proves \eqref{itm:mixedCorners}.
\end{proof}

\subsection{Suspension cobordisms of exact families}
\label{sec:suspension}
If $f:X \times Y \to Z$ is a map, we denote by $f_x:Y \to Z$, $f_y:X \to Z$ the maps $f(x,\cdot)$ and $f(\cdot,y)$. In what follows, we recall the notion of exact family of Lagrangian submanifolds from~\cite{Solomon13}.
Let $(X,\omega_X)$ be a symplectic $2n$-dimensional manifold, and let $\vlg \subset X$ be a smoothly embedded closed Lagrangian. Let $L$ be an $n$-dimensional compact connected manifold with non-empty boundary, and let $\cc{L} = \cc{L}(X,L; \vlg)$ denote the space of freely immersed Lagrangians $ W \subset X$ with boundary immersed in $\vlg$, such that for all $p \in \partial W$
\begin{equation}
	\label{eq:suspensionTransverse}
	T_p  W \neq T_{p_0} \vlg.
\end{equation}
where $p_0$ is the image of $p$ under the immersion of $W$.
\begin{definition}
	Let $B$ be a contractible manifold with corners. A map $W: B \to \cc{L}$ is called \emph{smooth} if there exists a smooth map $\psi:L \times  B \to X$ (called a \emph{lifting}), such that $\psi_b: L \to X$ is a Lagrangian immersion representing $W(b)$ for every $b\in B$.
\end{definition}

\begin{definition}
	\label{def:exactFamily}
	Let $W: B \to \cc{L}$ be a smooth family.  We say that $W$ is \emph{exact relative to $\vlg$} if for some (equivalently, any) lifting  $ \psi: L \times B \to X $, there exists a smooth function $h: L \times TB \to \bb{R}$ such that for all $(q,\xi)\in L \times TB$:
	\begin{enumerate}
		\item $h_q: TB \to \bb{R}$ is fiberwise linear.
		\item $h_\xi: L \to \bb{R}$ satisfies $h_\xi|_{\partial L}\equiv 0$.
		\item \label{itm:exactFamily} $ dh_\xi = \xi \lrcorner (\psi^\ast \omega_X) \in \Omega^1(L) $.
	\end{enumerate}
	$h$ is called the \emph{associated family of Hamiltonian functions}.
\end{definition}

\begin{definition}
	Let $h: L \times TB \to \bb{R}$ be a smooth map such that $h_q: TB \to \bb{R}$ is fiberwise linear. We denote by $h^\vee: L \times B \to T^*B$ the fiberwise dual given by $h^\vee(q,b)(\xi)=h(q,\xi)$ for $ \xi \in T_bB $.
\end{definition}

\begin{lemma}
	\label{lem:suspensionLagrangianImmersion}
	Let $B$ be a contractible manifold with corners. Let $W: B \to \cc{L}$ be a smooth family of immersed Lagrangians in $X$, exact relative to $\vlg$. Let $\psi: L \times B \to X  $ be a lifting of $W$, and  let $h: L \times TB \to \bb{R}$ be the associated family of Hamiltonian functions.  The map
	\begin{equation}
		\label{eq:suspensionImmersion}
		\Psi: L \times B \to  T^\ast B \times X,
		\quad \quad
		\Psi=(h^\vee, \psi)
	\end{equation}
	is a Lagrangian immersion.
\end{lemma}

\begin{proof}
	Consider the unique lift $\hat{h}: T(L \times B) \to \bb{R}$ of $h$ which is constant along the fibers of $TL$. Denote by $\hat{h}^\vee: L \times B \to T^\ast(L \times B)$ the corresponding $1$-form on $L\times B$. Let $\pi:L \times B \to B$ be the projection, then $h^\vee$, $\hat{h}^\vee$ factor canonically via a section $s$:
	\begin{equation}
		\begin{tikzcd}
			&
			L \times B
			\ar[dl, "\hat{h}^\vee"']
			\ar[dr, "h^\vee"]
			\ar[d, "s"]
			&
			\\
			T^\ast (L \times B)
			&
			\pi^\ast T^\ast B
			\ar[r,"\tilde \pi"]
			\ar[l, "d\pi^\ast"']
			&
			T^\ast B
		\end{tikzcd}
	\end{equation}
	We first note that $\Psi$ is an immersion. Let $0 \neq (v,u)\in T_{(q,b)}(L\times B) \simeq T_q L \oplus T_b B $. If $u \neq 0$ then $dh^\vee|_{(q,b)}(v,u) \neq 0$ since $\pi = \pi_B\circ h^\vee$ where $\pi_B:T^\ast B \to B$ is the projection. If $u =0$, then $d\psi_{(q,b)}(v,0)\neq0$ since $\psi_b: L \to X$ is an immersion.

	We now prove that $\Psi$ is Lagrangian. First, if $\omega_{T^\ast B}$ is the symplectic form on $T^\ast B$, we claim that
	\begin{equation}
		\label{eq:suspensionPullbackform1}
		d\hat{h}^\vee =
		(h^\vee)^\ast \omega_{T^\ast B}.
	\end{equation}
	Indeed, let $\theta_B$ be the canonical 1-form of $T^\ast B$. By equation \eqref{eq:canonical1Forms} with $ f=\pi $ and $ A=L \times B $ we get
	\begin{equation}
		(h^\vee)^\ast \theta_B= s^\ast \tilde \pi^\ast \theta_B = s^\ast (d\pi^\ast)^\ast \theta_{L\times B} =(\hat{h}^\vee)^\ast \theta_{L \times B}= \hat{h}^\vee,
	\end{equation}
	taking exterior derivative gives \eqref{eq:suspensionPullbackform1}. Thus, to show $\Psi$ is Lagrangian, it suffices to show that
	\begin{equation}
		\label{eq:suspensionPullbackform2}
		d\hat{h}^\vee = - \psi^\ast \omega_X
		\quad
		\in \Omega^2(L\times B).
	\end{equation}
	Let $\eta_1=(\eta_1',0)$ and $ \eta_2=(\eta_2',0) $ in $ \Gamma(T(L\times B))$, where $\eta_1',\eta_2' \in \Gamma(TL)$. We have
	\begin{equation}
		\label{eq:suspensionPullbackform4}
		d\hat{h}^\vee(\eta_1,\eta_2) =0 =-\psi^\ast \omega_X(\eta_1,\eta_2),
	\end{equation}
	the first equality follows since $\hat h$ is constant on the fibers of $TL$, and the second since $\psi_b$ is a Lagrangian immersion for every $b\in B$.
	Let $\eta=(\eta',0)$ and $ \xi=(0, \xi') $ in $ \Gamma(T(L\times B))$, where $\xi'\in \Gamma(TB)$ and $\eta'\in \Gamma(TL)$. We have
	\begin{align}
		\label{eq:suspensionPullbackform3}
		d\hat{h}^\vee (\eta,\xi)
		= \partial_\eta \hat{h}^\vee(\xi) - \partial_\xi \hat{h}^\vee(\eta)
		= \partial_\eta \hat{h}^\vee(\xi)
		= dh_\xi(\eta)
		= -\psi^\ast \omega_X(\eta, \xi)
	\end{align}
	where the second equation in \eqref{eq:suspensionPullbackform3} follows since $\hat{h}^\vee(\eta)$ is constant on $L\times B$ (since $\hat{h}$ is constant on the fibers of $TL$), and the last equation in \eqref{eq:suspensionPullbackform3} follows from condition \eqref{itm:exactFamily} in Definition \ref{def:exactFamily}.

	Let $\xi_1=(0,\xi_1')$ and $\xi_2=(0,\xi_2')$, where
	$\xi_1',\xi_2'\in\Gamma(TB)$. For $\eta=(\eta',0)$ as above,
	we have $[\eta,\xi_i]=0$.
	By \eqref{eq:suspensionPullbackform4} and \eqref{eq:suspensionPullbackform3},
	we have $(d\hat{h}^\vee+\psi^\ast\omega_X)(\eta,\cdot)=0$.
	Thus, when we evaluate $d(d\hat{h}^\vee+\psi^\ast\omega_X)=0$ on
	$(\eta,\xi_1,\xi_2)$, both derivatives in the $\xi_i$ directions
	vanish, as does the term pairing $[\xi_1,\xi_2]$ with $\eta$.
	Since $[\eta,\xi_i]=0$, only the derivative in the $\eta$ direction
	survive, giving
	\[
		\partial_\eta\bigl(d\hat{h}^\vee(\xi_1,\xi_2)
		+\psi^\ast\omega_X(\xi_1,\xi_2)\bigr)=0.
	\]
	Since $L$ is connected, the expression in parentheses is constant on
	$L\times\{b\}$ for each $b\in B$. Fix $q\in\partial L$. On
	$\{q\}\times B$, the pullback of $\hat{h}^\vee$ vanishes by the
	boundary condition on $h$, and the pullback of $\psi^\ast\omega_X$
	vanishes since $\psi(q,\cdot)$ takes values in the Lagrangian $\vlg$.
	Thus the constant is zero, so
	\begin{equation}
		\label{eq:suspensionPullbackform5}
		d\hat{h}^\vee(\xi_1,\xi_2)=-\psi^\ast\omega_X(\xi_1,\xi_2).
	\end{equation}
	Equations \eqref{eq:suspensionPullbackform4},
	\eqref{eq:suspensionPullbackform3} and \eqref{eq:suspensionPullbackform5}
	imply \eqref{eq:suspensionPullbackform2}, which proves the Lemma.
\end{proof}

\begin{prop}
	Let $B$ be a contractible manifold with corners. Let $W: B \to \cc{L}$ be a smooth family of immersed Lagrangians in $X$, exact relative to $\vlg$. Let $\psi: L \times B \to X  $ be a lifting of $W$, and  let $h: L \times TB \to \bb{R}$ be the associated family of Hamiltonian functions. Then \eqref{eq:suspensionImmersion} represents a Lagrangian cobordism $ SW \subset X_B $ relative to $\vlg$, the \emph{suspension of $ W $},  independent of $ \psi $.
\end{prop}

\begin{proof}
	Since $ \Psi $ is a Lagrangian immersion by Lemma \ref{lem:suspensionLagrangianImmersion}, it remains to verify \eqref{itm:bsubmersion}-\eqref{itm:verticalCorners} in Definition \ref{def:LagrangiankCobordism}. The map $ \Psi_B: L \times B \to B $ is the projection, so it satisfies Definition \ref{def:LagrangiankCobordism}\eqref{itm:bsubmersion}. By assumption, $ \psi_b $ is an immersion of  $ \partial L$  into $ \vlg $ and $ h_\xi $ vanishes on $ \partial L$. Therefore $\Psi$ maps
	\begin{equation}
		\tns{C}{_0^1}(L\times B) = \partial L \times B
		\longrightarrow
		0_B \times \vlg \subset X_B.
	\end{equation}
	so Definition \ref{def:LagrangiankCobordism}\eqref{itm:verticalBoundary} holds. By \eqref{eq:suspensionTransverse}, for every $ q\in \partial L $ and $ b\in B $  there exists $ v_b \in T_q L $ with $ d\psi_b(v_b) \notin T_{\psi(q)} \vlg $. Hence, with  $v=(v_b,0) \in T_{(q,b)} (L \times B) $,
	\begin{equation}
		d(\Psi_B)_{(q,b)}(v) =0,
		\quad \quad
		d(\Psi_X)_{(q,b)}(v) = d\psi_b(v_b) \notin T_{\psi_b(q)}\vlg,
	\end{equation}
	so Definition \ref{def:LagrangiankCobordism}\eqref{itm:verticalTransversality} holds.  Finally,
	\begin{equation}
		\tns{C}{_0^v}(L\times B) =  \tns{C}{_0^v}(L) \times B =\emptyset, \quad \quad \forall v\geq 2,
	\end{equation}
	since $ L $ is a manifold with (only) boundary; thus Definition \ref{def:LagrangiankCobordism}\eqref{itm:verticalCorners} holds.

	For independence of the lifting, let $\psi': L \times B \to X$ be another lifting of $W$. By freeness and
	\cite[proof of Theorem~2.12]{SolomonYuval20},  there exists a smooth map $\vp: L \times B \to L$, with $\vp_b: L \to L$ a diffeomorphism, such that $\psi_b' = \psi_b \circ \vp_b$.
	The associated Hamiltonians satisfy $h'_\xi = h_\xi \circ \vp_b$
	for $\xi \in T_bB$, since the two functions have the same differential
	by Definition~\ref{def:exactFamily}\eqref{itm:exactFamily} and both vanish on $\partial L$.  Hence $\Psi' = \Psi \circ (\vp \times \mathrm{id}_B) $, so $\Psi$ and $\widetilde \Psi$  represent the same immersed Lagrangian of type $L \times B$.
\end{proof}

\subsection{Generalized cobordisms: vertical systems}
\label{sec:verticalSystems}
There are at least three natural extensions of Definition \ref{def:LagrangiankCobordism}:
\begin{itemize}
	\item allow the vertical boundary condition $\vlg\subset X$ to be \emph{immersed} (not necessarily embedded);
	\item allow $\vlg$ to \emph{vary along the base}: i.e.\ the vertical boundary of $W$ lies in a nontrivial Lagrangian cobordism over $B$ rather than in the product $\vlg_B$;
	\item allow \emph{higher vertical corners}, permitting $\tensor{C}{_0^v}(W)\neq\varnothing$ for $v\ge 2$, so that we can consider cobordisms between Lagrangians in $X$ with corners.
\end{itemize}
The next definition packages these three features in a single framework.

\begin{definition}
	\label{def:verticalSystem}
	Let $ \ell \in \bb{Z}_{\geq 0} $. A \emph{vertical system $W$ of immersed Lagrangian cobordisms over $ B $ of length $\ell$} is an equivalence class of pairs $(\Psi,\beta)$, where
	\begin{enumerate}

		\item $\Psi=\bigl(\Psi^v: L^v \to X_B\bigr)_{0\le v\le \ell}$ is a sequence of Lagrangian immersions such that each base projection
		      \[
			      \Psi_B^v: L^v \longrightarrow B
		      \]
		      is a b\mbox{-}submersion.

		\item $\beta$ is a family of \emph{vertical boundary immersions},
		      \begin{equation}
			      \beta=\bigl\{\;\tensor{\beta}{^{\,vv'}}:\tensor{C}{_0^{v'}}(L^v)\to L^{\,v+v'}\;\bigm|\;0\le v, v',\,\, v+v'\le \ell\;\bigr\}
		      \end{equation}
		      where the corner $\tensor{C}{_0^{v'}}(L^v)$ is taken with respect to $\Psi_B^v$, such that the diagrams commute:
		      \begin{equation}
			      \begin{tikzcd}
				      \tensor{C}{_0^{v'}}(L^v)
				      \arrow[d, hook', "\jmath"']
				      \arrow[r, "\tensor{\beta}{^{\,vv'}}"]
				      &
				      L^{\,v+v'}
				      \arrow[d, "\Psi^{\,v+v'}"]
				      \\
				      L^v
				      \arrow[r, "\Psi^{\,v}"']
				      &
				      X_B
			      \end{tikzcd}
		      \end{equation}
		      For $v'=0$ we take $\beta^{v0}=\mathrm{id}_{L^v}$. For $v'>0$ and every $q\in \tensor{C}{_0^{v'}}(L^v)$, we require
		      \begin{equation}
			      d\Psi^{\,v+v'}_{\tensor{\beta}{^{\,vv'}}(q)}\!\bigl(T_{\tensor{\beta}{^{\,vv'}}(q)}L^{\,v+v'}\bigr)
			      \;\neq\;
			      d\Psi^{\,v}_{\jmath(q)}\!\bigl(T_{\jmath(q)}L^{\,v}\bigr).
		      \end{equation}
	\end{enumerate}
	Two pairs $ (\Psi, \beta) $ and $ (\widetilde\Psi, \tilde \beta )$ (with the same $\ell$) are \emph{equivalent} if for each $   0 \leq v \leq \ell$ there exists a diffeomorphism $  \vp^v: L^v \to L^v $ such that
    \begin{equation}
        \begin{aligned}
            \widetilde\Psi^{\,v}
            &=\Psi^{\,v}\circ\vp^{\,v},\\
            \vp^{\,v+v'}\circ\widetilde{\tensor{\beta}{^{\,vv'}}}
            &=\tensor{\beta}{^{\,vv'}}\circ\tensor{C}{_0^{v'}}(\vp^{\,v})
            \qquad\text{for all }v,v'\text{ with }v+v'\le\ell,
     \end{aligned}
    \end{equation}
	where $\tensor{C}{_0^{v'}}(\vp^{\,v}):\tensor{C}{_0^{v'}}(L^v)\to \tensor{C}{_0^{v'}}(L^v)$ is the diffeomorphism induced on the $v'$\nobreakdash-vertical corner by $\vp^{\,v}$, and $\jmath$ denotes the canonical corner inclusion. That is, the following diagram commutes:
	
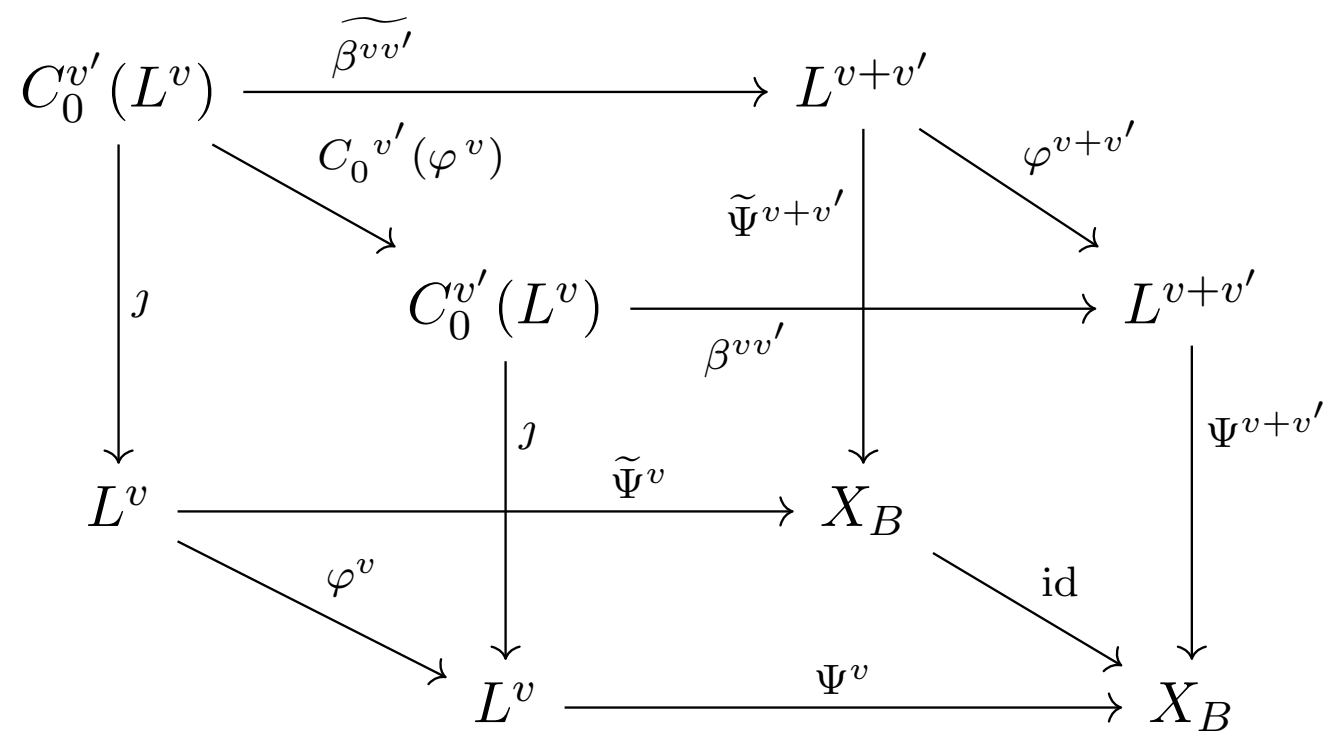
\begin{equation}
		\begin{tikzcd}
			C_0^{v'}(L^v)
			\arrow[dd, "\jmath"]
			\arrow[rr, "\widetilde{\beta^{vv'}}" near start]
			\arrow[rd, "\tensor{C}{_0^{v'}}(\vp^{\,v})" ]
			&
			&
			L^{v+v'}
			\arrow[dd, "\widetilde\Psi^{v+v'}"' near start]
			\arrow[rd, "\vp^{v+v'}"]
			&
			\\
			&
			C_0^{v'}(L^v)
			\arrow[rr, "\beta^{vv'}"' near start]
			\arrow[dd, "\jmath" near start]
			&
			&
			L^{v+v'}
			\arrow[dd, "\Psi^{v+v'}" near start]
			\\
			L^v
			\arrow[rr, "\widetilde\Psi^v" near end]
			\arrow[rd, "\vp^v"]
			&
			&
			X_B
			\arrow[rd, "\mathrm{id}"]
			&
			\\
			&
			L^v
			\arrow[rr, "\Psi^v"]
			&
			&
			X_B
		\end{tikzcd}
	\end{equation}
\end{definition}

Let $ W = [(\Psi,\beta)]$ be a vertical system of Lagrangian cobordisms over $ B $ of length $\ell$. A smooth map $ f: A \to B $ which is transverse to $ W $ (i.e. $f \pitchfork \Psi_B^v$ for all $0 \leq v \leq \ell$) determines a vertical system
\begin{equation}
	f^\ast W =[(f^*\Psi, f^*\beta)],
\end{equation}
where the immersions $ f^*\Psi $ are defined as in \eqref{eq:pullback}, and
\begin{equation}
	f^\ast \beta =
	\{
	\tns{f^\ast \beta}{^{vv'}}: \tns{C}{_0^{v'}}(\tns{f^\ast L}{^v}) \to \tns{f^\ast L}{^{v + v'}} \, : 0 \leq v,v'\,  \, 0 \leq v+v' \leq \ell
	\}
\end{equation}
are immersions canonically determined by $ \beta $.

\begin{remark}
	The cobordisms of a vertical system $ W $ can be conveniently organized in tensor notation (cf. Definition \ref{def:associatedTensors}).
	Let $ \iota_h: C_h(B) \to B $ be the canonical immersion, and denote by $ \tns{W}{_h^v} \subset X_{C_h(B)} $ the immersed Lagrangian represented by $ \iota_h^* \Psi^v $. The ``tensor'' $ \tns{W}{_h^v} $ packages all the Lagrangian cobordisms forming the vertical system $ W $; for fixed $ h_0$, the ``vector'' $ \tns{W}{_{h_0}^v} $ is the vertical system over $ C_h(B) $ obtained by restricting $ W $.

	If $ f:A \to B $ is a b-submersion transverse to $ W $, then, using the tensor-product notation of Corollary \ref{lem:fiberAndTensorProducts}, the pullback tensor is given by
	\begin{equation}
		f^\ast W = \Cd(A) \otimes_{C(B)}W.
	\end{equation}
\end{remark}

\section{Exact collars and compatible Weinstein neighborhoods}
\label{sec:collars}

\subsection{Setup and notation.}
In this section we study Lagrangians $W:S \to \cob(X,\vlg)$ over a simplicial subset $S \subset \Delta[n]$. To construct collars uniformly on all faces, we realize the bases of all simplices in $W$ inside a single ambient geometric simplex.

\label{sec:simplex-collars}
\begin{definition}
	\label{def:n-simplex}
	For $n\geq 0$ define the \emph{geometric $ n $-simplex}
	\begin{equation}
		\Delta^n := \left\{ (s_0,\dots,s_n) \in \bb{R}^{n+1} \mid s_0,\dots,s_n \geq 0 , \,\, s_0 + \dots +s_n = 1 \right\}.
	\end{equation}
	For an index set $\emptyset \neq I\subseteq [n]$, its $I$-face is
	\begin{equation}
		\Delta_{I}^n = \left\{ (s_0,\dots,s_n) \in \Delta^n \mid   s_j=0 \text{ for all } j\notin I \right\}.
	\end{equation}
	See Figure \eqref{fig: 3simplex}. Then $ \Delta_{I}^n  $ is a connected component of $ C_{n+1-|I|}(\Delta^n)$ and is canonically diffeomorphic to $\Delta^{|I|-1}$.
\end{definition}

\begin{definition}
	\label{def:family-of-faces}
	For $\emptyset \neq I \subseteq [n]$ with $I= \{i_0< i_1 < \dots < i_k \}$ define
	\begin{equation}
		\iota_I: [k] \to [n], \quad \quad
		\iota_I(j) = i_j.
	\end{equation}
	If $S \subset \Delta[n]$ is a simplicial subset, set
	\begin{equation}
		\cc{F}_S :=
		\left\{
		\emptyset\neq I \subseteq [n]
		\; | \;
		\iota_I \in S
		\right\}.
	\end{equation}
	Let $J \in \cc{F}_S$. If $I \subseteq J$, then $\iota_I$ is a face of $\iota_J$ in $\Delta[n]$. Since $S$ is a simplicial \emph{subset}, it is closed under faces, hence $\iota_I \in S$, so $\cc{F}_S$ is closed under taking subsets.
\end{definition}

Recall the notion of a Lagrangian $W$ over a simplicial set $S$ (Definition \ref{def:Lag-over-sSet}).
\begin{definition}
	\label{def:family-of-cobordisms}
	Let $S \subset \Delta[n]$ be a simplicial subset, and let $W: S \to \cob(X,\vlg)$ be a Lagrangian over $S$. For $I \in \cc{F}_S$, write $k=|I|-1$. Define
	\begin{equation}
		W_I  \subset X_{\Delta_I^n}
	\end{equation}
	to be the pullback of $W(\iota_I) \subset X_{\Delta^{k}}$ along the canonical isomorphism $\Delta_I^n \to \Delta^{k}$. Thus $W$ determines a family of cobordisms
	\begin{equation}
		\cc{F}_W :=
		\left\{
		W_I \subset X_{\Delta^n_I}
		\; | \;
		I \in \cc{F}_S
		\right\}.
	\end{equation}
\end{definition}

\begin{lemma}
	For every inclusion $I \subseteq J$ in $\cc{F}_S$ we have the face-compatibility
	\[
		W_J|_{\Delta_I^n} = W_I.
	\]
\end{lemma}

\begin{proof}
	Write $k = |I|-1$ and $l = |J|-1$. There is a unique order-preserving injection $\delta: [k] \to [l]$ such that $S(\delta)(\iota_J)=\iota_I$. Let $\delta: \Delta^k \to \Delta^l$ be the corresponding geometric face inclusion. We have a commutative diagram with vertical isomorphisms
	\begin{equation}
		\begin{tikzcd}
			\Delta_I^n
			\ar[d, "\phi_I"]
			\ar[r, "\iota"]
			&
			\Delta_J^n
			\ar[d, "\phi_J"]
			\\
			\Delta^k
			\ar[r, "\delta"]
			&
			\Delta^l
		\end{tikzcd}
	\end{equation}
	Since $W$ is a natural transformation, the diagram
	\begin{equation}
		\begin{tikzcd}
			S_l
			\ar[d, "S(\delta)"']
			\ar[r, "W_l"]
			&
			\cob(X,\vlg)_l
			\ar[d, "\delta^\ast"]
			\\
			S_k
			\ar[r, "W_k"]
			&
			\cob(X,\vlg)_k
		\end{tikzcd}
	\end{equation}
	commutes, and we get
	\begin{align}
		W_J|_{\Delta_I^n}
		=
		\iota^\ast \phi_J^\ast \big(W_l(\iota_J)\big)
		 & =
		\phi_I^\ast \delta^\ast \big(W_l(\iota_J)\big)
		\\
		 & =
		\phi_I^\ast \big( W_k(S(\delta)(\iota_J)) \big)
		=
		\phi_I^\ast \big( W_k(\iota_I) \big)
		=
		W_I.
		\qedhere
	\end{align}
\end{proof}

\begin{lemma}
	\label{lem:horn-filling}
	Let $W: \Lambda_i[n] \to \cob(X,\vlg)$ be a Lagrangian over the $(n,i)$-horn. Let $\widetilde{W} \in \cob(X,\vlg)_n$ be a Lagrangian cobordism over $\Delta^n$, such that for every $I \in \cc{F}_{\Lambda_i[n]}$ we have
	\begin{equation}
		\widetilde{W}|_{\Delta_I^n}
		=
		W_I.
	\end{equation}
	Then $\widetilde{W}$ determines a map $\widetilde{W}: \Delta[n] \to \cob(X,\vlg)$ extending $W$, i.e. \eqref{eq:hornFilling} commutes.
\end{lemma}

\begin{proof}
	Recall that Yoneda lemma gives an equality
	\begin{equation}
		\cob(X,\vlg)_n = \hom_{}(\Delta[n], \cob(X,\vlg)).
	\end{equation}
	Let  $\iota: [k] \to [n]$ be a non-degenerate simplex in $\Lambda_i[n]$ and set $I = \mathrm{im}(\alpha) \subset [n]$. We have
	\begin{equation}
		\widetilde{W}_k(\iota) = \widetilde{W}|_{\Delta_I^n} = W_I = W_k(\iota).
	\end{equation}
	If $\alpha$ is a degenerate simplex in $\Lambda_i[n]$, we can factor $\alpha = \iota \circ \rho$ with $\rho$ surjective and $\iota$ injective, and we get $\widetilde{W}_k(\alpha) = \rho^\ast \widetilde{W}_k(\iota) = \rho^\ast W(\iota) = W(\alpha)$.   Thus, $\widetilde{W}$ and $W$ agree on all simplices, and the diagram \eqref{eq:hornFilling} commutes.
\end{proof}

\begin{ntn}
	\label{ntn: coordinates_projection}
	Let $n \geq 1$. If $I = \{i_0 < \dots < i_k\} \subseteq [n]$ and $t=(t_0,\dots,t_n) \in \bb{R}^{n+1}$, write  $t_I = (t_{i_1},\dots , t_{i_k})$ for the projection to the coordinates in $I$. We denote points of $\Delta_I^n$ (Definition \ref{def:n-simplex}) by $s_{I}$, i.e. by dropping the vanishing coordinates. For example,  we write $(s_1,s_2) \in \Delta_{12}^2$.

	For $I \subsetneq  J \subsetneq K \subseteq [n]$ we will often write
	\begin{equation}
		(t_{J \setminus I}, t_{K \setminus J})
	\end{equation}
	instead of $t_{K \setminus I}$ to highlight the distinct roles of the coordinates indexed by $J \setminus I$ and $K \setminus J$, but the order of the coordinates is unchanged. For example, if
	\begin{equation}
		K = \{1,2,3,4\},
		\quad \quad
		J =\{1,3\},
		\quad \quad
		I = \{1\},
	\end{equation}
	then $t_{K\setminus J} = (t_2,t_4)$ and $t_{J \setminus I} = t_3$, so $(t_{J \setminus I}, t_{K \setminus J})$ denotes $(t_2,t_3,t_4)$  and not $(t_3,t_2,t_4)$.
	However, $(s_{J\backslash I}, t_{K \backslash J})$ still stands for $(s_3,t_2,t_4)$.
\end{ntn}

\subsection{Collars for Lagrangian cobordisms}
\label{sec:cobordism-collars}
For the remainder of this section, fix a simplicial subset $S \subset \Delta[n]$, a Lagrangian $W:S \to \cob(X,\vlg)$, and, for each $I \in \cc{F}_S$ a representing immersion $\Psi_I: L_I \to X_{\Delta_I^n}$ of $W_I$.

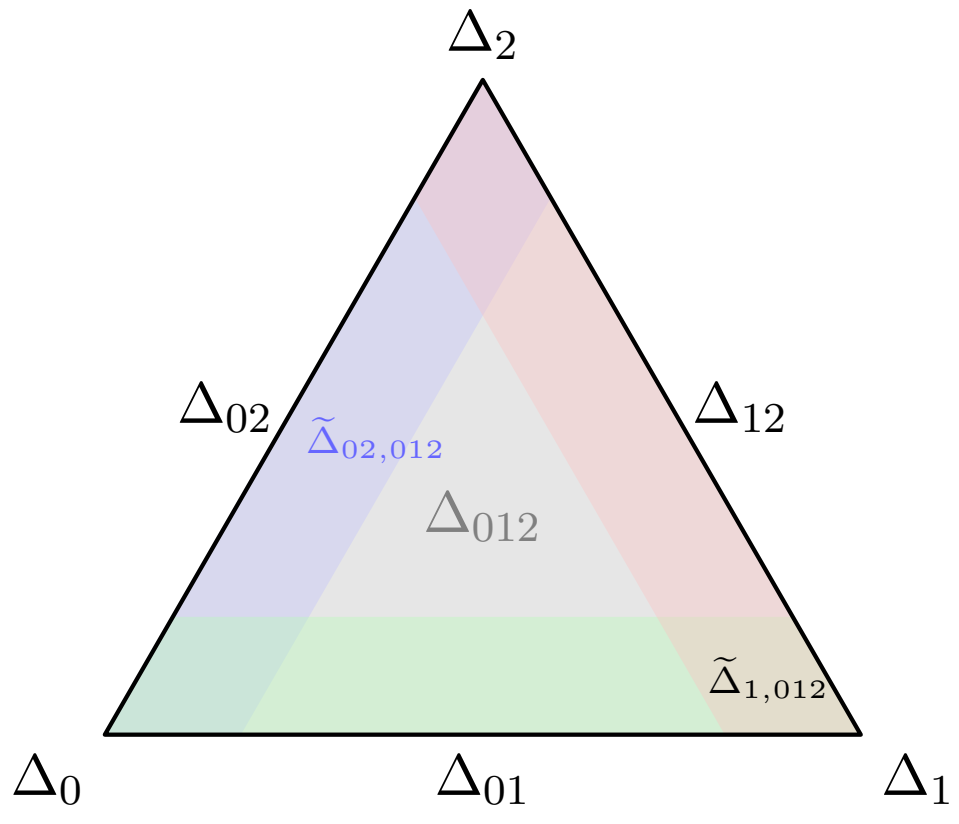
\begin{figure}[htbp]
	\centering
	\begin{tikzpicture}[scale=1.0, line join=round, line cap=round]
		\coordinate (A) at (0,2.887);
		\coordinate (B) at (-2.5,-1.443);
		\coordinate (C) at (2.5,-1.443);

		\def\t{0.18}

		\coordinate (AtoC) at ($(A)!\t!(C)$);
		\coordinate (BtoC) at ($(B)!\t!(C)$);
		\coordinate (BtoA) at ($(B)!\t!(A)$);
		\coordinate (CtoA) at ($(C)!\t!(A)$);
		\coordinate (CtoB) at ($(C)!\t!(B)$);
		\coordinate (AtoB) at ($(A)!\t!(B)$);

		\fill[gray!20] (A)--(B)--(C)--cycle;

		\fill[blue!25, opacity=0.35] (A)--(B)--(BtoC)--(AtoC)--cycle;  
		\fill[green!30, opacity=0.35] (B)--(C)--(CtoA)--(BtoA)--cycle; 
		\fill[red!25, opacity=0.35] (C)--(A)--(AtoB)--(CtoB)--cycle;   

		\draw[thick] (A)--(B)--(C)--cycle;
		\node[above] at (A) {$\Delta_2$};
		\node[below left] at (B) {$\Delta_0$};
		\node[below right] at (C) {$\Delta_1$};
		\node[above left=2pt, font=\tiny] at (C) {$\widetilde\Delta_{1,012}$};

		\node[left] at ($(A)!0.5!(B)$) {$\Delta_{02}$};
		\node[right=2pt, font=\tiny, text=blue!60] at ($(A)!0.55!(B)$) {$\widetilde\Delta_{02,012}$};
		\node[below] at ($(B)!0.5!(C)$) {$\Delta_{01}$};
		\node[right]  at ($(C)!0.5!(A)$) {$\Delta_{12}$};
		\node[text=gray] at (barycentric cs:A=1,B=1,C=1) {$\Delta_{012}$};

	\end{tikzpicture}
	\caption{Corners and collars of the 2-simplex.}
	\label{fig: 3simplex}
\end{figure}

\begin{definition}
	Let $ n>0 $. For $ \emptyset \neq I \subseteq J \subseteq [n]$ and $\varepsilon>0$:
	\begin{enumerate}
		\item The $\varepsilon$-\emph{collar of $\Delta_I^n$ in $\Delta_J^n$} is
		      \begin{equation}
			      \widetilde \Delta_{IJ} ^n(\varepsilon)
			      :=
			      \left\{
			      s_{J} \in \Delta_J^n
			      \,\, | \,\,
			      s_{J\setminus I} \in [0,\varepsilon)^{|J\setminus I|}
			      \right\}.
		      \end{equation}
		      See Figure \eqref{fig: 3simplex}. If $J=I$, then $\widetilde\Delta_{II}^n(\varepsilon) = \Delta_I^n$. The projection to the
		      \emph{collar coordinates} is
		      \begin{equation}
			      \widetilde{\tau}_{IJ}\colon
			      \widetilde{\Delta}_{IJ}^n(\varepsilon)\longrightarrow
			      [0,\varepsilon)^{|J\setminus I|},
			      \qquad
			      \widetilde{\tau}_{IJ}(s_{J})=s_{J\setminus I}.
		      \end{equation}
		\item The \emph{product collar} of $\Delta_I^n$ in $\Delta_J^n$ is
		      \begin{equation}
			      \Delta_{IJ}^n(\varepsilon)
			      := \Delta_I^n \times [0,\varepsilon)^{|J\setminus I|}.
		      \end{equation}
		      whose points we write as \((s_{I},t_{J\setminus I})\). We write $0_{J\setminus I}$ for the zero vector in $[0,\varepsilon)^{|J\setminus I|}$. If $J=I$, then $\Delta_{II}^n(\varepsilon) = \Delta_I^n$.
		      The projection is
		      \begin{equation}
			      \tau_{IJ}\colon
			      \Delta_{IJ}^n(\varepsilon)\longrightarrow
			      [0,\varepsilon)^{|J\setminus I|},
			      \qquad
			      \tau_{IJ}(s_{I},t_{J\setminus I})=t_{J\setminus I}.
		      \end{equation}
	\end{enumerate}

\end{definition}

\begin{remark}
	From now on we fix $n>0$ and $\varepsilon>0$ and drop them from the notation: $\Delta = \Delta^n$, $\Delta_{IJ} = \Delta_{IJ}^n(\varepsilon)$, $\widetilde\Delta_{IJ} = \widetilde\Delta_{IJ}^n(\varepsilon)$, etc.
\end{remark}

\begin{definition}
	\label{def:collarMap}
	Let $\emptyset \neq I \subseteq J \subseteq [n]$. A \emph{collar map for $\Delta_I$ in $\Delta_J$} is a diffeomorphism of manifolds with corners
	\begin{equation}
		\coll_{IJ}:  \Delta_{IJ} \to \widetilde\Delta_{IJ} \subset \Delta_J
	\end{equation}
	such that:
	\begin{enumerate}
		\item \label{itm:collarCoordinates}
		      $ \coll_{IJ} $ preserves the collar coordinates: $ \tilde\tau_{IJ} \circ \coll_{IJ}= \tau_{IJ} $.
		\item  \label{itm:faceIdentity}
		      $ \coll_{IJ} $ is the identity on the face: $\coll_{IJ}(s_{I},0_{J\setminus I})= (s_{I},0_{J\setminus I}) $ for all $ s_{I} \in \Delta_I $.
	\end{enumerate}
\end{definition}

\begin{remark}
	\label{rem:restriction_of_collar_maps}
	If $\coll_{IJ}: \Delta_{IJ}(\varepsilon) \to \widetilde\Delta_{IJ}(\varepsilon)$ is a collar map for $\Delta_I$ in $\Delta_J$, then for any $0<\varepsilon' < \varepsilon$, the restriction  $\coll_{IJ}|_{\Delta_{IJ}(\varepsilon')}: \Delta_{IJ}(\varepsilon') \to \widetilde\Delta_{IJ}(\varepsilon')$  is also a collar map for $\Delta_I$ in $\Delta_J$.
\end{remark}

\begin{definition}
	\label{def:collarEmbeddings1}
	Let $\emptyset \neq I  \subseteq  J \subseteq K \subseteq [n]$.
	\begin{enumerate}
		\item We have canonical embeddings
		      \begin{align}
			      \tilde\iota^{IK}_{JK}
			                                 & : \widetilde\Delta_{IK} \to \widetilde\Delta_{JK},
			                                 &
			      s_{ K}                     & \mapsto s_{K}
			      \\
			      \tilde\iota^{IJ}_{IK}
			                                 & : \widetilde\Delta_{IJ} \to \widetilde\Delta_{IK},
			                                 &
			      s_{ J}                     & \mapsto (s_{J}, 0_{K \setminus J}),
			      \\
			      \iota^{IJ}_{IK}
			                                 & : \Delta_{IJ} \to \Delta_{IK},
			                                 &
			      (s_{ I}, t_{J\setminus I}) & \mapsto (s_{ I},t_{J \setminus I} ,0_{K \setminus J}).
		      \end{align}
		      and projections
		      \begin{align}
			      \widetilde\pi^{IK}_{IJ}
			                                 & : \widetilde\Delta_{IK} \to \widetilde\Delta_{IJ},
			                                 &
			      s_{ K}                     & \mapsto s_{ J},
			      \\
			      \pi^{IK}_{IJ}
			                                 & : \Delta_{IK} \to \Delta_{IJ},
			                                 &
			      (s_{ I}, t_{K\setminus I}) & \mapsto (s_{ I},t_{J \setminus I}).
		      \end{align}
		\item
		      \label{itm:collarEmbeddings1}
		      Given a collar map $\coll_{IJ}$, we have an embedding:
		      \begin{align}
			      \coll^{IK}_{JK}
			                                 & : \Delta_{IK} \to \Delta_{JK},
			                                 &
			      (s_{ I}, t_{K\setminus I}) & \mapsto (\coll_{IJ}(s_{ I},t_{J \setminus I}) ,t_{K \setminus J}).
		      \end{align}
		\item For each map above, a hat denotes the induced map of symplectic manifolds; e.g
		      \begin{equation}
			      \widehat{\coll}_{IJ}: X_{\Delta_{IJ}} \to X_{\widetilde\Delta_{IJ}},
			      \quad \quad
			      \widehat{\tilde\iota}^{\ IK}_{JK}: X_{\widetilde\Delta_{IK}} \to X_{\widetilde\Delta_{JK}}.
		      \end{equation}
	\end{enumerate}
\end{definition}

\begin{definition}
	Let $I \subseteq J $ in $\cc{F}_S$. Denote by $\widetilde W_{IJ}$ the Lagrangian cobordism over $\widetilde\Delta_{IJ}$ (relative to $\vlg$) given by pullback of $W_J$ along the inclusion $\tilde\iota  : \widetilde\Delta_{IJ} \to \Delta_J$.
	Set
	\begin{equation}
		\widetilde L_{IJ} := \tilde\iota^\ast L_J = \widetilde{\Delta}_{IJ} \times_{\Delta_J} L_J,
		\quad \quad
		\widetilde\Psi_{IJ} := \tilde\iota^\ast \Psi_J: \widetilde L_{IJ} \to X_{\widetilde{\Delta}_{IJ}}.
	\end{equation}
	Write $\Psi_{\widetilde\Delta_{IJ}}: \widetilde L_{IJ} \to {\widetilde{\Delta}_{IJ}}$ for the composition of $\widetilde\Psi_{IJ}$ with the projection.  
\end{definition}

Our goal in this subsection is to replace the Lagrangian cobordism $\widetilde W_{IJ}$ over $\widetilde\Delta_{IJ}$ by a Lagrangian cobordism $W_{IJ}$ over the product collar $\Delta_{IJ}$, functorially in $I \subseteq J$.
\begin{definition}
	\label{def:LagrangianCollar}
	Let $I \subseteq J $ in $\cc{F}_S$. Set $L_{IJ} := L_I \times [0,\varepsilon)^{|J \setminus I|}$, with points written $(q_I,t_{J\setminus I})$.
	A \emph{collar map for $ L_I $ in $ L_J $} is a diffeomorphism of manifolds with corners
	\begin{equation}
		\lcol_{IJ}: L_{IJ} \to \widetilde L_{IJ}
	\end{equation}
    such that
    \begin{enumerate}
        \item $\lcol_{IJ}$ preserves the collar coordinates, i.e. the diagram
        \begin{equation}
            \label{eq:collars_over_corner}
		\begin{tikzcd}[column sep=large]
			L_{IJ}
			\arrow[r, "\lcol_{IJ}"]
			\arrow[dr, "\operatorname{pr}_2"']
			&
			\widetilde L_{IJ}
			\arrow[d, "\widetilde\tau_{IJ}
				\circ \Psi_{\widetilde\Delta_{IJ}}"]
			\\
			& {[0,\varepsilon)^{|J\setminus I|}}
		\end{tikzcd}
	   \end{equation}
       commutes, where $\operatorname{pr}_2$ is the projection to the second factor of $L_{IJ}$.
       \item $\lcol_{IJ}$ is the identity on the face: $ \lcol_{IJ}(q_I,0_{J\setminus I}) = \iota(q_I) $, where $ \iota: L_I \to \widetilde{L}_{IJ} $ is the canonical map induced by the face inclusion  $\Delta_I \to \tilde{\Delta}_{IJ} $.
    \end{enumerate}
    A \emph{collar pair for $W_I$ in
    $W_J$} is a pair $(\coll_{IJ},\lcol_{IJ})$ consisting of a collar
    map for $\Delta_I$ in $\Delta_J$ and a collar map for $L_I$ in $L_J$,
    respectively.
\end{definition}

Let $(\coll_{IJ},\lcol_{IJ})$ be a collar pair for $W_I$ in $W_J$.  Denote by $\Psi_{IJ}: L_{IJ} \to X_{\Delta_{IJ}}$  the composition
\begin{equation}
	\label{eq:correspondingRepresentation}
	L_{IJ}
	\xrightarrow{\lcol_{IJ}}
	\widetilde L_{IJ}
	\xrightarrow{\widetilde \Psi_{IJ}}
	X_{\widetilde\Delta_{IJ}}
	\xrightarrow{\widehat{\coll_{IJ}}^{-1}}
	X_{\Delta_{IJ}}.
\end{equation}

\begin{lemma}
	The map $\Psi_{IJ}: L_{IJ} \to X_{\Delta_{IJ}}$ represents a Lagrangian cobordism over $\Delta_{IJ}$ relative to $ \vlg $.
\end{lemma}

\begin{proof}
	Since $\coll_{IJ}$ is a diffeomorphism, the pullback
	$\coll_{IJ}^*\widetilde W_{IJ}$ is a Lagrangian cobordism
	over $\Delta_{IJ}$ relative to $\vlg$.
	By \eqref{eq:correspondingRepresentation}, the immersion
	\[
		\Psi_{IJ}
		= \widehat{\coll_{IJ}}^{-1}
		\circ \widetilde\Psi_{IJ} \circ \lcol_{IJ}
	\]
	represents this pullback.
\end{proof}

\begin{definition}
	Let $W_{IJ}$ be the Lagrangian cobordism represented by $\Psi_{IJ}: L_{IJ} \to X_{\Delta_{IJ}}$. We call $W_{IJ}$ the \emph{representation of $\widetilde W_{IJ}$ corresponding to $(\coll_{IJ}, \lcol_{IJ})$}.
\end{definition}

\begin{definition}
	\label{def:collarEmbeddings2}
	Let $ I \subseteq  J \subseteq K $ in $\cc{F}_S$.
	\begin{enumerate}
		\item We  have a natural inclusion and projection
		      \begin{align}
			      \iota^{IJ}_{IK}
			                              & : L_{IJ} \to L_{IK},
			                              &
			      (q_I, t_{J\setminus I}) & \mapsto (q_I,t_{J \setminus I} ,0_{K \setminus J})
			      \\
			      \pi^{IK}_{IJ}
			                              & : L_{IK} \to L_{IJ},
			                              &
			      (q_I, t_{K\setminus I}) & \mapsto (q_I,t_{J \setminus I}).
		      \end{align}
		\item The inclusion $\tilde\iota^{IK}_{JK}: \widetilde\Delta_{IK} \to \widetilde\Delta_{JK}$ induces a canonical inclusion $\tilde\iota^{IK}_{JK}: \widetilde L_{IK} \to \widetilde L_{JK}$.

		\item
		      \label{itm:collarEmbeddings2}
		      Given a collar pair $(\coll_{IJ}, \lcol_{IJ})$, we have an embedding:
		      \begin{align}
			      \lcol^{IK}_{JK}
			                              & : L_{IK} \to L_{JK},
			                              &
			      (q_I, t_{K\setminus I}) & \mapsto (\lcol_{IJ}(q_I,t_{J \setminus I}) ,t_{K \setminus J}).
		      \end{align}
	\end{enumerate}
\end{definition}

\begin{definition}
	Let $ I \subseteq  J \subseteq K $ in $\cc{F}_S$.  Let $(\coll_{IJ},\lcol_{IJ})$, $(\coll_{IK},\lcol_{IK})$ and $(\coll_{JK},\lcol_{JK})$ be collar pairs for $W_I$ in $W_J$, $W_I$ in $W_K$, and $W_J$ in $W_K$, respectively. We say these collar pairs are \textit{compatible}
	if for all $(s_{I}, t_{K\setminus I}) \in \Delta_{IK}$ and $(q_I, t_{K\setminus I}) \in L_{IK}$ we have
	\begin{equation}
		\label{eq:collars_relation}
		\coll_{IK}(s_{I},t_{K\setminus I})
		=
		\coll_{JK}(\coll_{IJ}(s_{I}, t_{J\setminus I}), t_{K \setminus J}),
		\tag{$R_{IJK}$}
	\end{equation}
	\begin{equation}
		\label{eq:collars_relation_2}
		\lcol_{IK}(q_I,t_{K\setminus I})
		=
		\lcol_{JK}(\lcol_{IJ}(q_I,t_{J\setminus I}), t_{K \setminus J}).
		\tag{$\hat R_{IJK}$}
	\end{equation}
	Equivalently, recalling Definitions \ref{def:collarEmbeddings1}\eqref{itm:collarEmbeddings1} and \ref{def:collarEmbeddings2}\eqref{itm:collarEmbeddings2}, the pairs are compatible if the following diagrams commute:
	\begin{equation}
		\label{eq:compatible_collars}
		\begin{tikzcd}
			{\Delta_{IK}}
			\arrow[d, "{\coll_{JK}^{IK}}"']
			\arrow[r, "{\coll_{IK}}"]
			&
			{\widetilde\Delta_{IK}}
			\arrow[d, "{\tilde\iota_{JK}^{IK}}"]
			\\
			{\Delta_{JK}}
			\arrow[r, "{\coll_{JK}}"]
			&
			{\widetilde\Delta_{JK}}
		\end{tikzcd}
		\quad \quad \quad
		\begin{tikzcd}
			{L_{IK}}
			\arrow[d, "{\lcol_{JK}^{IK}}"']
			\arrow[r, "{\lcol_{IK}}"]
			&
			{\widetilde L_{IK}}
			\arrow[d, "{\tilde\iota_{JK}^{IK}}"]
			\\
			{L_{JK}}
			\arrow[r, "{\lcol_{JK}}"]
			&
			{\tilde  L_{JK}}
		\end{tikzcd}
	\end{equation}
\end{definition}

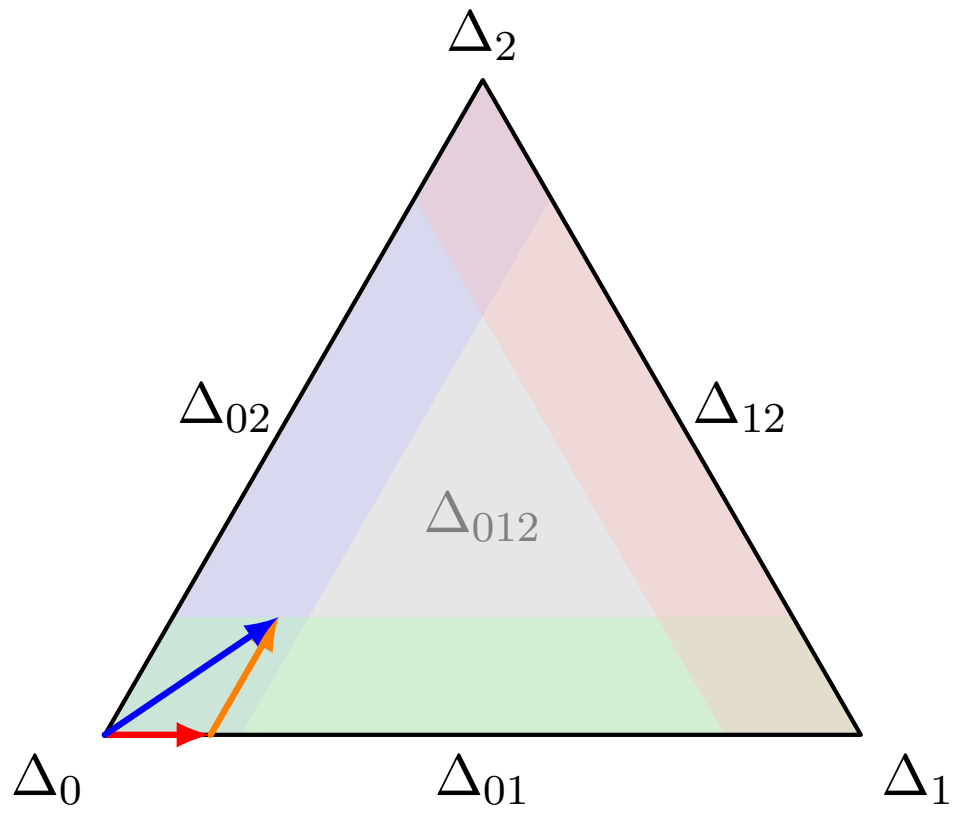
\begin{figure}[htbp]
	\centering
	\begin{tikzpicture}[scale=1.0, line join=round, line cap=round]
		\coordinate (A) at (0,2.887);
		\coordinate (B) at (-2.5,-1.443);
		\coordinate (C) at (2.5,-1.443);

		\def\t{0.18}

		\coordinate (AtoC) at ($(A)!\t!(C)$);
		\coordinate (BtoC) at ($(B)!\t!(C)$);
		\coordinate (BtoA) at ($(B)!\t!(A)$);
		\coordinate (CtoA) at ($(C)!\t!(A)$);
		\coordinate (CtoB) at ($(C)!\t!(B)$);
		\coordinate (AtoB) at ($(A)!\t!(B)$);

		\coordinate (target) at ($ (B)!\t!(A) + (0.7,0)$);

		\fill[gray!20] (A)--(B)--(C)--cycle;

		\fill[blue!25, opacity=0.35] (A)--(B)--(BtoC)--(AtoC)--cycle;  
		\fill[green!30, opacity=0.35] (B)--(C)--(CtoA)--(BtoA)--cycle; 
		\fill[red!25, opacity=0.35] (C)--(A)--(AtoB)--(CtoB)--cycle;   

		\draw[thick] (A)--(B)--(C)--cycle;
		\node[above] at (A) {$\Delta_2$};
		\node[below left] at (B) {$\Delta_0$};
		\node[below right] at (C) {$\Delta_1$};

		\node[left] at ($(A)!0.5!(B)$) {$\Delta_{02}$};
		\node[below] at ($(B)!0.5!(C)$) {$\Delta_{01}$};
		\node[right]  at ($(C)!0.5!(A)$) {$\Delta_{12}$};
		\node[text=gray] at (barycentric cs:A=1,B=1,C=1) {$\Delta_{012}$};

		\draw[-latex, very thick, red] (B) -- (-1.8,-1.443);
		\draw[-latex, very thick, orange] (-1.8,-1.443) -- (target);
		\draw[-latex, very thick, blue] (B) -- (target);

	\end{tikzpicture}
	\caption{Compatibility $R_{0,01,012}$ of collar maps: $ \coll_{0,012}(s_0,t_1,t_2) = \coll_{01,012}(\coll_{0,01}(s_0,t_1),t_2) $.}
	\label{fig: CollarsCompatible}
\end{figure}

\begin{lemma}
	\label{lem:zeroingCoordinates}
	Let $ I \subseteq  J \subseteq K $ in $\cc{F}_S$.  Suppose $(\coll_{IJ},\lcol_{IJ})$, $(\coll_{IK},\lcol_{IK})$ and $(\coll_{JK},\lcol_{JK})$ are compatible collar pairs. Then the following diagrams commute,
	\begin{equation}
		\begin{tikzcd}
			{\Delta_{IJ}}
			\arrow[d, "{\iota_{IK}^{IJ}}"']
			\arrow[r, "{\coll_{IJ}}"]
			&
			{\widetilde\Delta_{IJ}}
			\arrow[d, "{\tilde\iota_{IK}^{IJ}}"]
			\\
			{\Delta_{IK}}
			\arrow[r, "{\coll_{IK}}"]
			&
			{\widetilde\Delta_{IK}}
		\end{tikzcd}
		\quad \quad \quad
		\begin{tikzcd}
			{L_{IJ}}
			\arrow[d, "{\iota_{IK}^{IJ}}"']
			\arrow[r, "{\lcol_{IJ}}"]
			&
			{\widetilde L_{IJ}}
			\arrow[d, "{\tilde\iota_{IK}^{IJ}}"]
			\\
			{L_{IK}}
			\arrow[r, "{\lcol_{IK}}"]
			&
			{\tilde  L_{IK}}
		\end{tikzcd}
	\end{equation}
	Equivalently, for all $(s_{I}, t_{J\setminus I}) \in \Delta_{IJ}$ and $(q_I, t_{J\setminus I}) \in L_{IJ}$ we have
    \begin{equation}
        \label{eq:zeroingCoordinates}
        \begin{aligned}[t]
            \coll_{IJ}(s_I,t_{J\setminus I})
            &=\coll_{IK}(s_I,t_{J\setminus I},0_{K\setminus J}),\\
            \lcol_{IJ}(q_I,t_{J\setminus I})
            &=\lcol_{IK}(q_I,t_{J\setminus I},0_{K\setminus J}).
        \end{aligned}
    \end{equation}
\end{lemma}

\begin{proof}
	Write $t_{K\setminus I} = (t_{J\setminus I},0_{K \setminus J}) $. By \eqref{eq:collars_relation} and the face-identity property \eqref{itm:faceIdentity} in Definition \ref{def:collarMap},
	\begin{equation}
		\coll_{IK}(s_{I},t_{J\setminus I}, 0_{K \setminus J})
		=
		\coll_{JK}(\coll_{IJ}(s_{I}, t_{J\setminus I}), 0_{K \setminus J})
		=
		\coll_{IJ}(s_{I}, t_{J\setminus I}).
	\end{equation}
	The statement for $\lcol_{IJ}$ is analogous, using \eqref{eq:collars_relation_2}.
\end{proof}

\begin{prop}
	\label{prop:compatibleCollarsPullback}
	Let $ I \subseteq  J \subseteq K $ in $\cc{F}_S$.  Let $(\coll_{IJ},\lcol_{IJ})$, $(\coll_{IK},\lcol_{IK})$ and $(\coll_{JK},\lcol_{JK})$  be compatible collar pairs. Then $ (\coll_{JK}^{IK})^{\ast} W_{JK} = W_{IK} $.
	Moreover, $ \Psi_{IK} $ is the unique map making the following diagram commute:
	\begin{equation}
		\label{eq:restriction_of_collar_representation}
		\begin{tikzcd}
			{{L_{IK}} }
			\arrow[d,"\lcol_{JK}^{IK}"']
			\arrow[r, "\Psi_{IK}" ]
			&
			X_{\Delta_{IK}}
			\arrow[d, "\widehat{\coll_{JK}^{IK}}"]
			\\
			L_{JK}
			\arrow[r, "\Psi_{JK}" below]
			&
			X_{\Delta_{JK}}
		\end{tikzcd}
	\end{equation}
\end{prop}

\begin{proof}
	By Corollary \ref{cor:functorialPullback} on the functoriality of pullback, $(\tilde\iota_{JK}^{IK})^{\ast}\widetilde W_{JK} = \widetilde{W}_{IK} $.
	By Lemma \ref{lem:restriction_to_open_subset}, $ \widetilde\Psi_{IK} $ is the unique map making
	the diagram
	\begin{equation}
		\label{eq:cobordism_restricted_collar}
		\begin{tikzcd}
			\widetilde L_{IK}
			\arrow[r, "\widetilde \Psi_{IK}", dashed]
			\arrow[d, "\tilde\iota_{JK}^{IK}" left]
			&
			X_{\widetilde\Delta_{IK}}
			\arrow[d, "\widehat{\tilde\iota_{JK}^{IK}}"]
			\\
			\widetilde L_{JK}
			\arrow[r,  "\widetilde \Psi_{JK}" below]
			&
			X_{\widetilde\Delta_{JK}}
		\end{tikzcd}
	\end{equation}
	commute. Consider the diagram
	
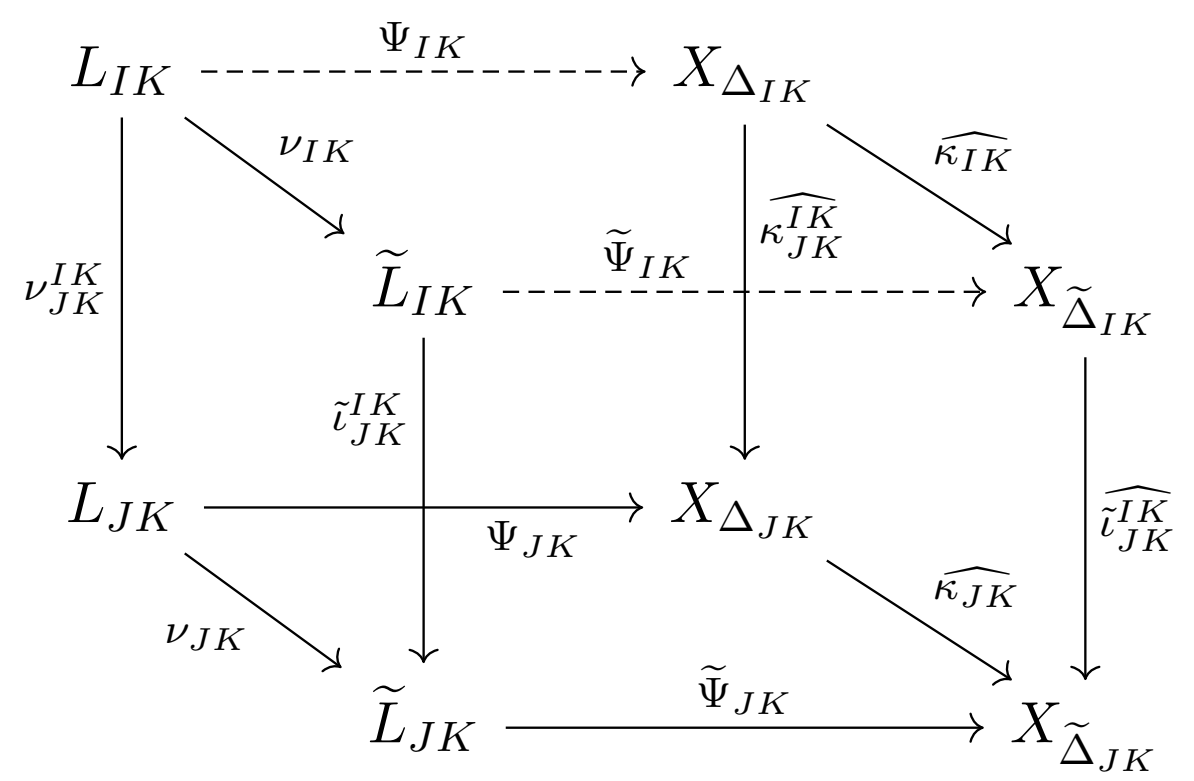
\begin{equation}
		\begin{tikzcd}
			L_{IK}
			\arrow[rd, "\lcol_{IK}"]
			\arrow[rr, "\Psi_{IK}", dashed]
			\arrow[dd, "\lcol_{JK}^{IK}" left]
			&
			&
			X_{\Delta_{IK}}
			\arrow[rd, "\widehat{\coll_{IK}}"]
			\arrow[dd, "\widehat{\coll_{JK}^{IK}}" pos=0.3]
			&
			\\
			&
			\widetilde L_{IK}
			\arrow[dd, "\tilde\iota_{JK}^{IK}"' pos=0.25]
			\arrow[rr, "\widetilde \Psi_{IK}" pos=0.3, dashed]
			&
			&
			X_{\widetilde\Delta_{IK}}
			\arrow[dd, "\widehat{\tilde\iota_{JK}^{IK}}"]
			\\
			L_{JK}
			\arrow[rd, "\lcol_{JK}"']
			\arrow[rr, "\Psi_{JK}"' pos=0.75]
			&
			&
			X_{\Delta_{JK}}
			\arrow[rd, "\widehat{\coll_{JK}}"]
			&
			\\
			&
			\widetilde L_{JK}
			\arrow[rr, "\widetilde \Psi_{JK}"]
			&
			&
			X_{\widetilde\Delta_{JK}}
		\end{tikzcd}
	\end{equation}
	The front face commutes by \eqref{eq:cobordism_restricted_collar}, the left and right faces commute by compatibility of collar maps \eqref{eq:compatible_collars}, and the top and bottom faces commute by definition of $ \Psi_{IK} $, $ \Psi_{JK} $ from the collar maps (equation \eqref{eq:correspondingRepresentation}). Since $\lcol_{IK}, \lcol_{JK}, \widehat{\coll_{IK}}, \widehat{\coll_{JK}} $ are diffeomorphisms, $ \Psi_{IK} $ is the unique map making the back face commute. By Lemma \ref{lem:restriction_to_open_subset}, this uniqueness implies $\Psi_{IK} = (\coll_{JK}^{IK})^\ast \Psi_{JK}$, hence $ (\coll_{JK}^{IK})^{\ast} W_{JK} = W_{IK} $.
\end{proof}

\subsection{Exact collar models}
\label{sec:exact-collar-models}

\begin{ntn}
	\label{ntn:graphs}
	Given $\alpha\in\Omega^1(L)$, set
	\[
		\Gamma(\alpha)\colon L\to T^*L,\qquad p\mapsto (p,\alpha_p),
	\]
	and write $\mathrm{graph}(\alpha):=\mathrm{im}\,\Gamma(\alpha)\subset T^*L$.
	For a function $f\in C^\infty(L)$ we abbreviate
	\[
		\Gamma f := \Gamma(df),\qquad \mathrm{graph}(f):=\mathrm{graph}(df).
	\]
	We also write $\cc{T}(\alpha)(p,\xi)=(p,\xi+\alpha_p)$ for the fiberwise translation and abbreviate
	$\cc{T}f := \cc{T}({df})$.
\end{ntn}

\begin{definition}
	\label{def:WeinsteinModel}
	Let $  I \in \cc{F}_S$. Let
	\begin{equation}
		N^*_I := N^*(\tns{C}{_0^1}(L_I)) \subset T^* L_I
	\end{equation}
	be the conormal bundle of the vertical boundary in $ L_I $. A \textit{Weinstein chart for $W_I$ relative to $\vlg$} consists of an open neighborhood $V_I \subset T^\ast L_I$ of the zero section and a local symplectomorphism $\vp_I : V_I \to X_{\Delta_I}$ such that:
	\begin{enumerate}
		\item \label{itm:zero-section}
		      identifying $L_I$ with the zero section in $T^\ast L_I$, we have $\vp_I|_{L_I} = \Psi_I$;
		\item \label{itm:conormalBundle} for $ q \in \tns{C}{_0^1}(L_I) $ and  $ \xi \in T_q^*L_I $, we have  
        \[
        \vp_I(q,\xi) \in \vlg_{\Delta_I} \iff (q,\xi) \in N^*_I.
        \]
	\end{enumerate}
    More generally, if $V\subset T^*L_I$ is open, a \emph{partial Weinstein chart for $W_I$ relative to $\vlg$} is a local symplectomorphism $\vp:V\to X_{\Delta_I}$ such that $\vp|_{V\cap L_I}=\Psi_I|_{V\cap L_I}$ and \eqref{itm:conormalBundle} holds for every $(q,\xi)\in V$ with $q\in\tns{C}{_0^1}(L_I)$.
\end{definition}

\begin{remark}
	\label{rem:WeinsteinChartSmall}
	Let $\vp_I: V_I \to X_{\Delta_I}$ be a Weinstein chart for $W_I$ relative to $\vlg$. By shrinking $V_I$ if necessary, we may assume that for every $p\in L_I$ there exists an open neighborhood $U\subset L_I$ of $p$ such that $\vp_I|_{\pi^{-1}(U)}$ is a diffeomorphism onto its image, where $\pi:V_I\to L_I$ is the bundle projection.
	Hereafter, we will always make this choice.
\end{remark}

\begin{definition}
	Let $I \in \cc{F}_S$ and $ I \subseteq  J \subseteq [n]$.  Given a relative Weinstein chart $\vp_I: V_I \to X_{\Delta_I}$ for $W_I$, we define an open neighborhood
	\begin{equation}
		V_{IJ} := V_I \times T^\ast \left([0,\varepsilon)^{|J\setminus I|}\right)
		\subset T^\ast L_{IJ},
	\end{equation}
	and a local symplectomorphism
	\begin{equation}
		\vp_{IJ}: V_{IJ} \to X_{\Delta_{IJ}},
		\quad \quad
		\vp_{IJ} := \vp_I \times \text{id}_{T^\ast \left([0,\varepsilon)^{|J \setminus I|}\right)}.
	\end{equation}
\end{definition}

\begin{lemma}
	\label{lem:exact-perturbations}
	Let $I \in \cc{F}_S$ and $ I \subseteq  J \subseteq [n]$. Let $\vp_I: V_I \to X_{\Delta_I}$ be a Weinstein chart for $W_I$ relative to $\vlg$. After possibly shrinking $V_I$, for every smooth function $F: L_{IJ} \to \bb{R}$ such that $\mathrm{graph}(F) \subseteq V_{IJ}$ and $F$ vanishes on the vertical boundary $\tns{C}{_0^1}(L_{IJ})$, the map
	\begin{equation}
		\vp_{IJ} \circ \Gamma F: L_{IJ} \to X_{\Delta_{IJ}}
	\end{equation}
	is a Lagrangian immersion representing a Lagrangian cobordism over $\Delta_{IJ}$.
\end{lemma}

\begin{remark}
	Taking $J=I$ recovers the standard {exact-graph perturbation} statement for a Lagrangian cobordism inside a Weinstein chart.
	If $I\subsetneq J$, then the resulting cobordism $\vp_{IJ}\circ \Gamma F$ can be viewed as the \emph{suspension} of an exact family of Lagrangian cobordisms, parametrized by $[0,\varepsilon)^{|J\setminus I|}$, given by the restrictions of $F$ to the slices $L_I\times\{t_{J\setminus I}\}$.
\end{remark}

\begin{proof}
	Since $\Gamma F$ is a Lagrangian immersion and $\vp_{IJ}$ is a local symplectomorphism, the map $\vp_{IJ} \circ \Gamma F$ is a Lagrangian immersion. We now verify \eqref{itm:bsubmersion}-\eqref{itm:verticalCorners} in Definition \ref{def:LagrangiankCobordism}.
	Recall Definition \ref{def:b_submersions} concerning $b$-submersions. The composition
	\begin{equation}
		\beta: V_{IJ} \xrightarrow{\vp_{IJ}} X_{\Delta_{IJ}} \xrightarrow{\pi_{\Delta_{IJ}}} \Delta_{IJ}
	\end{equation}
	coincides with $\Psi_{\Delta_I} \times \mathrm{id}_{[0,\epsilon)^{|J\setminus I|}}$ on the zero section by Definition \ref{def:WeinsteinModel}\eqref{itm:zero-section}, so it is a b-submersion on the zero section by Definition \ref{def:LagrangiankCobordism}\eqref{itm:bsubmersion}.
	Let $(q,\xi)\in S^l(T^* L_{IJ})$ with $\beta(q,\xi)\in S^k(\Delta_{IJ})$.
	Surjectivity of the induced linear map
	\begin{equation}
		(d\beta|_{(q,\xi)})_\ast :
		T_{(q,\xi)} (T^*L_{IJ})/T_{(q,\xi)}(S^l(T^*L_{IJ})) \to
		T_{\beta(q,\xi)}\Delta_{IJ} / T_{\beta(q,\xi)} S^k(\Delta_{IJ})
	\end{equation}
	is an open condition; by the product formula for $\vp_{IJ}$, we can shrink $V_{I}$ so that $\beta$ is a b-submersion on $V_{IJ}$. Note that
	\begin{equation}
		S^l(T^* L_{IJ}) = \{ (q,\xi) \in T^* L_{IJ} \, : \, q \in S^l(L_{IJ}) \}.
	\end{equation}
	In particular, the fiber directions in $T(T^*L_{IJ})$ are tangent to $S^l(T^* L_{IJ})$, so the map
	\begin{equation}
		(d(\Gamma F)|_q)_\ast : T_qL_{IJ} / T_q S^l(L_{IJ}) \to
		T_{(q,dF_q)} (T^*L_{IJ})/T_{(q,dF_q)}(S^l(T^*L_{IJ}))
	\end{equation}
    is an isomorphism. Set $b=\beta(q,dF_q)$. It follows that
    \begin{equation}
        (d\beta|_{(q,dF_q)})_\ast\circ(d(\Gamma F)|_q)_\ast
        \colon T_qL_{IJ}/T_qS^l(L_{IJ})
        \longrightarrow T_b\Delta_{IJ}/T_bS^k(\Delta_{IJ})
    \end{equation}
	is surjective, hence $\vp_{IJ} \circ \Gamma F$ satisfies Definition \ref{def:LagrangiankCobordism}\eqref{itm:bsubmersion}.

	Denote $N^*_{IJ}:=N^*(\tns{C}{_0^1}(L_{IJ}))$. We have $\tns{C}{_0^1}(L_{IJ}) = \tns{C}{_0^1}(L_{I}) \times [0,\varepsilon)^{|J \setminus I|}$, so
	\begin{equation}
		N^*_{IJ}
		=
		N^*_I \times [0,\varepsilon)^{|J \setminus I|} \times \{0\}.
	\end{equation}
	Let
	\begin{equation}
		(q_I,t_{J\setminus I}) \in   \tns{C}{_0^1}(L_I) \times [0,\varepsilon)^{|J\setminus I|},
		\quad
		(\xi,\zeta) \in T^*_{q_I} L_I \times T^*_{t_{J\setminus I}}[0,\varepsilon)^{|J\setminus I|}.
	\end{equation}
	Since $\vp_{I}$ satisfies Definition \ref{def:WeinsteinModel}\eqref{itm:conormalBundle}, the inclusion
	\begin{equation}
		\vp_{IJ}(q_I,t_{J\setminus I},\xi,\zeta) =
		(\vp_I(q_I,\xi),(t_{J\setminus I},\zeta))\in \vlg_{\Delta_I} \times [0,\varepsilon)^{|J\setminus I|} \times \{0\}= \vlg_{\Delta_{IJ}}
	\end{equation}
	is equivalent to
	\begin{equation}
		(q_I,\xi) \in N^*_{I}, \qquad \zeta=0.
	\end{equation}
	Since $F|_{\tns{C}{_0^1}(L_{IJ})} \equiv 0$, we have
	\begin{equation}
		\left( q_I,d(F|_{L_I \times \{t_{J \setminus I}\}})_{(q_I,t_{J \setminus I})} \right) \in N^*_I,
		\qquad
		d\left(F|_{\{q_I\} \times [0,\varepsilon)^{J \setminus I}}\right)_{(q_I,t_{J \setminus I})} =0.
	\end{equation}
	It follows that  $(\vp_{IJ} \circ \Gamma F)$ is an immersion of $\tns{C}{_0^1}(L_{IJ}) $ in $\vlg_{\Delta_{IJ}}$, proving Definition \ref{def:LagrangiankCobordism}\eqref{itm:verticalBoundary}.

	We verify Definition \ref{def:LagrangiankCobordism}\eqref{itm:verticalTransversality}.
	Let $(q_I,t_{J\setminus I}) \in   \tns{C}{_0^1}(L_I) \times [0,\varepsilon)^{|J\setminus I|}$. Since $W_I$ satisfies Definition \ref{def:LagrangiankCobordism}\eqref{itm:verticalTransversality},
	there exists $v\in T_{q_I}L_I$ such that $d(\Psi_I)_{q_I}(v)\notin T_{\Psi_I(q_I)}\vlg_{\Delta_I}$.
	By Definition \ref{def:LagrangiankCobordism}\eqref{itm:verticalBoundary} we have
	$v\notin T_{q_I}\tns{C}{_0^1}(L_I)$.
	Let \(\pi:T^*L_{IJ}\to L_{IJ}\) be the bundle projection. We have
	\begin{equation}
		d\pi_{((q_I,t_{J \setminus I}),dF_{(q_I,t_{J \setminus I})})}\!\big(d(\Gamma F)_{(q_I,t_{J \setminus I})}(v,0)\big)=(v,0)\notin T_{(q_I,t_{J \setminus I})}\tns{C}{_0^1}(L_{IJ}).
	\end{equation}
	Since a tangent vector to $N_{IJ}^*$ projects under $d\pi$ to a tangent vector to $\tns{C}{_0^1}(L_{IJ})$, we have
    \[
        d(\Gamma F)_{(q_I,t_{J\setminus I})}(v,0) \notin T_{(q_I,t_{J\setminus I})}N_{IJ}^*.
    \]
    Consequently,
    \[
      d(\vp_{IJ}\circ\Gamma F)_{(q_I,t_{J\setminus I})}(v,0)
      \notin T_{(\vp_{IJ}\circ\Gamma F)(q_I,t_{J\setminus I})}
                \vlg_{\Delta_{IJ}},
    \]
	proving Definition \ref{def:LagrangiankCobordism}\eqref{itm:verticalTransversality}.

	Finally, since $\tns{C}{_0^v}(L_I)=\emptyset$ for $v \geq 2$, we have $\tns{C}{_0^v}(L_{IJ})=\tns{C}{_0^v}(L_{I})\times [0,\varepsilon)^{|J\setminus I|}=\emptyset$ for $v \geq 2$, proving Definition \ref{def:LagrangiankCobordism}\eqref{itm:verticalCorners}.
\end{proof}

\begin{definition}
	\label{def:exactCollarModel}
	Let $ I \subseteq  J $ in $\cc{F}_S$.
	An \textit{exact collar model for $W_I$ in $W_J$} is a tuple
	\begin{equation}
		\cc{C}_{IJ} =
		\left(
		(\coll_{IJ},\lcol_{IJ}),
		\vp_I,
		g_{IJ}
		\right)
	\end{equation}
	consisting of:
	\begin{enumerate}
		\item a collar pair $(\coll_{IJ},\lcol_{IJ})$ for $W_I$ in $W_J$;
		\item a relative Weinstein chart $\vp_I: V_I \to X_{\Delta_I}$ for $W_I$;
		\item a smooth function $g_{IJ}: L_{IJ} \to \bb{R}$, the \emph{transition potential},
		      such that:
		      \begin{enumerate}
			      \item \label{itm:boundaryVanishing} $g_{IJ}$ vanishes on $ \tns{C}{_0^1}(L_{IJ}) \simeq  \tns{C}{_0^1}(L_I) \times [0,\varepsilon)^{|J\setminus I|} $ and on the corner $ L_I $;
			      \item \label{itm:suspensionRepresentation} $\Gamma g_{IJ} \subset V_{IJ} $, and the representation of $ \widetilde{W}_{IJ} $ corresponding to $(\coll_{IJ},\lcol_{IJ})$ is generated by $ g_{IJ} $, i.e. the following diagram commutes:
			            \begin{equation}
				            \begin{tikzcd}
					            {L_{IJ}}
					            \arrow[r, "{\Gamma g_{IJ}}"]
					            \arrow[dr, "{\Psi_{IJ}}"']
					            &
					            {V_{IJ}}
					            \arrow[d, "{\vp_{IJ}}"]
					            \\
					            &
					            X_{\Delta_{IJ}}
				            \end{tikzcd}
				            \label{diag:exact-collar-model}
			            \end{equation}
		      \end{enumerate}
	\end{enumerate}
\end{definition}

\begin{remark}
	Let $g_{IJ}: L_{IJ} \to \bb{R}$ be a smooth function satisfying Definition \ref{def:exactCollarModel}\eqref{itm:boundaryVanishing}. By Lemma \ref{lem:exact-perturbations}, $\vp_{IJ} \circ \Gamma g_{IJ}$ represents a Lagrangian cobordism over $\Delta_{IJ}$. Definition \ref{def:exactCollarModel}\eqref{itm:suspensionRepresentation} says this cobordism is exactly $W_{IJ}$.
\end{remark}

\begin{example}
	Let $I=\{0\}$ and $J=\{0,1\}$. Let $\pi: X_{\Delta_{0}} \times T^*[0,\varepsilon)   \to X_{\Delta_{0}} $ be the canonical projection. The Lagrangian immersion $\Psi_{0,01}: L_0 \times [0,\epsilon) \to X_{\Delta_0} \times T^*[0,\epsilon) $ defines a family of Lagrangian embeddings
	\begin{equation}
		t \mapsto (\pi \circ \Psi_{0,01}(-,t)): L_0 \to X_{\Delta_0} \simeq  X.
	\end{equation}
	Equation \eqref{diag:exact-collar-model} means that $g_t:=g_{0,01}(-,t): L_0 \to \bb{R}$ satisfies
	\begin{equation}
		(\pi \circ \Psi_{0,01}(-,t)) = \vp_0 \circ \Gamma g_t
	\end{equation}
	where $\vp_0: V_0 \to X$ is the relative Weinstein chart of $L_0.$
\end{example}

\begin{ntn}
	\label{rem:collarIdentifications}
	Let $I  \subseteq  J \subseteq K$ in $\cc{F}_S$. Moving forward, we will implicitly use the collar maps $\coll_{IJ}$ and $\lcol_{IJ}$ to regard $L_{IJ}$, $\Delta_{IJ}$ as open submanifolds of $L_J$ and $\Delta_J$, respectively. We denote
    \begin{equation}
		V_J|_{L_{IJ}} := V_J \cap T^\ast L_{IJ},
		\quad \quad
		\vp_J|_{L_{IJ}} := \vp_J|_{V_J \cap T^\ast L_{IJ}}.
	\end{equation}
    More generally, we use the
    (collar-dependent) open embeddings $ \coll^{IK}_{JK}  $ and $\lcol^{IK}_{JK}$ of Definitions \ref{def:collarEmbeddings1} and \ref{def:collarEmbeddings2} to regard $L_{IK}$ and $\Delta_{IK}$ as open submanifolds of $L_{JK}$ and $\Delta_{JK}$, respectively, and write
	\begin{equation}
		\label{eq:collarIdentifications}
		V_{JK}|_{L_{IK}} := V_{JK} \cap T^\ast L_{IK},
		\quad \quad
		\vp_{JK}|_{L_{IK}} := \vp_{JK}|_{V_{JK} \cap T^\ast L_{IK}}.
	\end{equation}
\end{ntn}

\begin{figure}
	\begin{tikzpicture}
		\coordinate (A) at (-0.5,1.5);
		\coordinate (B) at (2,2.5);
		\coordinate (C) at (4,2.5);
		\coordinate (D) at (1.5,1.5);
		\coordinate (E) at (-0.5,-1.7);
		\coordinate (F) at (2,-0.5);
		\coordinate (G) at (4,-0.5);
		\coordinate (H) at (1.5,-1.7);

		\draw[very thick] (0,0) -- (5,0);
		\draw[gray, line width=25pt, draw opacity=.15,
			line cap=round, line join=round] (0.2,0) -- (4.8,0);
		\draw[very thick, red] (0,0) -- (3,0);
		\draw[very thick, black] (0,-2) -- (0,2);
		\draw[line width=2pt, blue] (0,-1.5) -- (0,1.7);
		\draw[very thick, black] (-0.5,-0.2) -- (2,0.8);
		\draw[thick, black, smooth] plot coordinates {(1,0.4) (1.8,1) (2,0.5) (3,1.5)};
		\draw[olive, line width=10pt, draw opacity=.25,
			line cap=round, line join=round, smooth] plot coordinates {(1,0.4) (1.8,1) (2,0.5) (3,1.5)};

		\draw[fill=red!30, fill opacity=0.15, draw=red] (A) -- (B) -- (C) -- (D) -- cycle;
		\draw[fill=red!30, fill opacity=0.15, draw=red] (E) -- (F) -- (G) -- (H) -- cycle;
		\draw[fill=red!30, fill opacity=0.15, draw=red] (E) -- (A) -- (D) -- (H) -- cycle;
		\draw[fill=red!30, fill opacity=0.15, draw=red] (H) -- (D) -- (C) -- (G) -- cycle;
		\fill[red!30, opacity=0.15] (A) -- (B) -- (C) -- (D) -- cycle;
		\fill[red!30, opacity=0.15] (E) -- (F) -- (G) -- (H) -- cycle;
		\draw[red] (E) -- (A);
		\draw[red] (H) -- (D);
		\draw[red] (F) -- (B);
		\draw[red] (G) -- (C);

		\node[text=gray] at (4.5,0.25) {$ V_J $};
		\node at (4.5,-0.25) {$ L_J $};
		\node[text=red] at (6.2,2) {$ V_{IJ} = V_I \times T^*[0,\varepsilon)^{|J \setminus I|} $};
		\node [text=red, font=\tiny] at (2.25,-0.25) {$ L_{IJ} $};
		\node [text=blue] at (0.25,0.72) {$V_I$};
		\node[text=olive, font=\tiny] at (3.5,0.9) {$ \cc{T}g_{IJ}(V_{J}|_{L_{IJ}}) $};
		\node [text=blue] at (-0.3,0.2) {$L_I$};
		\node [text=black, font=\tiny] at (3,1.8) {$\mathrm{im}\Gamma g_{IJ}$};
		\fill[blue] (0,0) circle (2pt);
	\end{tikzpicture}
	\caption{Compatibility of relative Weinstein charts.}
	\label{fig:Weinstein-models}
\end{figure}
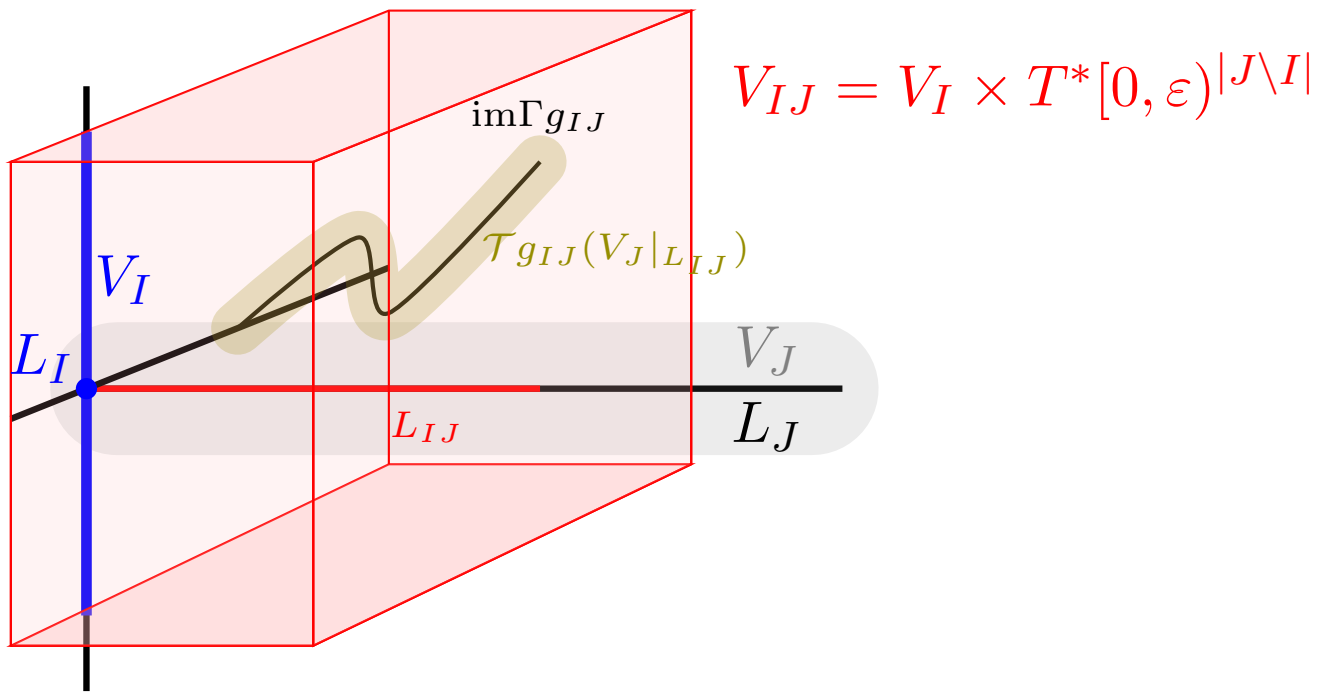

\begin{definition}
	\label{def:compatibleExactModels}
	Let $ I \subseteq  J \subseteq K $ in $\cc{F}_S$. We say the exact collar models
	\begin{gather}
		\cc{C}_{IJ}=\{(\coll_{IJ},\lcol_{IJ}), \vp_I, g_{IJ}\},
		\quad \quad
		\cc{C}_{IK}=\{(\coll_{IK},\lcol_{IK}), \vp_I, g_{IK}\},
		\\
		\cc{C}_{JK}=\{(\coll_{JK},\lcol_{JK}), \vp_J, g_{JK}\}
	\end{gather}
	are \textit{compatible} if:
	\begin{enumerate}
		\item \label{itm:compatible-collars}
		      the collar pairs $(\coll_{IJ},\lcol_{IJ}), (\coll_{IK},\lcol_{IK}), (\coll_{JK},\lcol_{JK})$ are compatible;
		\item the Weinstein charts $ \vp_I $ and $ \vp_J $ satisfy the following transition relation. We have
		      \begin{equation}
			      \label{eq:compatible_Weinstein_symplecto_inclusion}
			      V_J|_{L_{IJ}} \cup \cc{T}g_{IJ} \left(	V_J|_{L_{IJ}}  \right) \subset V_{IJ},
		      \end{equation}
		      and the following diagram commutes (see Figure \ref{fig:Weinstein-models}):
		      \begin{equation}
			      \begin{tikzcd}
				      V_J|_{L_{IJ}}
				      \arrow[r, "{\cc{T}g_{IJ}}"]
				      \arrow[dr, "{\vp_{J}}|_{L_{IJ}}"']
				      &
				      \cc{T}g_{IJ} \left(V_J|_{L_{IJ}}  \right)
				      \arrow[d, "{\vp_{IJ}}"]
				      \\
				      &
				      X_{\Delta_{IJ}}
			      \end{tikzcd}
			      \label{eq:compatible_Weinstein_symplecto}
		      \end{equation}
	\end{enumerate}
\end{definition}

\begin{remark}
	The maps $	\vp_J|_{L_{IJ}}$ and $\vp_{IJ}$ both land in $X_{\Delta_{IJ}}$, but are centered at different Lagrangians. As $\vp_J$ is a Weinstein chart for $W_J$, its restriction to the zero section in $T^*L_J$ parametrizes $W_J$, and over $L_{IJ}$ it parametrizes the collar piece $W_{IJ}\subset X_{\Delta_{IJ}}$. By contrast, $\vp_{IJ} =  \vp_I \times \text{id}_{T^\ast
		[0,\varepsilon)^{|J \setminus I|}}$ is built from the chart for $W_I$, so its restriction to the zero section in $T^*L_{IJ}\simeq T^*(L_I\times [0,\epsilon)^{|J\setminus I|})$ parametrizes the product
	\begin{equation}
		\Psi_I(L_I)\times [0,\epsilon)^{|J\setminus I|} \subset X_{\Delta_I} \times T^*([0,\epsilon)^{|J\setminus I|}) = X_{\Delta_{IJ}}.
	\end{equation}
	The compatibility condition \eqref{eq:compatible_Weinstein_symplecto} says that the fibers of the two charts are compatible, while the translation $\cc{T} g_{IJ}$ identifies the centers. Thus, in the product chart $\vp_{IJ},$ the collar $W_{IJ}$ is represented by $\Gamma g_{IJ}$ rather than by the zero section.
\end{remark}

\begin{remark}
	The compatibility condition of	Definition~\ref{def:compatibleExactModels} does not explicitly refer to $g_{IK}$ or $g_{JK}.$ However, Proposition~\ref{prop:compatible_generating_functions} below shows that in fact compatibility of exact collar models implies a cocycle condition \eqref{eq:compatible_generating_functions} on $g_{IJ},g_{IK}$ and $g_{JK}$. The reason, explained in detail in the proof of Proposition~\ref{prop:compatible_generating_functions}, is that the compatibility of $\vp_I$ and $\vp_J$ extends to compatibility of $\vp_{IK}$ and $\vp_{JK}$, with transition potential $g_{IJ} \circ \pi_{IJ}^{IK}$, by taking the product of \eqref{eq:compatible_Weinstein_symplecto} with $T^*[0,\varepsilon)^{|K\setminus J|}$. Combining this with the transition identities \eqref{diag:exact-collar-model} for $\vp_{IK}$ and $\vp_{JK}$ gives the cocycle relation.
\end{remark}

\begin{prop}
	\label{prop:compatible_generating_functions}
	Let $ I \subseteq  J \subseteq K $ in $\cc{F}_S$  and let $\cc{C}_{IJ},\cc{C}_{IK} , \cc{C}_{JK}$ be compatible exact collar models.  Then,  the transition potentials satisfy the following \emph{cocycle condition}:
	\begin{equation}
		\label{eq:compatible_generating_functions}
		g_{IK}(q_I,t_{K\setminus I})=g_{JK}(\lcol_{IJ}(q_I,t_{J\setminus I}), t_{K \setminus J})+ g_{IJ}(q_{I},t_{J \setminus I}).
		\tag{$G_{IJK}$}
	\end{equation}
	Equivalently, the following diagram commutes:
	\begin{equation}
		\begin{tikzcd}
			L_{IK}
			\arrow[d, "\lcol_{JK}^{IK}" left]
			\arrow[drr, "g_{IK}- g_{IJ} \circ \pi_{IJ}^{IK}"]
			\\
			L_{JK}
			\arrow[rr, "g_{JK}" below]
			&&
			\bb{R}
		\end{tikzcd}
	\end{equation}
\end{prop}

The proof of Proposition \ref{prop:compatible_generating_functions} is based on the following lemmas.

\begin{lemma}
	Let $ I \subseteq  J \subseteq K $ in $\cc{F}_S$  and let $\cc{C}_{IJ},\cc{C}_{IK} , \cc{C}_{JK}$ be compatible exact collar models. Then,
	\begin{equation}
    \label{eq:productNbhdRestriction}
        \begin{aligned}[t]
            V_J|_{L_{IJ}}\times T^\ast[0,\varepsilon)^{|K\setminus J|}
            &=V_{JK}|_{L_{IK}},\\
            \vp_J|_{L_{IJ}}\times\mathrm{id}_{T^\ast[0,\varepsilon)^{|K\setminus J|}}
            &=\vp_{JK}|_{L_{IK}}.
        \end{aligned}
    \end{equation}
\end{lemma}

\begin{proof}
	Recalling Notation \ref{rem:collarIdentifications}, the lemma follows from the natural identification
	\begin{align}
		V_J|_{L_{IJ}} \times T^\ast [0,\varepsilon)^{|K \setminus J|}
		 & =
		(V_J \cap T^\ast L_{IJ}) \times T^\ast [0,\varepsilon)^{|K \setminus J|}
		\\
		 & =
		\left( V_J \times  T^\ast [0,\varepsilon)^{|K \setminus J|}\right)\cap
		\left(
		T^\ast L_{IJ} \times T^\ast [0,\varepsilon)^{|K \setminus J|}
		\right)
		\\
		 & =
		V_{JK} \cap T^\ast L_{IK}
		\\
		 & =
		V_{JK}|_{L_{IK}}.
		\qedhere
	\end{align}
\end{proof}

\begin{lemma}
	Let $ I \subseteq  J \subseteq K $ in $\cc{F}_S$  and let $\cc{C}_{IJ},\cc{C}_{IK} , \cc{C}_{JK}$ be compatible exact collar models. Then,
	\begin{equation}
		\label{eq:extended_compatible_Weinstein_symplecto}
		\vp_{JK}|_{L_{IK}} = \vp_{IK} \circ \cc{T}(g_{IJ} \circ \pi_{IJ}^{IK}).
	\end{equation}
\end{lemma}

\begin{proof}
	Taking the product of  \eqref{eq:compatible_Weinstein_symplecto} with $T^\ast [0,\varepsilon)^{|K \setminus J|}$ and using \eqref{eq:productNbhdRestriction} we obtain
	\begin{align}
		\vp_{JK}|_{L_{IK}}
		 & =
		\vp_J|_{L_{IJ} } \times \text{id}_{T^\ast [0,\varepsilon)^{|K \setminus J|}}
		\\
		 & =
		(\vp_{IJ} \circ \cc{T}g_{IJ}) \times \text{id}_{T^\ast [0,\varepsilon)^{|K \setminus J|}}
		\\
		 & =(\vp_{IJ} \times \text{id}_{T^\ast [0,\varepsilon)^{|K \setminus J|}}) \circ \cc{T}(g_{IJ} \circ \pi_{IJ}^{IK})
		=
		\vp_{IK} \circ \cc{T}(g_{IJ} \circ \pi_{IJ}^{IK}).
		\qedhere
	\end{align}
\end{proof}

\begin{lemma}
	Let $ I \subseteq  J \subseteq K $ in $\cc{F}_S$  and let $\cc{C}_{IJ},\cc{C}_{IK} , \cc{C}_{JK}$ be compatible exact collar models. Then,
	\begin{equation}
		\label{eq:extended_compatible_Weinstein_symplecto2}
		\widehat{\coll_{JK}^{IK}} \circ \vp_{JK}|_{L_{IK}} \circ \Gamma (g_{JK}\circ \lcol_{JK}^{IK})
		=
		\vp_{JK} \circ \Gamma g_{JK} \circ \lcol_{JK}^{IK}.
	\end{equation}
\end{lemma}

\begin{proof}
	Let $\iota: V_{JK}|_{L_{IK}} \to V_{JK}$ be the inclusion determined by $ \lcol_{JK}^{IK} $ (see Notation \ref{rem:collarIdentifications}). It follows by definition that the diagram
	\begin{equation}
		\begin{tikzcd}
			L_{IK}
			\arrow[rr, "{\Gamma (g_{JK}\circ \lcol_{JK}^{IK})}"]
			\arrow[d, "{\lcol_{JK}^{IK}}"']
			&
			&
			V_{JK}|_{L_{IK}}
			\arrow[d, "\iota", hook']
			\arrow[rr, "{\vp_{JK}}|_{L_{IK}}"]
			&
			&
			X_{\Delta_{IK}}
			\arrow[d, "\widehat{\coll_{JK}^{IK}}"]
			\\
			L_{JK}
			\arrow[rr, "{\Gamma g_{JK}}"]
			&
			&
			V_{JK}
			\arrow[rr, "\vp_{JK}"]
			&
			&
			X_{\Delta_{JK}}
		\end{tikzcd}
	\end{equation}
	commutes, which proves the Lemma.
\end{proof}

\begin{proof}[Proof of Proposition \ref{prop:compatible_generating_functions}]
	We have
	\begin{align}
		\widehat{\coll_{JK}^{IK}} \circ \vp_{JK}|_{L_{IK}} \circ \Gamma (g_{JK}\circ \lcol_{JK}^{IK})
		 & =
		\vp_{JK} \circ \Gamma g_{JK} \circ \lcol_{JK}^{IK}
		 &                                             &
		\text{ by } \eqref{eq:extended_compatible_Weinstein_symplecto2}
		\\
		 & = \Psi_{JK} \circ \lcol_{JK}^{IK}
		 &                                             &
		\text{ by }
		\eqref{diag:exact-collar-model}
		\\
		 & = \widehat{\coll_{JK}^{IK}} \circ \Psi_{IK}
		 &                                             &
		\text{ by }
		\eqref{eq:restriction_of_collar_representation}.
	\end{align}
	It follows from the compatibility of the collar pairs in Definition \ref{def:compatibleExactModels}\eqref{itm:compatible-collars} and the uniqueness in Proposition \ref{prop:compatibleCollarsPullback} that
	\begin{equation}
		\Psi_{IK} = \vp_{JK}|_{L_{IK}} \circ \Gamma (g_{JK}\circ\lcol_{JK}^{IK}).
	\end{equation}
	Using \eqref{eq:extended_compatible_Weinstein_symplecto} gives
	\begin{equation}
		\Psi_{IK}= \vp_{IK}  \circ  \cc{T}(g_{IJ} \circ \pi_{IJ}^{IK}) \circ \Gamma (g_{JK}\circ \lcol_{JK}^{IK})
		= \vp_{IK} \circ \Gamma (g_{JK}\circ \lcol_{JK}^{IK} + g_{IJ}\circ \pi_{IJ}^{IK}).
	\end{equation}
	Thus, by \eqref{diag:exact-collar-model},
	\begin{equation}
		\vp_{IK} \circ \Gamma g_{IK}
		=
		\vp_{IK} \circ \Gamma (g_{JK}\circ \lcol_{JK}^{IK} + g_{IJ}\circ \pi_{IJ}^{IK}).
	\end{equation}
	Since $\vp_{IK}$ is injective on fibers by Remark \ref{rem:WeinsteinChartSmall}, we deduce that equation \eqref{eq:compatible_generating_functions} holds on $L_{IK}$ up to an additive constant. This constant must be zero. Indeed, for $ (q_I,t_{K\setminus I}) \in \tns{C}{_0^1}(L_{IK}) $ we have $\pi_{IJ}^{IK}(q_I, t_{K\setminus I}) \in \tns{C}{_0^1}(L_{IJ})$ and $ \lcol_{JK}^{IK}(q_I,t_{K\setminus I}) \in \tns{C}{_0^1}(L_{JK}) $, so each term vanishes by condition \eqref{itm:boundaryVanishing} in Definition \ref{def:exactCollarModel}.
\end{proof}

\begin{remark}
	Let $ I \subseteq  J \subseteq K $ in $\cc{F}_S$. Note that condition \eqref{itm:boundaryVanishing} in Definition \ref{def:exactCollarModel} together with the cocycle condition \eqref{eq:compatible_generating_functions} gives the zeroing identity
    \begin{equation}
        \label{eq:zeroingCoordinatesGenFunc}
        \begin{aligned}[t]
            g_{IK}(q_I,t_{J\setminus I},0_{K\setminus J})
            &=g_{JK}\bigl(\lcol_{IJ}(q_I,t_{J\setminus I}),0_{K\setminus J}\bigr)
            +g_{IJ}(q_I,t_{J\setminus I})\\
            &=g_{IJ}(q_I,t_{J\setminus I}).
        \end{aligned}
    \end{equation}
\end{remark}

\subsection{Compatible collar systems}
\label{sec:compatible-systems}

The following theorem is the main result of this section. Recall the family $\cc{F}_S$ determined by a Lagrangian $W$ over a simplicial set $S$ (Definition \ref{def:family-of-cobordisms}).

\begin{theorem}
	\label{thm:compatible_exact_models}
	For $\varepsilon>0$ sufficiently small there exist:

	\begin{enumerate}
		\item for every $J \in \cc{F}_S$ a relative Weinstein chart
		      $\vp_J:V_J\to X_{\Delta_J}$ for $W_J$; and

		\item for every inclusion $J\subseteq K$ in $\cc{F}_S$ a collar pair
		      $(\coll_{JK},\lcol_{JK})$ and a  function
		      $g_{JK}:L_{JK}\to\mathbb R$,
	\end{enumerate}
	such that each triple $\cc{C}_{JK} = ((\coll_{JK},\lcol_{JK}),\vp_J,g_{JK})$ is an exact
	collar model of $W_J$ in $W_K$, and for every chain
	$I\subsetneq J\subseteq K$ in $\cc{F}_S $ the collars $\cc{C}_{IJ}, \cc{C}_{IK}, \cc{C}_{JK}$ are compatible.
	Moreover, the collars $\{\coll_{JK} \mid J \subsetneq K,\; J,K \in \cc{F}_S\}$ extend to a compatible system $\{\coll_{JK} \mid J \subsetneq K \subseteq [n]\}$.
\end{theorem}

\begin{proof}
	The theorem follows from Propositions \ref{prop:extension_of_collar_pairs}, \ref{prop:inductive_construction_of_Weinstein} and \ref{prop:existence_of_generating_function} by induction on $|J|$.
\end{proof}
Fix $1\leq m \leq n+1$. If $m>1$, suppose that for some $\varepsilon>0$ we have constructed a family of compatible exact collar models
\begin{equation}
	\cc{C}_{m-1} = \{ \cc{C}_{JK} = ((\coll_{JK},\lcol_{JK}),
	\vp_J, g_{JK})
	\,\, | \,\,
	J \subsetneq K \text{ in } \cc{F}_S, \;\; |J| \leq m-1
	\}
\end{equation}
and an ambient system of compatible simplex collars
\begin{equation}
	\cc{C}_{m-1}(\Delta)=\left\{ \coll_{JK}: \Delta_{JK} \to \widetilde{\Delta}_{JK} \; \big | J \subsetneq K \subseteq [n], \;\; |J| \leq m-1 \right\}
\end{equation}
which extends the collars from $\cc{C}_{m-1}$ whenever $K \in \cc{F}_S$. The rest of this section proves the following propositions.

\begin{prop}
	\label{prop:extension_of_collar_pairs}
	Fix $J \subseteq [n]$ with $|J|=m$. Then, after possibly shrinking $\varepsilon$,
	\begin{enumerate}
		\item for every $K \subseteq [n]$ with $J \subseteq K$ there exists a collar map $\coll_{JK}: \Delta_{JK} \to \widetilde{\Delta}_{JK}$ compatible with the system $\cc{C}_{m-1}(\Delta)$.
		\item for every $K \in \cc{F}_S$ with $J \subseteq K$ there exists a Lagrangian collar map $\lcol_{JK}: L_{JK} \to \widetilde{L}_{JK}$, such that $(\coll_{JK}, \lcol_{JK})$ is a collar pair compatible with all collar pairs in $\cc{C}_{m-1}$.
	\end{enumerate}
\end{prop}

\begin{prop}
	\label{prop:inductive_construction_of_Weinstein}
	Fix $J \in \cc{F}_S$ with $|J|=m$. After possibly shrinking $\varepsilon$,  there exists a relative Weinstein chart $\vp_J: V_J \to X_{\Delta_J}$ for $W_J$ such that, for every $I \subsetneq J$, the chart $\vp_I$ from $\cc{C}_{m-1}$ is compatible with $\vp_J$: it satisfies \eqref{eq:compatible_Weinstein_symplecto_inclusion} and \eqref{eq:compatible_Weinstein_symplecto}.
\end{prop}

\begin{prop}
	\label{prop:existence_of_generating_function}
	Fix $J \in \cc{F}_S$ with $|J|=m$. Let $K \in \cc{F}_S$ with $J \subsetneq K$, let $\coll_{JK}$ as in Proposition \ref{prop:extension_of_collar_pairs} and $\vp_J: V_J \to X_{\Delta_J}$ as in Proposition \ref{prop:inductive_construction_of_Weinstein}.  Then, after possibly shrinking $\varepsilon$,  we can choose $\lcol_{JK}$ so that there exists a smooth function $ g_{JK}: L_{JK} \to \bb{R} $ with
	\begin{equation}
		\cc{C}_{JK} = ((\coll_{JK},\lcol_{JK}), \vp_J, g_{JK})
	\end{equation}
	an exact collar model for $W_J$ in $W_K$.
\end{prop}

\begin{remark}
	Since there are only finitely many pairs $ \emptyset \neq I \subseteq J \subseteq [n] $, the constructions require finitely many restrictions of $ \varepsilon $. After all choices are made, we fix a single $\varepsilon$ which is the minimum of those arising along the way; all collar maps and Weinstein charts are then understood to be restricted to this fixed $\varepsilon$.
\end{remark}

\subsection{Extension of collars}
\label{sec:extend-collars}

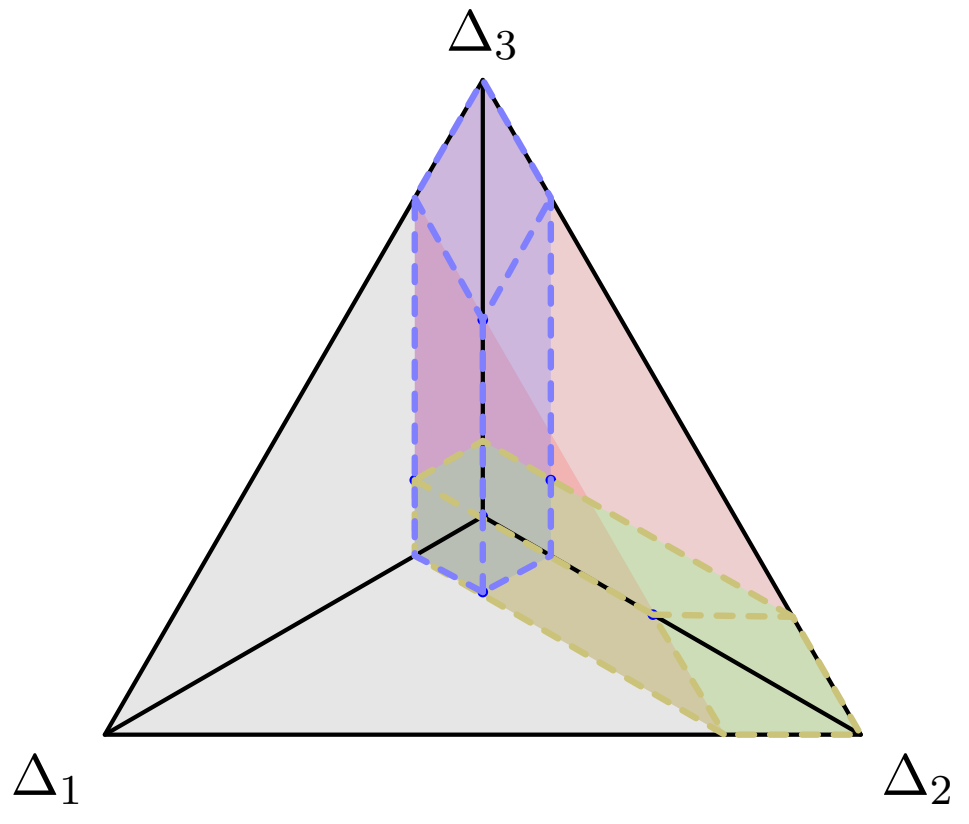
\begin{figure}[htbp]
	\centering
	\begin{tikzpicture}[scale=1.0, line join=round, line cap=round]
		\coordinate (A) at (0,2.887);
		\coordinate (B) at (-2.5,-1.443);
		\coordinate (C) at (2.5,-1.443);

		\def\t{0.18}

		\coordinate (AtoC) at ($(A)!\t!(C)$);
		\coordinate (BtoC) at ($(B)!\t!(C)$);
		\coordinate (BtoA) at ($(B)!\t!(A)$);
		\coordinate (CtoA) at ($(C)!\t!(A)$);
		\coordinate (CtoB) at ($(C)!\t!(B)$);
		\coordinate (AtoB) at ($(A)!\t!(B)$);

		\coordinate (center) at (barycentric cs:A=1,B=1,C=1);
		\coordinate (CentoB) at ($(center)!\t!(B)$);
		\coordinate (CentoC) at ($(center)!\t!(C)$);
		\coordinate (CentoB1) at ($(CentoB)+(0,0.5)$);
		\coordinate (CentoC1) at ($(CentoC)+(0,0.5)$);
		\coordinate (CentoC2) at ($(C)!0.55!(center)$);
		\coordinate (CentoA2) at ($(A)!0.55!(center)$);
		\coordinate (CentoA) at ($(center)+(0,0.5)$);
		\fill[gray!20] (A)--(B)--(C)--cycle;

		\fill[red!35, opacity=0.35] (C)--(A)--(AtoB)--(CtoB)--cycle;   
		\fill[red!50, opacity=0.5] (AtoB)--(CentoB)--(CtoB)--cycle;
		\fill[blue!50, opacity=0.3] (AtoB)--(CentoB)--(0,-0.5)--(CentoC)--(AtoC)--(A)--cycle;
		\fill[green!50, opacity=0.3] (CtoB)--(CentoB)--(CentoB1)--(CentoA)--(CtoA)--(C)--cycle;

		\draw[thick] (A)--(B)--(C)--cycle;
		\node[above] at (A) {$\Delta_3$};
		\node[below left] at (B) {$\Delta_1$};
		\node[below right] at (C) {$\Delta_2$};

		\fill[blue] (center) circle (1pt);
		\draw[thick] (A) -- (center);
		\draw[thick] (B) -- (center);
		\draw[thick] (C) -- (center);
		\fill[blue] (0,-0.5) circle (1pt);
		\fill[blue] (CentoB1) circle (1pt);
		\fill[blue] (CentoC1) circle (1pt);
		\fill[blue] (CentoC2) circle (1pt);
		\fill[blue] (CentoA2) circle (1pt);

		\draw[dashed, very thick, olive!50] (CentoB1)--(CentoC2)--(CtoB)--(CentoB)--cycle;
		\draw[dashed, very thick, olive!50] (CentoC2)--(CtoA)--(C)--(CtoB);
		\draw[dashed, very thick, olive!50] (CtoA)--(CentoA)--(CentoB1);

		\draw[dashed, very thick, blue!50] (0,-0.5)--(CentoB)--(AtoB)--(A)--(AtoC)--(CentoC)--cycle;
		\draw[dashed, very thick, blue!50] (AtoB) -- (CentoA2) -- (AtoC);
		\draw[dashed, very thick, blue!50] (0,-0.5) -- (CentoA2);
	\end{tikzpicture}
	\caption{The collar $\widetilde\Delta_{023,\top}$ (in red), $\widetilde\Delta_{02,\top}$ (in green) and $\widetilde\Delta_{03,\top}$ (in blue), in the 3-simplex with $\top = \{0123\}$. }
	\label{fig: collarExtension}
\end{figure}

In this subsection we prove Proposition~\ref{prop:extension_of_collar_pairs}.
Let $K \in \cc{F}_S$ and $J \subseteq K$. We construct the collars $(\coll_{JK},\lcol_{JK})$ as time-$t$ flows of vector fields on $\widetilde\Delta_{J K}$ and $\widetilde L_{J K}$ whose pushforwards along
\[
	\widetilde L_{JK}\xrightarrow{\ \Psi_{\widetilde\Delta_{JK}}\ }\widetilde\Delta_{JK}
	\xrightarrow{\ \tilde\tau_{JK}\ }[0,\varepsilon)^{|K \setminus J|}
\]
are the coordinate fields $\partial/\partial t_k$.
To ensure compatibility with previously constructed collars, we interpolate these lifted fields with the vector fields prescribed by these collars on codimension-$1$ faces $\widetilde\Delta_{IK},\widetilde L_{IK}$ (with $I=J\setminus\{j\}$); by prior compatibility they agree on overlaps, hence the flows glue to a single compatible collar pair. See Figure~\ref{fig: collarExtension} for an example of overlapping collars with $K=\{0,1,2,3\}$, $J=\{0,2,3\}$, and
$I=\{0,2\}$ or $\{0,3\}$.

\begin{lemma}
	\label{lem:zetaFields}
	Write $\top = [n]$. Let $J \subseteq [n]$ with $|J|=m$. If $J \in \cc{F}_S$, let $F \in \cc{F}_S$ be a facet such that $J \subseteq F$, i.e. there is no $G \in \cc{F}_S$ with $F \subsetneq G$.
	Define open subsets by
	\begin{equation}
		U_0 = \bigcup_{I =J\setminus \{j\}} \widetilde \Delta_{I\top} \subset \widetilde \Delta_{J\top},
		\quad \quad
		\widehat U_0 = \bigcup_{I =J\setminus \{j\}} \widetilde L_{I F} \subset \widetilde L_{JF}.
	\end{equation}
	\begin{enumerate}
		\item For each $k \in \top \setminus J$ there exists a vector field $\zeta_k$ on $U_0$ such that for every $ I=J\setminus\{j\} $
		      \begin{align}
			      (\coll_{I\top})_\ast \left( \frac{\partial}{\partial t_k} \right)
			       & = \zeta_k |_{\widetilde\Delta_{I \top}} ,
                   \label{eq:pushVector1}
			      \\
			      \label{eq:fields_near_corners1}
			      (\tilde\tau_{I\top})_\ast (\zeta_k)
			       & = \frac{\partial}{\partial t_k}.
		      \end{align}

		\item For each $k \in F\setminus J$ there exists a vector field $\hat\zeta_k$ on $\widehat U_0$ such that for every $ I = J\setminus \{j\} $
		      \begin{align}
			      \label{eq:pushVector2}
			      (\lcol_{I F})_\ast \left( \frac{\partial}{\partial t_k} \right)
			       & = \hat\zeta_k |_{\widetilde L_{I F}},
			      \\
			      \label{eq:fields_near_corners2}
			      d\bigl(\tilde\tau_{IF}\circ\Psi_{\widetilde\Delta_{IF}}\bigr)
			      \bigl(\hat\zeta_k|_{\widetilde L_{IF}}\bigr)
			       & = \frac{\partial}{\partial t_k}.
		      \end{align}
	\end{enumerate}
\end{lemma}

\begin{proof}
	Define $\zeta_k$ on each $ \widetilde\Delta_{I \top} $ by \eqref{eq:pushVector1} and $ \hat \zeta_k $  on each $ \widetilde L_{I F} $ by \eqref{eq:pushVector2}.  We first show that the definitions of $\zeta_k$ agree on overlaps. The argument for $\hat\zeta_k$ is analogous, so it is omitted.

	Let $j_1,j_2 \in J$ and denote $I_1 = J \setminus \{j_1\}$, $I_2 = J \setminus \{j_2\}$ and $H = J \setminus \{j_1,j_2\}$. For $s_{\top} \in \widetilde\Delta_{I_1\top} \cap \widetilde\Delta_{I_2\top} = \widetilde\Delta_{H\top} $ let
    \begin{equation}
		(s_{I_1},t_{j_1},t_{\top\setminus J})\in\Delta_{I_1\top},
		\quad
		(s_{I_2},t_{j_2},t_{\top\setminus J})\in\Delta_{I_2\top},
		\quad
		(s_H,t_{j_1},t_{j_2},t_{\top\setminus J})\in\Delta_{H\top},
	\end{equation}
	with
	\begin{equation}
		\coll_{I_1\top}(s_{I_1}, t_{j_1},t_{\top \setminus J})
		=
		s_{\top}
		=
		\coll_{I_2\top}(s_{I_2}, t_{j_2},t_{\top \setminus J}).
	\end{equation}
	and
	\begin{equation}
		\coll_{HI_1}(s_{H}, t_{j_2}) = s_{I_1},
		\quad \quad
		\coll_{HI_2}(s_{H}, t_{j_1}) = s_{I_2}.
	\end{equation}
	By assumption, $ \coll_{I_1\top}, \coll_{I_2\top} $ satisfy  $\hyperref[eq:collars_relation]{(R_{H I_1 \top})}$ and $\hyperref[eq:collars_relation]{(R_{H I_2 \top})}$ because $ |I_1|=|I_2|=m-1 $. Therefore,
	\begin{align}
		\coll_{I_1\top}(s_{I_1}, t_{j_1},t_{\top\setminus J})
		 & =
		\coll_{I_1\top}(\coll_{HI_1}(s_{H}, t_{j_2}), t_{j_1},t_{\top\setminus J})
		\\
		 & =
		\coll_{H\top}(s_{H}, t_{j_1}, t_{j_2}, t_{\top\setminus J})
		\\
		 & =
		\coll_{I_2\top}(\coll_{HI_2}(s_{H}, t_{j_1}), t_{j_2},t_{\top\setminus J})
		\\
		 & =
		\coll_{I_2\top}(s_{I_2}, t_{j_2},t_{\top\setminus J})
	\end{align}
	for all $t_{\top\setminus J}$. In particular,
	\begin{equation}
		\frac{\partial }{\partial t_k}\coll_{I_1\top}(s_{I_1}, t_{j_1},t_{\top\setminus J}) = \frac{\partial}{\partial t_k} \coll_{I_2\top}(s_{I_2}, t_{j_2},t_{\top\setminus J}),
	\end{equation}
	so the two local definitions of $\zeta_k$ agree on the overlap $ \widetilde{\Delta}_{H \top} $.

	Equation \eqref{eq:fields_near_corners1} follows from  \eqref{itm:collarCoordinates} in Definition \ref{def:collarMap}.
	By \eqref{eq:collars_over_corner}, we have
	\[
		(\tilde\tau_{IF}\circ\Psi_{\widetilde\Delta_{IF}}
		\circ\lcol_{IF})(q_I,t_{F\setminus I}) = t_{F\setminus I}.
	\]
	Differentiating with respect to $t_k$ and using
	\eqref{eq:pushVector2} gives \eqref{eq:fields_near_corners2}.
\end{proof}

\begin{lemma}
	\label{lem:simplex_field_extension}
	Let $ J,F $ be as in Lemma \ref{lem:zetaFields}. 	Choose $0< \varepsilon_J < \varepsilon$, and open subsets $ U_1 \subset \widetilde\Delta_{J\top}$ and $\widehat U_1 \subset \widetilde L_{JF}$ such that
	\begin{equation}
		U_0 \cup U_1 = \widetilde\Delta_{J \top},
		\quad \quad
		\widehat U_0 \cup \widehat U_1 = \widetilde L_{JF} ,
	\end{equation}
	\begin{equation}
		\bigcup_{I=J\setminus \{j\}}
		\left(
		\widetilde\Delta_{I\top}
		\left( \varepsilon_J \right)
		\right)
		\subseteq
		U_0 \setminus U_1,
		\quad \quad \quad
		\bigcup_{I = J \setminus \{j\}}
		\left(
		\widetilde L_{IF}
		\left( \varepsilon_J \right)
		\right)
		\subseteq
		\widehat U_0 \setminus \widehat U_1.
	\end{equation}
	\begin{enumerate}
		\item For each $k \in \top\setminus J$ there exists a vector field $Z_k$ on $\widetilde\Delta_{J\top}$ such that
		      \begin{equation}
			      \label{eq:simplex_field_extension1}
			      \left(\tilde\tau_{J\top}\right)_\ast Z_k = \frac{\partial}{\partial t_k},
			      \quad \quad
			      {Z}_k|_{U_0 \setminus  U_1} = \zeta_k.
		      \end{equation}
		\item For each $ k \in F\setminus J $ there exists a vector field $\widehat{Z}_k$ on $\widetilde L_{J F}$ such that
		      \begin{equation}
			      \label{eq:simplex_field_extension2}
			      d\bigl(\tilde\tau_{JF}\circ\Psi_{\widetilde\Delta_{JF}}\bigr)
			      (\widehat Z_k)
			      = \frac{\partial}{\partial t_k},
			      \quad\quad
			      \widehat{Z}_k|_{\widehat U_0 \setminus \widehat U_1} = \widehat\zeta_k.
		      \end{equation}
	\end{enumerate}
\end{lemma}

\begin{proof}
	The composition $\widetilde L_{JF} \xrightarrow{{\Psi}_{\widetilde\Delta_{JF}}} \widetilde \Delta_{JF} \xrightarrow{\tilde \tau_{JF}} [0,\varepsilon)^{|F\setminus J|} $ is b-submersive map, since $\tilde \tau_{JF}$ is a submersion and ${\Psi}_{\widetilde\Delta_{JF}}$ is a b-submersion. After possibly shrinking $\varepsilon$, Lemma~\ref{lem:b-submersive_to_submersion} implies this composition is a submersion near $ L_J $. By Proposition \ref{prop:submersion_local_product}, for every $q \in \widehat{U}_1$ and $k \in F\setminus J$ there exist a neighborhood $\widehat{U}_q$ of $q$ and a vector field $\widehat\zeta_{k,q}$ such that
	\begin{equation}
		\left(\tilde\tau_{JF} \circ {\Psi}_{\widetilde\Delta_{JF}}\right)_\ast \widehat\zeta_{k,q} = \frac{\partial}{\partial t_k}.
	\end{equation}
	Using a partition of unity for the cover $\{\widehat{U}_q\}_{q\in \widehat U_1} \cup \widehat{U}_0$ of $\widetilde L_{J F}$, define a global vector field $\widehat{Z}_k$ on $\widetilde L_{J F}$ such that
	\begin{equation}
		\left(\tilde\tau_{JF} \circ {\Psi}_{\widetilde\Delta_{JF}}\right)_\ast \widehat Z_{k} = \frac{\partial}{\partial t_k},
		\quad \quad
		\widehat{Z}_k|_{\widehat U_0 \setminus \widehat U_1} = \widehat\zeta_k.
	\end{equation}
	For each $k\in\top\setminus J$, apply the same construction on
	$\widetilde\Delta_{J\top}$, using local product lifts of
	$\partial/\partial t_k$ under $\tilde\tau_{J\top}$ and the
	prescribed field $\zeta_k$ on $U_0$. This gives a vector field
	$Z_k$ satisfying \eqref{eq:simplex_field_extension1}.
\end{proof}

\begin{proof}
	[Proof of Proposition \ref{prop:extension_of_collar_pairs}]
	Write $\top\setminus J= \{k_0,\dots,k_{n-m}\}$. If $J \in \cc{F}_S$, let $F$ be a facet as in Lemma \ref{lem:zetaFields}, and suppose that $ F\setminus J =\{k_0,\dots,k_s\} $.
	Let $\left\{\widehat{Z}_k\right\}_{k \in F\setminus J}$, $\left\{Z_k\right\}_{k \in \top\setminus J}$ and $ \varepsilon_J $ be as in Lemma \ref{lem:simplex_field_extension}. Shrinking $ \varepsilon_J $ if necessary, the flows
	\begin{equation}
		\left\{\Phi_k^t: \widetilde \Delta_{J\top} \to \widetilde \Delta_{J\top}\right\}_{k \in \top\setminus J}
		\quad \quad
		\left\{\widehat\Phi_k^t: \widetilde L_{JF} \to \widetilde L_{J F} \right\}_{k \in F\setminus J}
	\end{equation}
	of $Z_k$ and $\widehat Z_k$ are defined for all $ t \in [0,\varepsilon_J) $, and remain inside $ \widetilde \Delta_{J\top} $ and $ \widetilde L_{JF}  $, respectively. Define
	\begin{align}
		\lcol_{JF}
		 & : L_{JF} \to \widetilde L_{JF},
		 &
		\lcol_{JF}(q_J, t_{F\setminus J})
		 & =  \widehat \Phi_{k_s}^{t_{k_s}} \circ  \dots \circ \widehat \Phi_{k_0}^{t_{k_0}} (q_J)
		\\
		\coll_{J\top}
		 & : \Delta_{J\top} \to \widetilde \Delta_{J\top}
		 &
		\coll_{J\top}(s_{J}, t_{\top\setminus J})
		 & =  \Phi_{k_{n-m}}^{t_{k_{n-m}}} \circ  \dots \circ \Phi_{k_0}^{t_{k_0}} (s_{J}).
	\end{align}
	Note that the vector fields $Z_k$, $\widehat Z_k$ do not necessarily commute. Nonetheless, the maps $\coll_{J\top}$, $\lcol_{JF}$ are diffeomorphisms whose inverse is given by flowing in the opposite order with negative time.

	We verify that $\coll_{J\top} $ and $  \lcol_{JF}$ are collar maps. First, by \eqref{eq:simplex_field_extension1} the function $ \tilde\tau_{J\top} \circ \coll_{J\top} $ satisfies the ordinary differential equation
	\begin{equation}
		\frac{\partial}{\partial t_k} (\tilde\tau_{J\top} \circ \coll_{J\top}) = \frac{\partial}{\partial t_k},
		\quad \quad
		(\tilde\tau_{J\top} \circ \coll_{J\top})|_{t=0} =0,
	\end{equation}
	hence $ \tilde\tau_{J\top} \circ \coll_{J\top}= \tau_{J\top} $, proving \eqref{itm:collarCoordinates} in Definition \ref{def:collarMap}. Second, $ \coll_{J\top}(s_{J},0)= s_{J} $ and $ \lcol_{JF}(q_J,0) =q_J $  since the flows at time 0 are the identity, proving \eqref{itm:faceIdentity} in Definition \ref{def:collarMap} and Definition \ref{def:LagrangianCollar}.

	For $ J \subseteq K \subseteq F $, define the restricted collars
	\begin{equation}
		\label{eq:restrictedCollar1}
		\lcol_{JK}:= \lcol_{JF} \circ \iota^{JK}_{JF}: L_{JK} \to \widetilde L_{JK},
		\quad \quad
		\lcol_{JK}(q_{J}, t_{K \setminus J})
		=
		\lcol_{JF}(q_{J}, t_{K \setminus J},0_{F\setminus K})
	\end{equation}
	similarly, for $ K \subseteq \top $ such that $ J \subseteq K $, define
    \begin{equation}
		\label{eq:restrictedCollar2}
		\begin{aligned}[t]
			\coll_{JK}
			&:=\coll_{J\top}\circ\iota^{JK}_{J\top}
				\colon\Delta_{JK}\to\widetilde\Delta_{JK},\\
			\coll_{JK}(s_J,t_{K\setminus J})
			&=\coll_{J\top}(s_J,t_{K\setminus J},0_{\top\setminus K}).
		\end{aligned}
	\end{equation}
	We claim $ (\coll_{JF},\lcol_{JF}) $ is a collar pair.
    By \eqref{eq:simplex_field_extension2}, for $k\in F\setminus J$ and $q \in \widetilde{L}_{JF}$, we have
\[
    (\tilde\tau_{JF}\circ\Psi_{\widetilde\Delta_{JF}})
    (\widehat\Phi_k^t(q))
    =
    (\tilde\tau_{JF}\circ\Psi_{\widetilde\Delta_{JF}})(q)+t e_k.
    \]
    Hence,
    \begin{align}
        (\tilde\tau_{JF}\circ\Psi_{\widetilde\Delta_{JF}}\circ\lcol_{JF})
        (q_J,t_{F\setminus J})
        &=
        (\tilde\tau_{JF}\circ\Psi_{\widetilde\Delta_{JF}})(q_J)
        +\sum_{k\in F\setminus J}t_k e_k =t_{F\setminus J},
    \end{align}
    which is \eqref{eq:collars_over_corner}.
    Since $ (\coll_{JF}, \lcol_{JF}) $ is a collar pair, so is $(\coll_{JK},\lcol_{JK})$ for every $ J \subseteq K \subseteq F $.

    It remains to verify the compatibility relation $\hyperref[eq:collars_relation]{(R_{IJK})}$.   Let $j \in J$ and set $I= J \setminus \{j\}$. We first prove $\hyperref[eq:collars_relation]{(R_{IJ\top})}$.
    By \eqref{eq:pushVector1} and \eqref{eq:simplex_field_extension1},
    for $k\in\top\setminus J$ we have
    \[
        Z_k=(\coll_{I\top})_\ast\left(\frac{\partial}{\partial t_k}\right)
        \qquad\text{on }\widetilde\Delta_{I\top}(\varepsilon_J).
    \]
    Thus
    \begin{align}
        \coll_{I\top}(s_I,t_j,t_{\top\setminus J})
        &=
        \Phi_{k_{n-m}}^{t_{k_{n-m}}}\circ\dots\circ
        \Phi_{k_0}^{t_{k_0}}
        \bigl(\coll_{I\top}(s_I,t_j,0_{\top\setminus J})\bigr)\\
        &=
        \Phi_{k_{n-m}}^{t_{k_{n-m}}}\circ\dots\circ
        \Phi_{k_0}^{t_{k_0}}\bigl(\coll_{IJ}(s_I,t_j)\bigr)\\
        &=
        \coll_{J\top}\bigl(\coll_{IJ}(s_I,t_j),t_{\top\setminus J}\bigr),
    \end{align}
    where the second equality follows from \eqref{eq:zeroingCoordinates}.
	Now let $ H \subsetneq J \subseteq K $, fix $j \in J \setminus H$, and set $I=J \setminus \{j\}$. Then
	\begin{align}
		\MoveEqLeft[1] \coll_{JK}(\coll_{HJ}(s_{H}, t_{J\setminus H}), t_{K \setminus J}) \notag\\
		 & =
		\coll_{J\top}(\coll_{HJ}(s_{H}, t_{J\setminus H}), t_{K \setminus J}, 0_{\top\setminus K})
		 &                 &
		\text{ by } \eqref{eq:restrictedCollar2}
		\\
		 & =
		\coll_{J\top}
		\left(
		\coll_{IJ}
		(
		\coll_{HI}(s_{H}, t_{I\setminus H}), t_j ), t_{K \setminus J}, 0_{\top\setminus K}
		\right)
		 &                 &
		\text{ by } \hyperref[eq:collars_relation]{(R_{HIJ})}
		\\
		 & = \coll_{I\top}
		\left(
		\coll_{HI}(s_{H}, t_{I \setminus H}), t_j ,t_{K \setminus J}, 0_{\top\setminus K}
		\right)
		 &                 &
		\text{ by } \hyperref[eq:collars_relation]{(R_{IJ\top})}
		\\
		 & =
		\coll_{IK}
		\left(
		\coll_{HI}(s_{H}, t_{I \setminus H}), t_j ,t_{K \setminus J}
		\right)
		 &                 &
		\text{ by } \eqref{eq:zeroingCoordinates}
		\\
		 & =
		\coll_{HK}(s_{H}, t_{K\setminus H})
		 &                 &
		\text{ by } \hyperref[eq:collars_relation]{(R_{HIK})},
	\end{align}
	which proves $\hyperref[eq:collars_relation]{(R_{HJK})}$. An analogous argument, replacing $ \coll $ by $ \lcol $ and $ \top $ by $ F $,  yields $\hyperref[eq:collars_relation_2]{(\hat R_{HJK})}$.
\end{proof}

\subsection{Extension of Weinstein charts}
\label{sec:extend-weinstein}
In this subsection we prove Proposition~\ref{prop:inductive_construction_of_Weinstein}.
Fix $J\in \cc{F}_S$. From the already-constructed relative Weinstein charts $ \vp_{IJ} $  on the codimension-$1$ faces (with $ I=J\setminus \{j\} $), we first build a partial Weinstein chart on $L_J$ compatible with the face charts. An extension argument (Lemma~\ref{lem:Weinstein_Extension}) then upgrades this partial chart to a full relative Weinstein chart for $W_J$, preserving compatiblity on overlaps and the vertical boundary condition along $\vlg$.

\begin{lemma}
	\label{lem:Weinstein_Extension}
	Identify $L_J$ with the zero section in $T^\ast L_J$.
	Let $V_0,V_1 \subset T^\ast L_J$ be open sets with $L_J \subset V_0 \cup V_1$.  Let $\vp_0: V_0 \to X_{\Delta_J}$ be a partial relative Weinstein chart for $ W_J $.
	Then, there exists a neighborhood $V_J \subset V_0 \cup V_1$ of $L_J$ and a relative Weinstein chart $\vp_J: V_J \to X_{\Delta_J}$ for $ W_J $ such that $\vp_J = \vp_0$ on $ V_J \setminus V_1$.
\end{lemma}

\begin{proof}
	For brevity set $ L=L_J $, $ \Psi=\Psi_J $, $ X = X_{\Delta_J} $, $ \vlg=  \vlg_{\Delta_J}  $,  $L_0 = L_J \cap V_0$ and $L_1 =L_J \cap V_1$.
	By~\cite[Lemma 2.8]{SolomonYuval20} choose an arbitrary partial relative Weinstein chart $\vp_1: V_1 \to X$. Moreover, by~\cite[Proposition 2.2]{Solomon14}, after shrinking $V_0$ we can choose $\vp_1$ so that
	\begin{equation}
		\label{eq:Weinstein_differentials}
		d\vp_1|_q=d\vp_0|_q,
		\quad \quad
		\forall
		q \in L_0 \cap L_1.
	\end{equation}
	Write $\pi_L: T^\ast L \to L$ for the canonical projection and write points of $T^*L$ as pairs $(q,\xi) \in T^\ast L$ with $q \in L$ and $\xi \in T_q^\ast L$. Fix a connection $ \nabla $  on $X$ for which $ \vlg $ is totally geodesic.  Choose an open neighborhood of the zero section $U \subset \Psi^*TX$ such that the restriction of $\exp^\nabla$ to $U$ is a diffeomorphism onto its image, and 
    \begin{equation}
        \label{eq:V-inverse-image}
        (\exp^\nabla|_U)^{-1}(\vlg) 
        \cap \Psi^*TX|_{\tns{C}{_0^1}(L)}
        \subset \Psi^*T\vlg.
    \end{equation}
    After shrinking $V_i$ we can assume $\vp_i(V_i) \subset  \exp^\nabla(U)$. Then, there exist unique sections $Y_i \in \Gamma( V_i,  \pi_L^\ast \Psi ^\ast TX)$ such that
	\begin{equation}
		\label{eq:Weinstein_sections}
		\vp_i (q,\xi )= \exp_{\Psi(q)}^{\nabla}( Y_i(q,\xi) ),
		\quad \quad
		Y_i|_{L_i} \equiv 0.
	\end{equation}
	Let $\chi_0, \chi_1$ be a partition of unity for the cover $L_0,L_1$ of $L$. 	Define 
    \[
    \vp:V_0\cup V_1\to X
    \]
    by
	\begin{equation}
		\label{eq:Weinstein_interpolation}
		\vp(q,\xi)=\exp_{\Psi(q)}^\nabla
			\bigl(\chi_0(q)Y_0(q,\xi)+\chi_1(q)Y_1(q,\xi)\bigr).
	\end{equation}
	After shrinking $V_0$ and $V_1$, we may assume $V_i \subset \pi_L^{-1}(L_i)$. 
    Then
	\begin{equation}
		\vp|_{V_0\setminus V_1} \equiv \vp_0, \quad \quad
		\vp|_{V_1 \setminus V_0} = \vp_1, \quad \text{ and }\quad
		\vp|_L=\Psi.
	\end{equation}
	Let $ N^\ast_v $ be the conormal bundle of the vertical boundary of $ L $. Since $ \vlg $ is totally geodesic, using the relative condition \eqref{itm:conormalBundle} in Definition \ref{def:WeinsteinModel} and inclusion~\eqref{eq:V-inverse-image}, we have
	\begin{equation}
		(q,\xi) \in N_v^* \cap V_i
        \quad 
		\Longrightarrow
        \quad 
		Y_i(q,\xi) \in T_{\Psi(q)}\vlg.
	\end{equation}
    It follows that $\vp(N_v^*\cap(V_0\cup V_1))\subset\vlg$.

	We now show that
    \begin{equation}
        \label{eq:equality-of-Wein-differentials}
        d\vp|_q = d\vp_0|_q = d\vp_1|_q,
        \qquad 
        \text{ for all }
        \quad 
        q \in L_0 \cap L_1.
    \end{equation}
    For $q \in L_0 \cap L_1$ and $w \in T_q (T^\ast L)$, decompose $w=w_{TL}+ w_{T^\ast L}$ via the canonical direct sum decomposition $TT^\ast L|_L \simeq T L \oplus T L^\ast $. Since $\vp|_L = \Psi$,
	\begin{equation}
		d\vp|_q(w_{TL})= d\Psi|_{q}(w_{TL})= d\vp_0|_q(w_{TL}).
	\end{equation}
    For the vertical part, put $u=w_{T^\ast L}$. Then
	\[
		\begin{aligned}
			d\vp|_q(u)
			&=\left. d\exp_{\Psi(q)}^\nabla\right|_0
				\left(\sum_{i=0}^1\chi_i(q)\,dY_i|_q(u)\right)
            &&
    		\text{by } \eqref{eq:Weinstein_interpolation}
    		\\
			&=\sum_{i=0}^1\chi_i(q)\,
				\left. d\exp_{\Psi(q)}^\nabla\right|_0\bigl(dY_i|_q(u)\bigr)\\
			&=\sum_{i=0}^1\chi_i(q)\,d\vp_i|_q(u)
            		 &         &
		      \text{by  }\eqref{eq:Weinstein_sections}
            \\
			     &=d\vp_0|_q(u),
            		 &         &
		      \text{by } \eqref{eq:Weinstein_differentials}
		\end{aligned}
	\]
    proving equation~\eqref{eq:equality-of-Wein-differentials}.
	
    Since $\vp_0$, $\vp_1$ are local diffeomorphisms, $d\vp|_L$ is invertible; by the inverse function theorem, there exists a neighborhood $V \subset V_0\cup V_1$ of $L$ such that $\vp|_V$ is a local diffeomorphism.     Since $\vp(N_v^*\cap V) \subset \vlg$ and $\dim N_v^*= \dim \vlg$, the restriction of $\vp$ to $N_v^*$ is locally a diffeomorphism onto a neighborhood in $\vlg$ at each vertical boundary point $(q,0)$. By local injectivity of $\vp$, after shrinking $V$ we obtain
    \begin{equation}
        \label{eq:conormal-preserved}
        \vp(q,\xi)\in\vlg \quad\Longleftrightarrow\quad (q,\xi)\in N_v^*,
        \qquad 
        (q,\xi) \in V, \quad q \in \tns{C}{_0^1}(L).
    \end{equation}  
    Replacing $V_i$ by $V_i \cap V$ we now have a local diffeomorphism $\vp: V \to X$ which additionally satisfies
    \begin{equation}
		\label{eq:modelProperties}
		\begin{aligned}[t]
			&\vp|_L=\Psi,\qquad \vp|_{V_0\setminus V_1}=\vp_0,\qquad
				\vp|_{V_1\setminus V_0}=\vp_1.
		\end{aligned}
	\end{equation}
    
	To turn $ \vp $ into a local \emph{symplectomorphism} on $V$, we now apply a Moser-type adjustment.
	Let $\omega_0$ be the canonical form on $T^\ast L$ and set
	\begin{equation}
		\tau = \vp^*\omega_X-\omega_0.
	\end{equation}
	By \eqref{eq:modelProperties} we have $ \tau \equiv 0 $ on   $V_0\setminus V_1$, $V_1 \setminus V_0.$ It follows from~\eqref{eq:equality-of-Wein-differentials} that  
    \begin{equation}\label{eq:res_tau_L}
    \tau|_L \equiv 0 \in \Gamma(L,T^*L|_L).
    \end{equation}
    Moreover, the pullback of $ \tau $ to  $ N_v^* \cap V $ vanishes since $ \vp(N_v^\ast \cap V) \subset \vlg $.
	Define
	\begin{equation}
		\phi:  [0,1] \times V \to V, \quad \quad \phi_t(q,\xi) := (q,t\xi).
	\end{equation}
	Then $\phi_0^\ast \tau = 0$ because $\phi_0(V) \subset L,$ and $\phi_1^\ast \tau = \tau$ since $\phi_1 = \text{id}_V$.
	For $t>0$, $\phi_t$ is a diffeomorphism; set
	\begin{equation}
		X_t := \left( \frac{d}{dt} \phi_t \right) \circ \phi_t^{-1},
		\quad \quad
		\sigma_t:= \phi_t^\ast (\iota_{X_t}\tau),
		\quad \quad
		\sigma := \intop_0^1 \sigma_t dt.
	\end{equation}
	Note that $\sigma \equiv \sigma_t \equiv 0$ on $V_0 \setminus V_1$, $V_1 \setminus V_0$ and $\sigma|_L \equiv \sigma_t|_L \equiv  0 \in \Gamma(L,T^*L|_L)$. Since $ \phi $ preserves $ N^*_v $, the pullback of $ \sigma $ to $ N_v^* \cap V $ vanishes.     We have
	\begin{align}
		\tau
		 & = \phi_1^\ast \tau - \phi_0^\ast \tau
		= \intop_0^1  \frac{d}{dt} (\phi_t^\ast \tau)\, dt
		=  \intop_0^1 \phi_t^\ast (\cc{L}_{X_t}\tau)\, dt
		= \intop_0^1 d\sigma_t dt =
		d\sigma.
	\end{align}
	For $ t\in [0,1] $ define  $\omega_t:= \omega_0 + t \cdot \tau $. By~\eqref{eq:res_tau_L}, shrinking $ V $ if necessary, we can assume $\omega_t$ is non-degenerate on $ V $. Thus, there exists a family of vector fields $Y_t$ on $V$ such that:
	\begin{enumerate}
		\item $\sigma + \iota_{Y_t} \omega_t  =0$,
		\item $Y_t \equiv 0$  on  $V_0 \setminus V_1$, $V_1 \setminus V_0$ and $L$,
		\item $ Y_t $ is tangent to $ N_v^*\cap V $.
	\end{enumerate}
    
    If $\tns{C}{_0^1}(L)\ne\emptyset$, the field $Y_t$ may point outward
    along the vertical boundary of~$V$, so its flow need not be defined.  
    To deal with this problem, let $\tilde L$ be an extension of $L$ obtained by gluing a collar $\tns{C}{_0^1}(L) \times [0,\varepsilon)$ to $L$ along $\tns{C}{_0^1}(L)$, so that $L \subset \tilde L$.
    By Seeley's Extension Theorem~\cite{Seeley64}, extend $\vp$ and
    $\sigma$ to smooth $\tilde \vp$ and $\hat\sigma$
    on a neighborhood of the zero section $\tilde V \subset T^* \tilde L$. Let $\tilde \omega_0$ be the canonical symplectic form on $T^*\tilde L$, and set 
\[
\tilde \tau = \tilde \vp^*\omega_X-\tilde \omega_0, \qquad
\tilde \omega_t := \tilde\omega_0 + t\cdot \tilde \tau.
\]
The form
    \[
        \alpha=\tilde \tau -d\hat\sigma
    \]
    is closed and vanishes on $V$. By the relative Poincaré lemma, there is a one-form $\beta$ with $d\beta=\alpha$ and
    $\beta|_V=0$. Thus $\tilde \sigma := \hat\sigma+\beta$ extends $\sigma$ and satisfies
    $d\tilde \sigma =\tilde \tau$. Define $\tilde Y_t$ by
    \[
        \tilde\sigma + \iota_{\tilde Y_t} \tilde\omega_t =0.
    \]
    Since $\tilde Y_t$ extends $Y_t$, it vanishes
    on $L$, hence its flow $\tilde\psi_t$ is defined on an open $\tilde V' \subset \tilde V$ such that $L \subset \tilde V'$
    for all $t\in[0,1]$.
	Moreover, $\tilde\psi_t$ is the identity on $(V_0 \setminus V_1) \cap \tilde V'$, $(V_1 \setminus V_0) \cap \tilde V'$ and $L$, and preserves $ N_v^* \cap V $.  We have
	\begin{equation}
		\frac{d}{dt} \tilde\psi_t^\ast \tilde\omega_t = \tilde\psi_t^\ast \left( \frac{d}{dt} \tilde\omega_t + \iota_{\tilde Y_t} d\tilde\omega_t + d \iota_{\tilde Y_t} \tilde\omega_t \right) = \tilde\psi_t^\ast (d\tilde\sigma + d\iota_{\tilde Y_t}\tilde\omega_t) = 0,
	\end{equation}
	so $\tilde\psi_1^\ast \tilde\vp^* \omega_X = \tilde\psi_1^\ast \tilde\omega_1 = \tilde\psi_0^\ast \tilde\omega_0 = \tilde\omega_0$. 
    
    Finally, set $V_J = V \cap \tilde V'$ and $\vp_J : =\tilde\vp \circ \tilde\psi_1|_{V_J}: V_J \to X$. Then $\vp_J$ is a local symplectomorphism, and it still satisfies~\eqref{eq:conormal-preserved}. It follows from~\eqref{eq:modelProperties} that $\vp_J = \vp_0$ on $ V_J \setminus V_1.$
\end{proof}

\begin{proof}[Proof of Proposition \ref{prop:inductive_construction_of_Weinstein}]
	For each $j \in J$, set $I = J \setminus \{j\}$. Choose open neighborhoods of the zero section
	\begin{equation}
		V_{IJ}' \subset V_{IJ} \subset T^\ast L_{IJ},
		\quad
		\text{ with }
		\quad
		\cc{T}g_{IJ}(V_{IJ}' ) \subseteq V_{IJ}.
	\end{equation}
	Using the collar embedding $ L_{IJ} \hookrightarrow L_J $ in Remark \ref{rem:collarIdentifications}, we regard $ V_{IJ}' $ as a subset of $ T^\ast L_J $ and define
	\begin{equation}
		V_J' = \bigcup_{I=J\setminus \{j\}} V_{IJ}' \qquad  \subseteq T^\ast L_J.
	\end{equation}
	Define $ \vp_J':V_J' \to X_{\Delta_J} $ piecewise by
	\begin{equation}
		\vp_J'|_{V_{IJ}'} = \vp_{IJ} \circ \cc{T}g_{IJ}.
	\end{equation}
	We first claim that $\vp_J'$ is well defined. Let $j_1, j_2 \in J$, and set
	\begin{equation}
		I_1 = J \setminus \{j_1\}, \quad \quad
		I_2 = J \setminus \{j_2\}, \quad \quad
		H = J \setminus \{j_1,j_2\}.
	\end{equation}
	We can assume that $V_{I_1J}' \cap V_{I_2J}' \subset V_{HJ}$. Since $\vp_{I_1}$ and $\vp_{I_2}$ are compatible with $\vp_H$ we have:
	\begin{align}
		(\vp_{I_1J} \circ \cc{T}g_{I_1J})|_{V_{HJ}}
		 & = \vp_{HJ} \circ \cc{T}(g_{HI_1} \circ \pi_{HI_1}^{HJ}) \circ \cc{T}(g_{I_1J} \circ \lcol_{I_1J}^{HJ})
		 &                                                                                                        &
		\text{ by }\eqref{eq:extended_compatible_Weinstein_symplecto}
		\\
		 & = \vp_{HJ} \circ \cc{T}(g_{HJ})
		 &                                                                                                        &
		\text{ by }\hyperref[eq:compatible_generating_functions]{(G_{HI_1J})}
		\\
		 & =
		\vp_{HJ} \circ \cc{T}(g_{HI_2} \circ \pi_{HI_1}^{HJ}) \circ \cc{T}(g_{I_2J} \circ \lcol_{I_2J}^{HJ})
		 &                                                                                                        &
		\text{ by }\hyperref[eq:compatible_generating_functions]{(G_{HI_2J})}
		\\
		 & = (\vp_{I_2J} \circ \cc{T}g_{I_2J})|_{V_{HJ}}
		 &                                                                                                        &
		\text{ by }\eqref{eq:extended_compatible_Weinstein_symplecto}
	\end{align}
	so $\vp_J'$ is well defined on $V_{I_1J}' \cap V_{I_2J}'$.

    The map $\vp_J'$ is a partial relative Weinstein chart on $V_J'$. Choose
    $0<\varepsilon'<\varepsilon$ and let $V_1$ be the open set of
    covectors based outside the closures of the
    $\widetilde L_{IJ}(\varepsilon')$ for $|J\setminus I|=1$.
    By Lemma~\ref{lem:Weinstein_Extension} applied with $V_0=V_J'$,
    there exists a neighborhood $V_J\subseteq T^\ast L_J$ of the zero
    section and a relative Weinstein chart $\vp_J:V_J\to X_{\Delta_J}$
    agreeing with $\vp_J'$ on $V_J\setminus V_1$.
    After replacing $\varepsilon$ by $\varepsilon'$, the charts $\vp_J$ and $\vp_{J}'$ agree
    on $V_J|_{L_{IJ}}$ for each $I$ with $|J\setminus I|=1$.
    Shrink $V_J$ so that $V_J|_{L_{IJ}}\subset V_{IJ}'$ for each such $I$.

	Fix $ H \subseteq I \subsetneq J$ with $|J\setminus I|=1$. It is left to prove \eqref{eq:compatible_Weinstein_symplecto_inclusion} and \eqref{eq:compatible_Weinstein_symplecto} for $H,J$. We have:
	\begin{align}
		V_J|_{L_{HJ}}
		 & =
		\left(V_J|_{L_{IJ}} \right)|_{L_{HJ}}
		\\
		 & \subseteq V_{IJ} |_{L_{HJ}}
		 &                                                                  &
		\text{since } V_J|_{L_{IJ}} \subseteq V_{IJ}' \subseteq V_{IJ}
		\\
		 & = V_I|_{L_{HI}} \times T^\ast [0,\varepsilon)^{|J\setminus I|}
		 &                                                                  &
		\text{by } \eqref{eq:productNbhdRestriction}
		\\
		 & \subseteq V_{HI} \times  T^\ast [0,\varepsilon)^{|J\setminus I|}
		 &                                                                  &
		\text{by } \eqref{eq:compatible_Weinstein_symplecto_inclusion} \text{ for } H,I
		\\
		 & = V_{HJ}
	\end{align}
	and
	\begin{align}
		\MoveEqLeft[1] \cc{T} g_{HJ}(V_J|_{L_{HJ}}) \notag\\
		 & =
		\cc{T}(g_{HI} \circ \pi_{HI}^{HJ})
		\left(
		\cc{T}(g_{IJ} \circ \lcol_{IJ}^{HJ})
			(\left(V_J|_{L_{IJ}} \right)|_{L_{HJ}})
		\right)
		 &                                                                  &
		\text{by } \hyperref[eq:compatible_generating_functions]{(G_{HIJ})}
		\\
		 & \subseteq
		\cc{T}(g_{HI} \circ \pi_{HI}^{HJ}) (V_{IJ}|_{L_{HJ}})
		 &                                                                  &
		\text{by } \eqref{eq:compatible_Weinstein_symplecto_inclusion} \text{ for } I,J
		\\
		 & =
		\cc{T}(g_{HI} \circ \pi_{HI}^{HJ}) \left( V_I|_{L_{HI}} \times T^\ast [0,\varepsilon)^{|J\setminus I|} \right)
		 &                                                                  &
		\text{by } \eqref{eq:productNbhdRestriction}
		\\
		 & \subseteq V_{HI} \times  T^\ast [0,\varepsilon)^{|J\setminus I|}
		 &                                                                  &
		\text{by } \eqref{eq:compatible_Weinstein_symplecto_inclusion} \text{ for } H,I
		\\
		 & = V_{HJ}
	\end{align}
	which proves \eqref{eq:compatible_Weinstein_symplecto_inclusion} for $H,J$. 	Finally, we have
	\begin{align}
		\vp_J|_{L_{HJ}}
		 & = (\vp_J|_{L_{IJ}})|_{L_{HJ}}
		\\
		 & = (\vp_{IJ} \circ \cc{T}g_{IJ})|_{L_{HJ}}
		 &                                                                                    &
		\text{since }  \vp_J|_{L_{IJ}} = \vp_J'|_{L_{IJ}}
		\\
		 & = \vp_{HJ} \circ \cc{T}(g_{HI} \circ \pi_{HI}^{HJ}) \circ \cc{T}( g_{IJ} \circ \lcol_{IJ}^{HJ})
		 &                                                                                    &
		\text{by }\eqref{eq:extended_compatible_Weinstein_symplecto}
		\\
		 & = \vp_{HJ} \circ  \cc{T}g_{HJ},
		 &                                                                                    &
		\text{by }\hyperref[eq:compatible_generating_functions]{(G_{HIJ})}
	\end{align}
	proving \eqref{eq:compatible_Weinstein_symplecto} for $H,J$, which completes the proof.
\end{proof}

\subsection{Generating functions on collars}
\label{sec:gen-fn}

In this subsection we prove Proposition \ref{prop:existence_of_generating_function}. After precomposing the immersion $ \Psi_{JK} $ with a diffeomorphism $\mu_{JK} : L_{JK} \to L_{JK}$,  we factor it as
\begin{equation}
	L_{JK} \xrightarrow{\sigma_{JK}} V_{JK} \xrightarrow{\vp_{JK}} X_{\Delta_{JK}}
\end{equation}
where $ \sigma_{JK} = dg_{JK} $ is an exact section of $ T^\ast L_{JK} $ vanishing on the vertical boundary. Prior compatibility gives $\mu_{JK}|_{L_{IK}} = \mathrm{id}_{L_{IK}}$ for every $I \subsetneq J$, so replacing $ (\coll_{JK}, \lcol_{JK}) $ by $ (\coll_{JK}, \lcol_{JK} \circ \mu_{JK}) $ gives a compatible collar pair, and $((\coll_{JK},\lcol_{JK} \circ \mu_{JK}), \vp_{JK}, g_{JK})$ is an exact collar model for $ W_J$ in $ W_K $.

\begin{lemma}
	\label{lem:almostSection}
	Let $J \subsetneq K$ in $\cc{F}_S$.
	Let $(\coll_{JK},\lcol_{JK})$ be a collar pair as in Proposition~\ref{prop:extension_of_collar_pairs}, and let
	$\Psi_{JK}\colon L_{JK}\to X_{\Delta_{JK}}$ be the associated representation of $W_{JK}$.
	Let $\vp_J\colon V_J\to X_{\Delta_J}$ be a relative Weinstein chart for $W_J$ as in Proposition~\ref{prop:inductive_construction_of_Weinstein}, and set
	\[
		V_{JK}:=V_J\times T^*\!\big([0,\varepsilon)^{|K\setminus J|}\big)\subset T^*L_{JK},
		\qquad
		\vp_{JK}:=\vp_J\times\mathrm{id}.
	\]
	For $\varepsilon>0$ sufficiently small, there exists a smooth map
	$\theta_{JK}\colon L_{JK}\to V_{JK}$ such that:
	\begin{enumerate}
		\item \label{itm:sectionRepresentation}
		      $\vp_{JK}\circ \theta_{JK}=\Psi_{JK}$.
		\item \label{itm:conormalSection}
		      Along the vertical boundary, $\theta_{JK}$ takes values in the conormal:
		      \[
			      \theta_{JK}\big(\,\tns{C}{_0^1}(L_{JK})\,\big)\ \subset\ N^{*}_{JK}
			      \quad
			      \text{(the conormal of $\tns{C}{_0^1}(L_{JK})\subset L_{JK}$).}
		      \]
		\item \label{itm:almostSection}
		      The base projection is a diffeomorphism: with $\pi_{JK}\colon T^*L_{JK}\to L_{JK}$,
		      \[
			      \pi_{JK}\circ\theta_{JK}\colon L_{JK}\longrightarrow L_{JK}
			      \quad\text{is a diffeomorphism.}
		      \]
		\item \label{itm:restrictionSection}
		      For every $\emptyset \neq I\subsetneq J$, under the inclusion $L_{IK}\subset L_{JK}$,
		      \[
			      (\pi_{JK}\circ\theta_{JK})\big|_{L_{IK}}=\mathrm{id}_{L_{IK}}.
		      \]
	\end{enumerate}
\end{lemma}

\begin{proof}
	Recall the product decomposition
	\begin{equation}
		X_{\Delta_{JK}} \simeq X_{\Delta_{J}} \times T^* [0,\varepsilon)^{|K \setminus J|}.
	\end{equation}
	Let
	\begin{equation}
		\Psi_1 : L_{JK} \to X_{\Delta_{J}},
		\quad \quad
		\Psi_2 : L_{JK} \to T^* [0,\varepsilon)^{|K \setminus J|}.
	\end{equation}
	be the components of $\Psi_{JK}$ under this product. Since $L_J$ is compact and $\vp_J$ is a local diffeomorphism, we can cover the zero section in $T^\ast L_J$ by open sets $U_1,\dots,U_n \subset T^\ast L_J$ on which $\vp_J$ is an embedding. Write
	\begin{equation}
		\vp_i = \vp_J|_{U_i}: U_i \xrightarrow{\simeq} \vp_J(U_i).
	\end{equation}
	Identify $ L_J $ with the zero section in $ T^\ast L_J $. Define
	\begin{equation}
		\theta_1: L_{JK} \to V_J,
		\quad \quad
		\theta_1(q_J,t_{K \setminus J}) = \vp^{-1}_i(\Psi_1(q_J, t_{K\setminus J})),
		\quad \quad \text{if } q_J \in U_i.
	\end{equation}
	We claim that $ \theta_1 $ is well defined for all $ \|t_{K \setminus J}\| $ small enough.
	Fix $q_J \in L_J$ and suppose $q_J \in U_{i_1} \cap \dots \cap U_{i_k}$. By Lemma \ref{lem:cobordismOverProduct} we have
	\begin{equation}
		\Psi_1(q_J,0_{K \setminus J}) = \Psi_J(q_J) \in \vp_J(U_{i_1} \cap \dots \cap U_{i_k}).
	\end{equation}
	We claim there exists $\varepsilon_{q_J}>0$ such that, for
    $q'_J$ near $q_J$ and $t=t_{K\setminus J}\in [0,\varepsilon_{q_J})^{|K\setminus J|}$,
    \begin{equation}
        \label{eq:psi-one-inclusion}
        \Psi_1(q'_J,t)\in \vp_J(U_{i_1}\cap\dots\cap U_{i_k}).
    \end{equation}
    Away from the vertical boundary this follows directly by continuity. If $q_J$ lies on the vertical boundary, let $\rho$ be a signed boundary defining function for the vertical face of the local chart image,
    so that the image is locally given by $\rho\ge 0$, and set $r=\rho\circ\Psi_J$. Since $\varphi_J$ maps the conormal bundle to $\vlg_{\Delta_J}$, $\rho\circ\Psi_1$ vanishes at $r=0$. At $t=0$ we have $\rho\circ\Psi_1=r$, so its inward derivative is $1$. By continuity this derivative stays positive nearby for small $t$. Thus $\rho\circ\Psi_1\ge 0$ for $r\ge 0$. Consequently, for $q'_J$ near $q_J$ and $t_{K\setminus J}$ sufficiently small we get~\eqref{eq:psi-one-inclusion}.
    
	Compactness of $ L_J $ yields a uniform $\varepsilon > 0$ with the same property for all $q_J \in L_J$. Hence if $q_J \in U_i \cap U_j$ and $t_{K\setminus J} \in  [0,\varepsilon)^{|K \setminus J|}$
	\begin{equation}
		\Psi_1(q_J, t_{K \setminus J}) \in \vp_J(U_i \cap U_j) \subset \vp_J(U_i) \cap \vp_J(U_j)
	\end{equation}
	so $ \vp_i^{-1}(\Psi_1(q_J, t_{K \setminus J})) = \vp_j^{-1}(\Psi_1(q_J, t_{K \setminus J})) $. This shows $\theta_1$ is well defined.	Set
	\begin{equation}
		\theta_{JK}: L_{JK} \to V_{JK}, \quad \quad
		\theta_{JK} = (\theta_1, \Psi_2).
	\end{equation}
	
    We now verify properties \eqref{itm:sectionRepresentation}-\eqref{itm:restrictionSection}.
	We have
	\begin{equation}
		\vp_{JK} \circ \theta_{JK} = (\vp_J \circ \theta_1, \Psi_2) = (\Psi_1,\Psi_2) = \Psi_{JK}
	\end{equation}
	which proves property \eqref{itm:sectionRepresentation}. Let $(q_J, t_{K \setminus J}) \in \tns{C}{_0^1}(L_{JK}) \simeq \tns{C}{_0^1}(L_{J}) \times [0,\varepsilon)^{|K \setminus J|}$.
	Since $\Psi_{JK}(q_J, t_{K \setminus J}) \in \vlg_{\Delta_{JK}}$ we have
	\begin{equation}
		\Psi_1(q_J, t_{K \setminus J})
		= \mathrm{pr}_{X_{\Delta_J}}(\Psi_{JK}(q_J, t_{K \setminus J})) \in \mathrm{pr}_{X_{\Delta_J}} (\vlg_{\Delta_{JK}})
		=  \vlg_{\Delta_J},
	\end{equation}
	and $\Psi_2(q_J, t_{K\setminus J})$ is the zero covector in $T^*[0,\varepsilon)^{|K\setminus J|}$.
	By the relative condition \eqref{itm:conormalBundle} of the Weinstein chart in Definition \ref{def:WeinsteinModel}, we have  $\theta_1(q_J, t_{K\setminus J}) \in N^*_J$, hence
	\begin{equation}
		\theta_{JK}(q_J, t_{K\setminus J}) \in N^*_J \times 0_{T^*[0,\varepsilon)^{|K\setminus J|}} \simeq N^*_{JK}.
	\end{equation}
	which proves \eqref{itm:conormalSection}. By Lemma \ref{lem:cobordismOverProduct},
	\begin{equation}
		\theta_1(q_J, 0_{K\setminus J}) = \vp_i^{-1}(\Psi_1(q_J, 0_{K\setminus J})) = \vp_i^{-1}(\Psi_J(q_J)) =q_J,
	\end{equation}
	Let $\pi_J:T^*L_J\to L_J$
	be the bundle projection. By \eqref{eq:collars_over_corner},
	the base coordinate of $\Psi_2(q_J,t_{K\setminus J})$ is $t_{K\setminus J}$. Since
	$\theta_{JK}=(\theta_1,\Psi_2)$, it follows that
	\[
		(\pi_{JK}\circ\theta_{JK})(q_J,t_{K\setminus J})=(h_{t_{K\setminus J}}(q_J),t_{K\setminus J}),
		\quad 
        h_{t_{K\setminus J}}(q_J):=\pi_J\bigl(\theta_1(q_J,t_{K\setminus J})\bigr).
	\]
	The preceding computation gives $h_0=\mathrm{id}_{L_J}$.
	To prove \eqref{itm:restrictionSection}, we have
	\begin{align}
		(\vp_{JK} \circ \theta_{JK})|_{L_{IK}}
		 & =
		\Psi_{JK}|_{L_{IK}}
		 &                                                                          &
		\text{by } \eqref{itm:sectionRepresentation}
		\\
		 & = \Psi_{IK}
		 &                                                                          &
		\text{by } \eqref{eq:restriction_of_collar_representation}
		\\
		 & = \vp_{IK} \circ \Gamma g_{IK}
		 &                                                                          &
		\text{by } \eqref{diag:exact-collar-model}
		\\
		 & = \vp_{JK}|_{L_{IK}} \circ \Gamma( g_{IK} - g_{IJ} \circ \pi_{IJ}^{IK}).
		 &                                                                          &
		\text{by }
		\eqref{eq:compatible_Weinstein_symplecto},
		\eqref{eq:extended_compatible_Weinstein_symplecto}
	\end{align}
	Since $ \vp_{JK}|_{L_{IK}} $ is a local diffeomorphism, this shows  $\theta_{JK}|_{L_{IK}}$ is a section of $V_{JK}|_{L_{IK}}$, and in particular $(\pi_{JK} \circ \theta_{JK})|_{L_{IK}} = \mathrm{id}_{L_{IK}}$, proving \eqref{itm:restrictionSection}.

	By \eqref{itm:conormalSection}, $h_{t_{K\setminus J}}$ preserves the vertical
	boundary, and by \eqref{itm:restrictionSection}, it is the identity on
	$\widetilde L_{IJ}\subset L_J$ for every nonempty $I\subsetneq J$. Thus $h_{t_{K\setminus J}}$ is a diffeomorphism after shrinking $\varepsilon$.
	The formula above for $\pi_{JK}\circ\theta_{JK}$ proves
	\eqref{itm:almostSection}.
\end{proof}

\begin{proof}[Proof of Proposition \ref{prop:existence_of_generating_function}]
	Let $\theta_{JK}$ be as in Lemma \ref{lem:almostSection}. By Lemma \ref{lem:almostSection}\eqref{itm:almostSection}, $ \pi_{JK} \circ \theta_{JK} $ is a diffeomorphism,  so set
	\begin{equation}
		\mu_{JK} = (\pi_{JK} \circ \theta_{JK})^{-1} : L_{JK} \to L_{JK}.
	\end{equation}
	Define
	\begin{equation}
		\sigma_{JK}: L_{JK} \to V_{JK},
		\quad \quad
		\sigma_{JK} = \theta_{JK} \circ \mu_{JK}
	\end{equation}
	Then $ \sigma_{JK} $ is a smooth section of $T^* L_{JK}$, since
	\begin{equation}
		\pi_{JK} \circ \sigma_{JK} = \pi_{JK} \circ \theta_{JK} \circ (\pi_{JK} \circ \theta_{JK})^{-1} = \mathrm{id}_{L_{JK}}.
	\end{equation}
	Moreover, by  Lemma \ref{lem:almostSection}\eqref{itm:conormalSection} the restriction of $\sigma_{JK}$ to $\tns{C}{_0^1}(L_{JK})$ is a section of the conormal bundle $N^*_{JK}$.
	We claim that $\sigma_{JK}$ represents an exact form on $L_{JK}$. First,
	\begin{equation}
		\Psi_{JK} \circ \mu_{JK} = \vp_{JK} \circ \theta_{JK} \circ \mu_{JK} = \vp_{JK} \circ \sigma_{JK}.
	\end{equation}
	Since $\Psi_{JK}$ is a Lagrangian immersion and $ \vp_{JK} $ is a local symplectomorphism,  $\sigma_{JK}$ is closed. Let $\iota: L_J \to L_{JK}$ be the canonical inclusion $q_J \mapsto (q_J, 0_{K \setminus J})$.  Since $\mu_{JK}(q_J,0)=(q_J,0)$, we have
	$\sigma_{JK}(q_J,0)=\theta_{JK}(q_J,0)$.
	Its $T^*L_J$ component is $\theta_1(q_J,0)=(q_J,0)$.
	Thus $\iota^*\sigma_{JK}=0$,
	and $\sigma_{JK}$ defines a class in $H^1(L_{JK}, L_J)$. Since $L_{JK} = L_J \times [0,\varepsilon)^{|K\setminus J|}$ deformation retracts to $L_J$ we have $ H^1(L_{JK}, L_J)=0 $. Therefore, there exists a smooth function $g_{JK}: L_{JK} \to \bb{R}$ with
	\begin{equation}
		dg_{JK}  = \sigma_{JK},
		\quad
		\text{and}
		\quad
		g_{JK}|_{L_J} \equiv 0.
	\end{equation}
	Since $\sigma_{JK}|_{\tns{C}{_0^1}(L_{JK})} \in \Gamma(N^*_{JK}) $,    $g_{JK}$ is constant on $\tns{C}{_0^1}(L_{JK})$. But
	\begin{equation}
		\tns{C}{_0^1}(L_{JK}) \simeq \tns{C}{_0^1}(L_{J}) \times [0,\varepsilon)^{|K \setminus J|}
	\end{equation}
	so each point in $ \tns{C}{_0^1}(L_{JK}) $ is joined by a path to a point in $ L_J $ where $ g_{JK}  \equiv 0$. Therefore $ g_{JK}|_{\tns{C}{_0^1}(L_{JK})} \equiv 0$.

	Set $\lcol_{JK}' := \lcol_{JK} \circ \mu_{JK}$. By \eqref{eq:correspondingRepresentation}, the representation $\widetilde{W}_{JK}$ corresponding to $(\coll_{JK},\lcol_{JK}')$  is
	\begin{equation}
		\Psi_{JK}' := \Psi_{JK} \circ \mu_{JK}.
	\end{equation}
	By Lemma \ref{lem:almostSection}\eqref{itm:sectionRepresentation} we have $\Psi_{JK}' = \vp_{JK} \circ \Gamma{g_{JK}}$. Thus $((\coll_{JK},\lcol_{JK}'), \vp_J, g_{JK})$ will be a compatible exact collar model once we verify that $(\coll_{JK},\lcol_{JK}')$ is a collar pair (Definition \ref{def:LagrangianCollar}) and is compatible with the previously constructed collars.

	First, since $\mu_{JK}(q_J,0)=(q_J,0)$, we have
	\[
		\lcol'_{JK}(q_J,0)=\lcol_{JK}(q_J,0)=\iota(q_J).
	\]
	Second, by \eqref{eq:collars_over_corner} for $\lcol_{JK}$
	and the formula for $\mu_{JK}$ above,
	\[
		\widetilde\tau_{JK}\circ\Psi_{\widetilde\Delta_{JK}}
		\circ\lcol'_{JK}
		=\operatorname{pr}_2\circ\mu_{JK}
		=\operatorname{pr}_2.
	\]
	Thus $(\coll_{JK},\lcol'_{JK})$ is a collar pair.

	Regarding compatibility, for every $\emptyset\neq I \subsetneq J $, Lemma \ref{lem:almostSection}\eqref{itm:restrictionSection} is equivalent to $ \mu_{JK} \circ \iota_{JK}^{IK} = \iota_{JK}^{IK} $. Therefore
	\begin{equation}
		\lcol'_{JK}\circ \iota_{JK}^{IK} = \lcol_{JK} \circ \mu_{JK} \circ \iota_{JK}^{IK} = \lcol_{JK} \circ\iota_{JK}^{IK},
	\end{equation}
	so $(\coll_{JK},\lcol_{JK}')$ satisfy the same compatibility relation \eqref{eq:compatible_collars} as $(\coll_{JK},\lcol_{JK})$. Consequently, $((\coll_{JK},\lcol_{JK}'), \vp_J, g_{JK})$ is an exact collar model for $ W_J$ in $ W_K $.

	Finally, the next induction step uses the fact that
	for $J\subsetneq K\subsetneq F$ we have
	\[
		\lcol'_{JF}(q_J,t_{K\setminus J},0_{F\setminus K})
		= \lcol'_{JK}(q_J,t_{K\setminus J}).
	\]
	Indeed, $\coll_{JF}$ restricts to $\coll_{JK}$ by
	\eqref{eq:restrictedCollar2}, and the graph of $dg_{JF}$ pulls back
	to the graph of $d(g_{JF}|_{L_{JK}})$. Thus both collar maps
	parametrize $W_{JK}$ as a graph in the same chart $\vp_{JK}$.
	In either parametrization, the parameter of a point is recovered
	by lifting its image on the branch determined along $L_J$ and
	applying $\pi_{JK}$. Hence the two collar maps agree.
\end{proof}

\section{Composition of Lagrangian cobordisms}
\label{sec:composition}
In this section we prove Theorem~\ref{thm:hornFilling}.
Let $\Lambda_i[n]$ be the $(n,i)$-horn as in Definition~\ref{def:n-i-horn}, and let $\cc{F}_{\Lambda_i[n]}$ be the corresponding family of faces as in Definition~\ref{def:family-of-faces}.
Fix a Lagrangian $W$ over ${\Lambda_i[n]}$ as in Definition~\ref{def:Lag-over-sSet}, and let $\cc{F}_W$ be the corresponding family of cobordisms as in Definition~\ref{def:family-of-cobordisms}.
Write $\top=[n]$. By Theorem~\ref{thm:compatible_exact_models} we obtain a compatible family of exact collar models $\cc{C}_{IJ} = \big((\coll_{IJ},\lcol_{IJ}), \vp_I, g_{IJ}\big) $ for $I \subsetneq J $ in $\cc{F}_{\Lambda_i[n]}$.
In Subsection~\ref{sec:extend-gen-empty} we use the cocycle condition \eqref{eq:compatible_generating_functions} of Proposition \ref{prop:compatible_generating_functions} on the transition potentials $g_{JK}$ to construct, in two steps, a family of smooth functions
\begin{equation}
	\label{eq:the-FJ-family}
	\left\{F_J: L_{J\top} \to \bb{R}\,\,|\,\, J \in \cc{F}_{\Lambda_i[n]} \right\}
\end{equation}
extending the transition potentials $g_{JK}$ to the product collar $L_{J\top}$ and satisfying the compatibility condition \eqref{eq:gen_function_extension_compatibility} whenever $I \subseteq J$.
In Subsection~\ref{sec:ext_Lag_horns} we use the functions $F_J$ to generate Lagrangian cobordisms $W_{J\top}$ over $\Delta_{J\top}$, the collar of $\Delta_J$ in $\Delta^n=\Delta_\top$. Finally, in Lemma~\ref{lem:globalLag} we show that the compatibility \eqref{eq:gen_function_extension_compatibility} implies the cobordisms $W_{J\top}$ agree on overlaps; their gluing allows us to to produce the desired
horn extension and complete the proof of Theorem~\ref{thm:hornFilling}.

\subsection{Extension of generating functions}
\label{sec:extend-gen-empty}
The first step of in the construction of the family \eqref{eq:the-FJ-family} is carried out in the following
Lemma, where we smoothly extend the potentials $g_{JK}$,
defined on the boundary faces of $L_{J\top}$, to a function on all of
$L_{J\top}$. The proof is based on a standard inclusion--exclusion argument for extending
smooth compatible functions from the faces of a manifold with corners to its interior~\cite[Proposition 12.3.7]{Hofer}.

\begin{lemma}
	\label{lem:inclusionExclusion}
	Write $\top = [n]$. Let $ J \in \cc{F}_{\Lambda_i[n]}$ with $i \in J$. 
    For $x=(q_J,t_{\top\setminus J})\in L_{J\top}$, define
	$\widetilde F_J:L_{J\top}\to\bb{R}$ by
	\begin{equation}
		\label{eq:inclusionExclusion}
		\widetilde F_J(x)=
		\begin{cases}
			(-1)^n\sum_{\substack{J\subsetneq K\\ K\in\cc{F}_{\Lambda_i[n]}}}
				(-1)^{|K|}g_{JK}(q_J,t_{K\setminus J}), & |J|\le n-1,\\
			0, & |J|=n.
		\end{cases}
	\end{equation}
	Then $\widetilde F_J$ is smooth, vanishes on ${\tns{C}{_0^1}(L_{J\top})} $ and the corner $L_J$ and restricts to the transition potentials on faces: for all $J \subsetneq K$ in $\cc{F}_{\Lambda_i[n]}$
	\begin{equation}
		\label{eq:gen_function_extension_zeroing}
		\widetilde F_J|_{L_{JK}}=g_{JK},
		\quad \text{i.e.} \quad
		\widetilde F_J(q_J,t_{K\setminus J}, 0_{\top\setminus K}) = g_{JK}(q_J, t_{K \setminus J}).
	\end{equation}
\end{lemma}

\begin{figure}[htbp]
	\centering
	\begin{tikzpicture}[scale=1.0, line join=round, line cap=round]
		\coordinate (A) at (0,2.887);
		\coordinate (B) at (-2.5,-1.443);
		\coordinate (C) at (2.5,-1.443);

		\def\t{0.4}

		\coordinate (AtoC) at ($(A)!\t!(C)$);
		\coordinate (BtoC) at ($(B)!\t!(C)$);
		\coordinate (BtoA) at ($(B)!\t!(A)$);
		\coordinate (CtoA) at ($(C)!\t!(A)$);
		\coordinate (CtoB) at ($(C)!\t!(B)$);
		\coordinate (AtoB) at ($(A)!\t!(B)$);

		\coordinate (center) at (barycentric cs:A=1,B=1,C=1);
		\coordinate (CentoB) at ($(center)!\t!(B)$);
		\coordinate (CentoC) at ($(center)!\t!(C)$);
		\coordinate (CentoB1) at ($(CentoB)+(0,0.5)$);
		\coordinate (CentoC1) at ($(CentoC)+(0,0.5)$);
		\coordinate (CentoC2) at ($(C)!0.55!(center)$);
		\coordinate (CentoA2) at ($(A)!0.55!(center)$);
		\coordinate (CentoA) at ($(center)!\t!(A)$);
		\filldraw[gray!20, draw=black, thick] (A)--(B)--(C)--cycle;


		\node[above] at (A) {$-g_{0,03}(q_0,t_3)$};
		\node[below left] at (B) {$-g_{0,01}(q_0,t_1)$};
		\node[below right] at (C) {$-g_{0,02}(q_0,t_2)$};
		\node[gray, left] at ($(A)!0.5!(B)$) {$g_{0,013}(q_0,t_1,t_3)$};
		\node[below, gray] at ($(B)!0.5!(C)$) {$g_{0,012}(q_0,t_1,t_2)$};
		\node[gray, right]  at ($(C)!0.5!(A)$) {$g_{0,023}(q_0,t_2,t_3)$};

		\fill[blue] (center) circle (1pt);
		\draw[thick] (center) -- (B);
		\draw[thick] (center) -- (C);
		\draw[thick] (center) -- (A);
	\end{tikzpicture}
	\caption{The function $\widetilde F_0= -g_{0,01}-g_{0,02}-g_{0,03}+g_{0,012}+g_{0,013}+g_{0,023}$ is defined by inclusion-exclusion principle on the faces of the 3-horn $\Lambda_0^3$.}
	\label{fig: InclusionExclusion}
\end{figure}
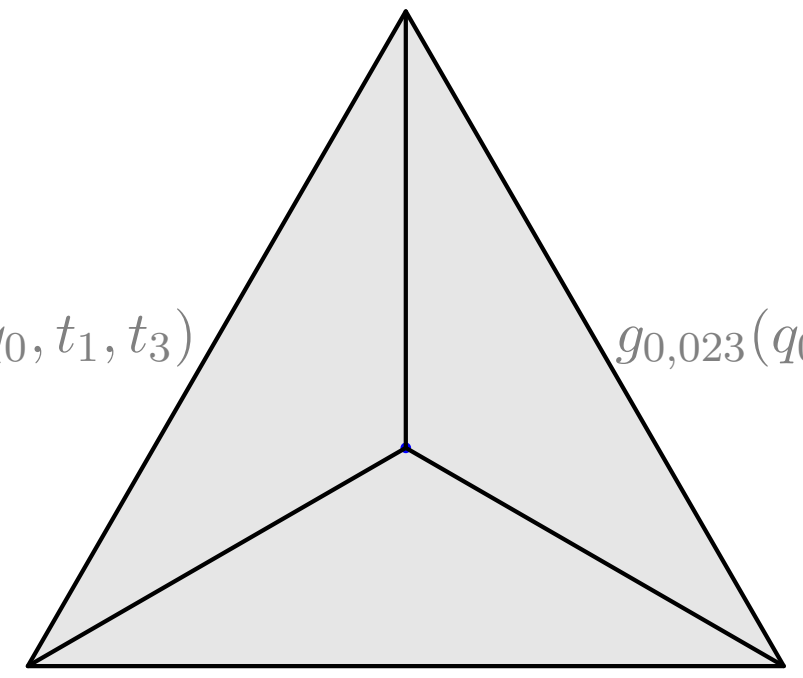
\begin{example}
	Let $n=3$ and consider the 3-horn $\Lambda_0^3$ (Figure \ref{fig: InclusionExclusion}). Let $J=\{0\}$, and define
	\begin{align}
		\widetilde F_0(q_0,t_1,t_2,t_3) = -g_{0,01}(q_0,t_1)-g_{0,02}(q_0,t_2)-g_{0,03}(q_0,t_3)
		\\
		+g_{0,012}(q_0,t_1,t_2)+g_{0,013}(q_0,t_1,t_3)+g_{0,023}(q_0,t_2,t_3).
	\end{align}
	Using the zeroing identity \eqref{eq:zeroingCoordinatesGenFunc}, we get:
	\begin{align}
		\widetilde F_0 |_{L_{0,012}} (q_0,t_1,t_2)
		 & = \widetilde{F}_0(q_0,t_1,t_2,0)
		\\
		 & = g_{0,012}(q_0,t_1,t_2) + (g_{0,013}(q_0,t_1,0) - g_{0,01}(q_0,t_1))
		\\
		 & \quad \quad + (g_{0,023}(q_0,t_2,0)-g_{0,02}(q_0,t_2))
		\\
		 & = g_{0,012}(q_0,t_1,t_2).
	\end{align}
\end{example}

\begin{proof}[Proof of Lemma \ref{lem:inclusionExclusion}]
	If $|J|=n$ there is nothing to prove, so assume $|J|\leq n-1$. Fix \(k\in \top\setminus J\) with $k\neq i$ and let $N = \top\setminus \{k\}$ be a facet containing $J$.  Note
	\[
		(-1)^n \cdot (\widetilde F_J(q_J,t_{N\setminus J},0_k)-g_{J N}(q_J,t_{N\setminus J}))
		=\sum_{\substack{J \subsetneq K \\ K \in \cc{F}_{\Lambda_i[n]}\setminus N}}(-1)^{|K|}\,g_{JK}(q_J,t_{K\setminus J}).
	\]
	Split the index set into
	\[
		\mathcal R=\{J\subsetneq R:\ k \notin R\},
		\qquad
		\mathcal S=\{J\subsetneq S:\ k\in S\}.
	\]
	For \(S\in\mathcal S\) write \(R_S:=S\setminus\{k\}\). Using the zeroing identity \eqref{eq:zeroingCoordinatesGenFunc} we get
	\begin{align}
		 & \
		(-1)^n \cdot (\widetilde F_J(q_J,t_{N\setminus J},0_k)-g_{JN}(q_J,t_{N\setminus J})) \\
		 & =\sum_{R\in\mathcal R}(-1)^{|R|}g_{JR}(q_J,t_{R\setminus J})
		\;+\;\sum_{S\in\mathcal S}(-1)^{|S|}g_{JS}(q_J,t_{R_S\setminus J},0_k)               \\
		 & =\sum_{R\in\mathcal R}(-1)^{|R|}g_{JR}(q_J,t_{R\setminus J})
		\;-\;\sum_{R_S\in\mathcal R}(-1)^{|R_S|}g_{JR_S}(q_J,t_{R_S\setminus J})
		\\
		 & =0,
	\end{align}
	since the two sums cancel term-by-term. Thus, for every facet $N = \top \backslash \{k\}$ we have
	\begin{equation}
		\label{eq:facetInclusionExclusion}
		\widetilde F_J|_{L_{JN}}=g_{JN}.
	\end{equation}
	Now let \( J\subsetneq K\), pick \(i \neq k\notin K\), and set \(N:=\top\setminus\{k\}\).
	We have
	\begin{align}
		\widetilde F_J(q_J,t_{K\setminus J}, 0_{\top\setminus K})
		 & =
		\widetilde F_J(q_J,t_{K\setminus J},0_k,0_{N\setminus K})
		\\
		 & =g_{JN}(q_J,t_{K\setminus J},0_{N\setminus K})
		 &                                                &
		\text{ by } \eqref{eq:facetInclusionExclusion}
		\\
		 & =g_{JK}(q_J,t_{K\setminus J})
		 &                                                &
		\text{ by } \eqref{eq:zeroingCoordinatesGenFunc}
	\end{align}
	which is \eqref{eq:gen_function_extension_zeroing}. Finally, suppose
	that either $t_{\top\setminus J}=0$ or
	\[
		(q_J,t_{\top\setminus J})\in\tns{C}{_0^1}(L_{J\top})
			\simeq\tns{C}{_0^1}(L_J)\times[0,\varepsilon)^{n+1-m}.
	\]
	By Definition~\ref{def:exactCollarModel}\eqref{itm:boundaryVanishing},
	every $g_{JK}(q_J,t_{K\setminus J})$ vanishes. Hence
	$\widetilde F_J|_{\tns{C}{_0^1}(L_{J\top})}\equiv 0$.
\end{proof}

\begin{prop}
	\label{prop:extension_of_generating_functions}
	After possibly shrinking $\varepsilon$, there exists a family of smooth functions
	\begin{equation}
		\left\{
		F_J: L_{J\top} \to \bb{R}
		\,\,|\,\,
		J \in \cc{F}_{\Lambda_i[n]}, \, i \in J
		\right\}
	\end{equation}
	such that:
	\begin{enumerate}
		\item \label{itm:boundaryVanishing2}
		      $F_J$ vanishes on $\tns{C}{_0^1}(L_{J\top})$ and the corner $L_J$.
		\item \label{itm:faceGenerating}
		      For all $ J \subsetneq K$ in $\cc{F}_{\Lambda_i[n]}$
		      \begin{equation}
			      F_J|_{L_{JK}}=g_{JK},
			      \quad \text{i.e.} \quad
			      F_J(q_J,t_{K\setminus J}, 0_{\top\setminus K}) = g_{JK}(q_J, t_{K \setminus J}).
		      \end{equation}
		\item \label{itm:0-cochain} For all $I\subsetneq J$, writing $t_{\top\setminus I}=(t_{\top\setminus J},t_{J\setminus I})$, we have
		      \begin{equation}
			      \label{eq:gen_function_extension_compatibility}
			      F_I(q_I,t_{\top\setminus I})
			      =
			      F_J(\lcol_{IJ}(q_I,t_{J\setminus I}), t_{\top \setminus J}) + g_{IJ}(q_{I},t_{J \setminus I}).
			      \tag{$ F_{IJ} $}
		      \end{equation}
	\end{enumerate}
\end{prop}

\begin{proof}
	For $|J|=1$ define $F_{J} = \widetilde{F}_J$ as in \eqref{eq:inclusionExclusion}. By Lemma \ref{lem:inclusionExclusion} it satisfies (1)-(2), and \eqref{eq:gen_function_extension_compatibility} holds trivially. Proceed by induction on $|J|$. Let $2 \leq m \leq n$. Suppose that for some $\varepsilon >0$ we have already constructed
	\begin{equation}
		F_{J}: L_{J\top} \to \bb{R}, \quad \quad |J| \leq m-1,
	\end{equation}
	satisfying \eqref{itm:boundaryVanishing2}-\eqref{itm:0-cochain}.
	Let $J \subsetneq \top$ with $|J|=m$.  Define an open subset:
	\begin{equation}
		\widehat U_0 = \bigcup_{I=J\setminus \{j\}}  L_{I\top} \, \, \subset L_{J\top},
	\end{equation}
	where $j\neq i$; see Figure \eqref{fig:union-of-collars}. For $q_J = \lcol_{IJ}(q_{I}, t_j)$, $ I= J\setminus \{j\}$, define $F^0_J: \widehat U_0 \to \bb{R} $ by
	\begin{equation}
		\label{eq:gen_ext_0}
		F_J^0(q_J,t_{\top\setminus J})
		=
		F_{I}(q_{I}, t_{\top\setminus J}, t_j)
		-
		g_{IJ}(q_{I}, t_j)
	\end{equation}

	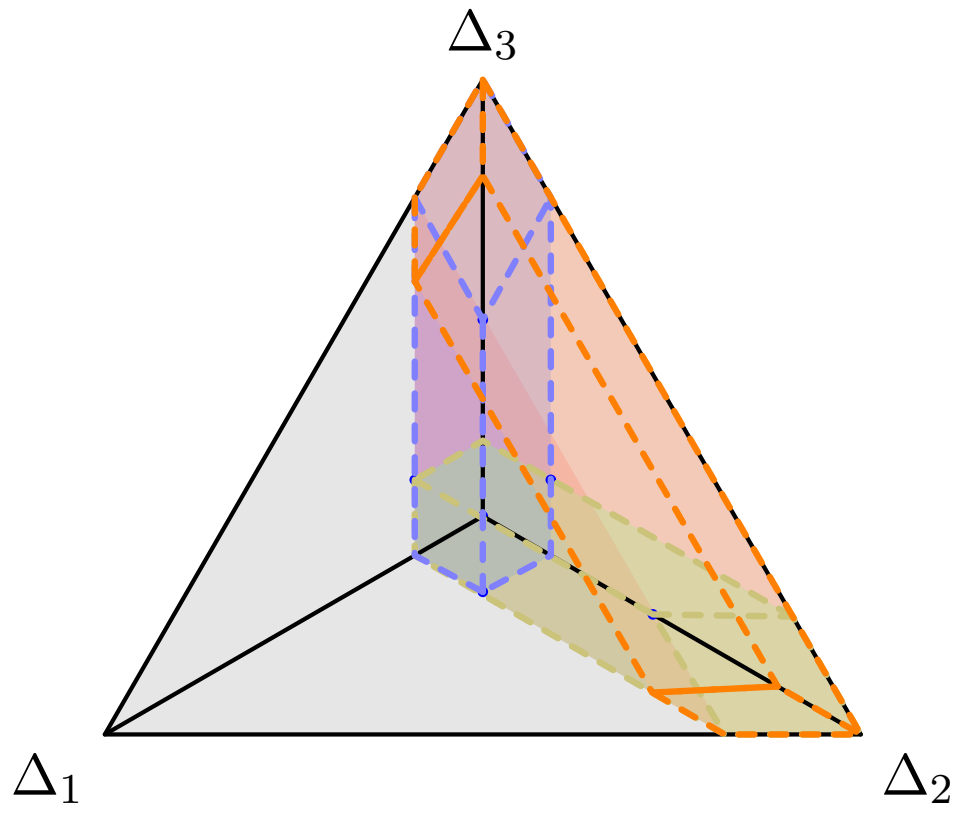
\begin{figure}[htbp]
		\centering
		\begin{tikzpicture}[scale=1.0, line join=round, line cap=round]
			\coordinate (A) at (0,2.887);
			\coordinate (B) at (-2.5,-1.443);
			\coordinate (C) at (2.5,-1.443);

			\def\t{0.18}

			\coordinate (AtoC) at ($(A)!\t!(C)$);
			\coordinate (BtoC) at ($(B)!\t!(C)$);
			\coordinate (BtoA) at ($(B)!\t!(A)$);
			\coordinate (CtoA) at ($(C)!\t!(A)$);
			\coordinate (CtoB) at ($(C)!\t!(B)$);
			\coordinate (AtoB) at ($(A)!\t!(B)$);

			\coordinate (center) at (barycentric cs:A=1,B=1,C=1);
			\coordinate (CentoB) at ($(center)!\t!(B)$);
			\coordinate (CentoC) at ($(center)!\t!(C)$);
			\coordinate (CentoB1) at ($(CentoB)+(0,0.5)$);
			\coordinate (CentoC1) at ($(CentoC)+(0,0.5)$);
			\coordinate (CentoC2) at ($(C)!0.55!(center)$);
			\coordinate (CentoA2) at ($(A)!0.55!(center)$);
			\coordinate (CentoA) at ($(center)+(0,0.5)$);
			\coordinate (CtoCent) at ($(center)!0.78!(C)$);
			\coordinate (AtoCent) at ($(center)!0.78!(A)$);
			\coordinate (AtoBLow) at ($(AtoB)!1.3*\t!(CentoB)$);
			\coordinate (CtoBLow) at ($(CtoB)!1.3*\t!(CentoB)$);
			\fill[gray!20] (A)--(B)--(C)--cycle;

			\fill[red!35, opacity=0.35] (C)--(A)--(AtoB)--(CtoB)--cycle;   
			\fill[red!50, opacity=0.5] (AtoB)--(CentoB)--(CtoB)--cycle;
			\fill[blue!50, opacity=0.3] (AtoB)--(CentoB)--(0,-0.5)--(CentoC)--(AtoC)--(A)--cycle;
			\fill[green!50, opacity=0.3] (CtoB)--(CentoB)--(CentoB1)--(CentoA)--(CtoA)--(C)--cycle;
			\fill[orange!50, opacity=0.3] (A)--(AtoB)--(AtoBLow)--(CtoBLow)--(CtoB)--(C)--cycle;

			\draw[thick] (A)--(B)--(C)--cycle;
			\node[above] at (A) {$\Delta_3$};
			\node[below left] at (B) {$\Delta_1$};
			\node[below right] at (C) {$\Delta_2$};

			\fill[blue] (center) circle (1pt);
			\draw[thick] (A) -- (center);
			\draw[thick] (B) -- (center);
			\draw[thick] (C) -- (center);
			\fill[blue] (0,-0.5) circle (1pt);
			\fill[blue] (CentoB1) circle (1pt);
			\fill[blue] (CentoC1) circle (1pt);
			\fill[blue] (CentoC2) circle (1pt);
			\fill[blue] (CentoA2) circle (1pt);

			\draw[dashed, very thick, olive!50] (CentoB1)--(CentoC2)--(CtoB)--(CentoB)--cycle;
			\draw[dashed, very thick, olive!50] (CentoC2)--(CtoA)--(C)--(CtoB);
			\draw[dashed, very thick, olive!50] (CtoA)--(CentoA)--(CentoB1);

			\draw[dashed, very thick, blue!50] (0,-0.5)--(CentoB)--(AtoB)--(A)--(AtoC)--(CentoC)--cycle;
			\draw[dashed, very thick, blue!50] (AtoB) -- (CentoA2) -- (AtoC);
			\draw[dashed, very thick, blue!50] (0,-0.5) -- (CentoA2);

			\draw[dashed, very thick, orange!] (AtoBLow)--(AtoCent)--(CtoCent)--(CtoBLow)--cycle;
			\draw[dashed, very thick, orange!] (AtoB)--(AtoBLow)--(AtoCent)--(A)--cycle;
			\draw[dashed, very thick, orange!] (CtoCent)--(CtoBLow)--(CtoB)--(C)--cycle;
			\draw[dashed, very thick, orange!] (A)--(C);

		\end{tikzpicture}
		\caption{The 3-simplex with $\top = \{0123\}$. The subset $\widehat{U}_0$ of the collar $\widetilde\Delta_{023,\top}$ is the union of the collars $\widetilde\Delta_{02,\top}$ (green), $\widetilde\Delta_{03,\top}$ (blue), and $\widetilde\Delta_{23,\top}$ (orange). }
		\label{fig:union-of-collars}
	\end{figure}

	We first show $ F_J^0$ is well defined. Let $j_1,j_2 \in J$ and set $I_1 = J \setminus \{j_1\}$, $I_2 = J \setminus \{j_2\}$ and $H = J \setminus \{j_1,j_2\}$.
	Let $q_J \in  {L}_{I_1J} \cap  {L}_{I_2J} = {L}_{HJ}$. There exist
	\begin{equation}
		(q_{I_1}, t_{j_1}) \in L_{I_1J},
		\quad \quad
		(q_{I_2}, t_{j_2}) \in L_{I_2J},
		\quad \quad
		(q_{H}, t_{j_1}, t_{j_2}) \in L_{HJ},
	\end{equation}
	such that
	\begin{equation}
		\begin{gathered}
			\lcol_{I_1J}(q_{I_1},t_{j_1})=q_J=\lcol_{I_2J}(q_{I_2},t_{j_2}),\\
			\lcol_{HI_1}(q_H,t_{j_2})=q_{I_1},\qquad
			\lcol_{HI_2}(q_H,t_{j_1})=q_{J_2}.
		\end{gathered}
	\end{equation}
	We get for all $t_J$:
	\begin{align}
		 &
		F_{I_1}(q_{I_1}, t_{\top\setminus J},t_{j_1})
		-
		g_{I_1J}(q_{I_1}, t_{j_1})
		\\
		 & =
		F_{H}(q_{H}, t_{\top\setminus J},t_{j_1},t_{j_2})
		-
		g_{HI_1}(q_{H}, t_{j_2})
		-
		g_{I_1J}(q_{I_1}, t_{j_1})
		 &   &
		\text{ by } \hyperref[eq:gen_function_extension_compatibility]{(F_{H I_1})}
		\\
		 & =
		F_{H}(q_{H}, t_{\top\setminus J},t_{j_1},t_{j_2})
		-
		g_{HJ}(q_{H},t_{j_1},t_{j_2})
		 &   &
		\text{ by } \hyperref[eq:compatible_generating_functions]{(G_{H I_1 J})}
		\\
		 & =
		F_{H}(q_{H}, t_{\top\setminus J},t_{j_1},t_{j_2})
		-
		g_{HI_2}(q_{H}, t_{j_1})
		-
		g_{I_2J}(q_{I_2}, t_{j_2})
		 &   &
		\text{ by } \hyperref[eq:compatible_generating_functions]{(G_{H I_2 J})}
		\\
		 & =
		F_{I_2}(q_{I_2}, t_{\top\setminus J},t_{j_2})
		-
		g_{I_2J}(q_{I_2}, t_{j_2}).
		 &   &
		\text{ by }\hyperref[eq:gen_function_extension_compatibility]{(F_{H I_2})}
	\end{align}
	Hence the value of $ F^0_J(q_J,t_{\top\setminus J})$ defined by either $ I_1 $ or $ I_2 $ agrees, so $ F^0_J $ is well defined.

	Choose $0< \varepsilon_J < \varepsilon$ and open subset $\widehat U_1 \subset  L_{J\top}$ such that
	\begin{equation}
		\widehat U_0 \cup \widehat U_1 =  L_{J\top} ,
		\quad \quad
		\bigcup_{I = J\setminus \{j\}}
		\left(
		L_{I\top}
		\left( \varepsilon_J \right)
		\right)
		\subseteq
		\widehat U_0 \setminus \widehat U_1.
	\end{equation}
	Let $\psi_{0}, \psi_1$ be a partition of unity for $\widehat U_0$, $\widehat U_1$, and define
	\begin{equation}
		F_J: L_{J\top} \to \bb{R},
		\quad \quad
		F_J
		=
		\psi_0 \cdot F^0_{J} +  \psi_1 \cdot F^1_{J},
	\end{equation}
	where $F^1_J = \widetilde{F}_J: L_{J \top} \to \bb{R}$ is defined by   \eqref{eq:inclusionExclusion}.

	We now prove that $F_J$ satisfies property \eqref{itm:faceGenerating}.
	Fix $ j\in J $ and set $ I = J \setminus \{j\} $. For $q_I \in  L_{IJ}$ and for every $J \subsetneq K $ we have
	\begin{align}
		\MoveEqLeft[1] F^0_J(\lcol_{IJ}(q_{I}, t_j), t_{K \setminus J},0_{\top\setminus K}) \notag\\
		 & =
		F_{I}(q_{I}, t_{K \setminus J},t_j,0_{\top\setminus K})
		-
		g_{IJ}(q_{I}, t_j)
		 &                                                     &
		\text{by }\eqref{eq:gen_ext_0}
		\\
		 & =
		g_{IK}(q_{I}, t_{K \setminus J},t_j)
		-
		g_{IJ}(q_{I}, t_j)
		 &                                                     &
		\text{by property } \eqref{itm:faceGenerating}
		\\
		 & = g_{JK}(\lcol_{IJ}(q_{I}, t_j), t_{K \setminus J})
		 &                                                     &
		\text{by } \hyperref[eq:compatible_generating_functions]{(G_{IJK})}.
		\label{eq:functionsZeroingcollar}
	\end{align}
	For every $J \subsetneq K$, property \eqref{itm:faceGenerating} holds for $ F_{J} $ because $F_{J}^0$  satisfies \eqref{eq:functionsZeroingcollar} and $F_{J}^1$ satisfies \eqref{eq:gen_function_extension_zeroing}.

	Next, we prove property \eqref{itm:boundaryVanishing2} for $F_J$. By Lemma \ref{lem:inclusionExclusion}, $F_J^1$ vanishes on ${\tns{C}{_0^1}(L_{J\top})}$ and on the corner $L_J$. For $F_{J}^0$ let $(q_J,t_{\top\setminus J})\in \tns{C}{_0^1}(L_{J\top})\simeq \tns{C}{_0^1}(L_J)\times [0,\varepsilon)^{n+1-m}$, such that $q_J= \lcol_{IJ}(q_I,t_j)$. Then $q_I \in \tns{C}{_0^1}(L_I)$, and
	\begin{equation}
		F_J^0(q_J,t_{\top\setminus J})
		=
		F_{I}(q_{I}, t_{\top\setminus J}, t_j)
		-
		g_{IJ}(q_{I}, t_j)
		=
		0-0=0,
	\end{equation}
	by the induction hypothesis for $F_I$ and Definition \ref{def:exactCollarModel}\eqref{itm:boundaryVanishing} for $g_{IJ}$.
	Finally, for $q_J = \lcol_{IJ}(q_I,t_j)$, property \eqref{itm:faceGenerating} of $F_I$ gives
	\begin{equation}
		F_J^0(q_J,0_{\top\setminus J})
		=
		F_{I}(q_{I}, 0_{\top\setminus J}, t_j)
		-
		g_{IJ}(q_{I}, t_j)
		=
		g_{IJ}(q_{I}, t_j) - g_{IJ}(q_{I}, t_j)
		=0.
	\end{equation}
	Thus, both $F_J^0$ and $F_J^1$ vanish on $\tns{C}{_0^1}(L_{J\top})$ and on $L_J$, so $F_J=\psi_0F_J^0+\psi_1F_J^1$ satisfies~\eqref{itm:boundaryVanishing2}.

	Now let $\emptyset \neq H \subsetneq J $, pick $j \in J \setminus H$ and set $I = J \setminus	 \{j\}$.
	Let $(q_H, t_{\top \setminus H}) \in L_{H\top}(\varepsilon_J)$. Then
	\begin{equation}
		(\lcol_{HI}(q_H,t_{I \setminus H}), t_j,t_{\top \setminus J}) \in L_{I \top }(\varepsilon_J),
	\end{equation}
	and by \hyperref[eq:collars_relation_2]{$(\widehat{R}_{HIJ})$},
	\begin{equation}
		(\lcol_{HJ}(q_H,t_{J\setminus H}), t_{\top\setminus J})
		=
		(\lcol_{IJ}(\lcol_{HI}(q_{H}, t_{I\setminus H}),t_j), t_{\top\setminus J}) \in
		\widehat{U}_0 \setminus \widehat{U}_1
		\subseteq L_{J \top}
	\end{equation}
	which is the region where $ \psi_0 \equiv 1, \psi_1 \equiv 0 $.
	Thus,
	\begin{align}
		 & F_{J}(\lcol_{HJ}(q_H,t_{J\setminus H}), t_{\top\setminus J}) + g_{HJ}(q_{H},t_{J \setminus H})
		\\
		 & =
		F_{J}^0(\lcol_{IJ}(\lcol_{HI}(q_{H}, t_{I\setminus H}),t_j), t_{\top\setminus J})
		 &                                                                                                &
		\text{by } \hyperref[eq:collars_relation_2]{(\widehat{R}_{HIJ})}
		\\
		 & \quad \quad +
		g_{IJ}(\lcol_{HI}(q_H,t_{I\setminus H}),t_j)
		+
		g_{HI}(q_H, t_{I\setminus H})
		 &                                                                                                &
		\text{by }\hyperref[eq:compatible_generating_functions]{(G_{HIJ})}
		\\
		 & =
		F_{I}( \lcol_{HI}(q_{H}, t_{I\setminus H}), t_{\top\setminus J},t_j)
		- g_{IJ}(\lcol_{HI}(q_H,t_{I\setminus H}),t_j)
		 &                                                                                                &
		\text{by }\eqref{eq:gen_ext_0}
		\\
		 & \quad \quad +
		g_{IJ}(\lcol_{HI}(q_H,t_{I\setminus H}),t_j)
		+
		g_{HI}(q_H, t_{I\setminus H})
		\\
		 & =
		F_{I}( \lcol_{HI}(q_{H}, t_{I\setminus H}), t_{\top\setminus J},t_j) +
		g_{HI}(q_H, t_{I\setminus H})
		\\
		 & = F_{H}(q_H, t_{\top\setminus H}),
		 &                                                                                                &
		\text{by }\hyperref[eq:gen_function_extension_compatibility]{(F_{HI})}
	\end{align}
	which proves property \eqref{itm:0-cochain} for $H,J$ on $L_{H \top}(\varepsilon_J)$. Finally, shrink $\varepsilon$ to  
    \[
        \min \{\varepsilon_J \,\,|\,\, |J|=m\}
    \] 
    This yields functions $F_{J}$ for all $|J| \leq m$ satisfying \eqref{itm:boundaryVanishing2}-\eqref{itm:0-cochain}, completing the induction step.
\end{proof}

\subsection{Extension of Lagrangian horns}
\label{sec:ext_Lag_horns}

Fix a Lagrangian $W$ over ${\Lambda_i[n]}$ and write $\top = [n]$. For $J \in \cc{F}_{\Lambda_i[n]}$, let $\vp_J:V_J \to X_{\Delta_J}$ be a Weinstein chart and let $\coll_{J\top}: \Delta_{J\top} \to \widetilde{\Delta}_{J\top}$ be a collar map as in Theorem~\ref{thm:compatible_exact_models}. Let $F_J: L_{J\top} \to \bb{R}$ be a smooth function as in Proposition \ref{prop:extension_of_generating_functions}.

\begin{definition}
	Let $\widehat\Lambda$ be the space obtained by gluing of the facet collars along their codimension-2 collars:
	\begin{equation}
		\widehat\Lambda := \bigsqcup_{\substack{K \in \cc{F}_{\Lambda_i[n]}\\ |K|=n}} \Delta_{K\top} \bigg/ \sim
	\end{equation}
	where for $K_1 = \top \setminus \{i_1\}$, $K_2 = \top \setminus \{i_2\}$ and $J = \top \setminus \{i_1,i_2\}$, we declare  $(s_{K_1},t_{i_1}) \sim (s_{K_2},t_{i_2})$ if and only if there exists $s_{J} \in \Delta_J$ such that
	\begin{equation}
		\coll^{J \top}_{K_1\top}(s_J,t_{i_1},t_{i_2}) = (s_{K_1},t_{i_1}) , \quad \quad
		\coll^{J \top}_{K_2\top}(s_J,t_{i_1},t_{i_2}) = (s_{K_2},t_{i_2}).
	\end{equation}
	Similarly, define $\widehat L$ by gluing of the Lagrangian facet collars along their codimension-2 Lagrangian collars:
	\begin{equation}
		\widehat L := \bigsqcup_{\substack{K \in \cc{F}_{\Lambda_i[n]}\\ |K|=n}} L_{K\top} \bigg/ \sim
	\end{equation}
	where $(q_{K_1},t_{i_1}) \sim (q_{K_2},t_{i_2})$ if and only if there exists $q_J \in L_J$ such that
	\begin{equation}
		\lcol^{J \top}_{K_1\top}(q_J,t_{i_1},t_{i_2}) = (q_{K_1},t_{i_1}) , \quad \quad
		\lcol^{J \top}_{K_2\top}(q_J,t_{i_1},t_{i_2}) = (q_{K_2},t_{i_2}).
	\end{equation}
\end{definition}

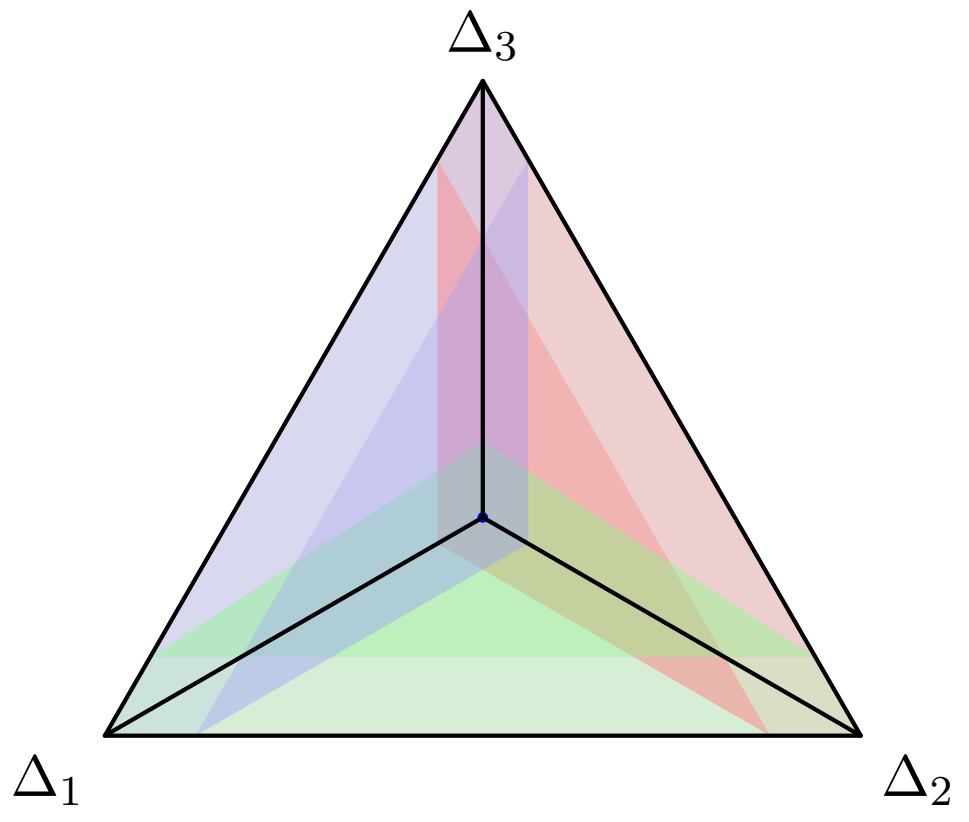
\begin{figure}[htbp]
	\centering
	\begin{tikzpicture}[scale=1.0, line join=round, line cap=round]
		\coordinate (A) at (0,2.887);
		\coordinate (B) at (-2.5,-1.443);
		\coordinate (C) at (2.5,-1.443);

		\def\t{0.12}

		\coordinate (AtoC) at ($(A)!\t!(C)$);
		\coordinate (BtoC) at ($(B)!\t!(C)$);
		\coordinate (BtoA) at ($(B)!\t!(A)$);
		\coordinate (CtoA) at ($(C)!\t!(A)$);
		\coordinate (CtoB) at ($(C)!\t!(B)$);
		\coordinate (AtoB) at ($(A)!\t!(B)$);

		\coordinate (center) at (barycentric cs:A=1,B=1,C=1);
		\coordinate (CentoB) at ($(center)!\t!(B)$);
		\coordinate (CentoC) at ($(center)!\t!(C)$);
		\coordinate (CentoB1) at ($(CentoB)+(0,0.5)$);
		\coordinate (CentoC1) at ($(CentoC)+(0,0.5)$);
		\coordinate (CentoC2) at ($(C)!0.55!(center)$);
		\coordinate (CentoA2) at ($(A)!0.55!(center)$);
		\coordinate (CentoA) at ($(center)+(0,0.5)$);
		\fill[gray!20] (A)--(B)--(C)--cycle;

		\fill[red!35, opacity=0.35] (C)--(A)--(AtoB)--(CtoB)--cycle;   
		\fill[blue!25, opacity=0.35] (A)--(B)--(BtoC)--(AtoC)--cycle;  
		\fill[green!30, opacity=0.3] (B)--(C)--(CtoA)--(BtoA)--cycle; 
		\fill[red!50, opacity=0.5] (AtoB)--(CentoB)--(CtoB)--cycle;
		\fill[green!50, opacity=0.4] (BtoA)--(CentoA)--(CtoA)--cycle;
		\fill[blue!40, opacity=0.4] (BtoC)--(CentoC)--(AtoC)--cycle;
		\draw[thick] (A)--(B)--(C)--cycle;
		\node[above] at (A) {$\Delta_3$};
		\node[below left] at (B) {$\Delta_1$};
		\node[below right] at (C) {$\Delta_2$};

		\fill[blue] (center) circle (1pt);
		\draw[thick] (A) -- (center);
		\draw[thick] (B) -- (center);
		\draw[thick] (C) -- (center);
	\end{tikzpicture}
	\caption{The space $\widehat\Lambda$ obtained by gluing the facet collars $\Delta_{023,\top}$ (in red), $\Delta_{013,\top}$ (in blue) and $\Delta_{012,\top}$ (in green).}
\end{figure}

Recall Notation \ref{ntn:graphs}.
\begin{lemma}
	\label{lem:globalLag}
	The Lagrangian immersions
	\begin{equation}
		\vp_{K\top} \circ \Gamma F_K :\
		L_{K\top}
		\to
		X_{\Delta_{K\top}},
		\quad \quad
		K \in \cc{F}_{\Lambda_i[n]}, \,\, |K|=n,
	\end{equation}
	glue to a global immersed Lagrangian  $\widehat{W}  = \left[\widehat{\Psi} : \widehat{L} \to X_{\widehat{\Lambda}} \right]$ such that $\widehat{W}|_{\Delta_K} = W_K$.
\end{lemma}
\begin{proof}
	Define $\widehat{\Psi}|_{L_{K\top}} := \vp_{K\top} \circ \Gamma F_K$. To show it is well defined, it suffices to verify that for every $J \subsetneq K$, the following diagram commutes:
	\begin{equation}
		\begin{tikzcd}
			{L_{J\top} }
			\arrow[d,"\lcol_{K\top}^{J\top}" left]
			\arrow[rr, "\vp_{J\top} \circ \Gamma F_{J}" ]
			&&
			X_{\Delta_{J\top}}
			\arrow[d, "\widehat{\coll}_{K\top}^{J\top}"]
			\\
			L_{K\top}
			\arrow[rr, "\vp_{K\top} \circ \Gamma F_{K}" below]
			&&
			X_{\Delta_{K\top}}
		\end{tikzcd}
	\end{equation}
	Indeed,
	\begin{align}
		\MoveEqLeft[1] (\vp_{K\top} \circ \Gamma F_{K}) \circ \lcol_{K\top}^{J\top} \notag\\
		 & = \widehat{\coll_{K\top}^{J\top}} \circ \vp_{K\top}|_{L_{J\top}} \circ \Gamma(F_{K} \circ \lcol^{J\top}_{K\top})
		 &                                                                                                                                                 &
		\text{ by } \eqref{eq:extended_compatible_Weinstein_symplecto2}
		\\
		 & = \widehat{\coll_{K\top}^{J\top}} \circ\vp_{J\top} \circ \cc{T}(g_{JK} \circ \pi_{JK}^{J\top}) \circ  \Gamma (F_{K}\circ \lcol_{K\top}^{J\top})
		 &                                                                                                                                                 &
		\text{ by }\eqref{eq:extended_compatible_Weinstein_symplecto}
		\\
		 & = \widehat{\coll_{K\top}^{J\top}} \circ\vp_{J\top} \circ \Gamma (F_{K}\circ \lcol_{K\top}^{J\top} + g_{JK} \circ \pi_{JK}^{J\top})
		\\
		 & =\widehat{\coll_{K\top}^{J\top}} \circ \vp_{J\top} \circ \Gamma F_{J}
		 &                                                                                                                                                 &
		\text{ by } \hyperref[eq:gen_function_extension_compatibility]{(F_{JK})}.
	\end{align}
	Denote by $W_{K\top}$ the immersed Lagrangian represented by
	\begin{equation}
		\vp_{K\top} \circ \Gamma F_K :\
		L_{K\top}
		\to
		X_{\Delta_{K\top}}.
	\end{equation}
	Let $\iota: \Delta_{K \top} \to \widehat\Lambda$ be the inclusion. After shrinking $\varepsilon$ if necessary, continuity and compactness ensure that any point of $\widehat{L}$ projecting into $\Delta_{K\top}$ from a different face collar lies in the original overlap with $L_{K\top}$ before the shrinking. Since
    \[
        \tau_{K\top}\circ\pi_{\Delta_{K\top}}\circ
        \vp_{K\top}\circ\Gamma F_K=\operatorname{pr}_2,
    \]
    the representative in the original $L_{K\top}$ has collar coordinate in $[0,\varepsilon)$. Thus the overlap identification already takes place at the shrinked collar, giving $L_{K\top}\simeq\iota^*\widehat L$.
    By Lemma \ref{lem:restriction_to_open_subset}, the commutative diagram
	\begin{equation}
		\begin{tikzcd}
			{L_{K\top}}
			\arrow[d, hook']
			\arrow[rr, "\vp_{K\top} \circ \Gamma F_{K}" ]
			&&
			X_{\Delta_{K\top}}
			\arrow[d, hook']
			\\
			\widehat{L}
			\arrow[rr, "\widehat\Psi"]
			&&
			X_{\widehat\Lambda}
		\end{tikzcd}
	\end{equation}
	implies that $\iota^* \widehat{\Psi} = \vp_{K\top} \circ \Gamma F_{K}$, so $\widehat{W}|_{\Delta_{K\top}} = W_{K \top}$. Let $i: \Delta_K \hookrightarrow \Delta_{K\top}$ be the face inclusion. By Proposition \ref{prop:extension_of_generating_functions}\eqref{itm:boundaryVanishing2}, $F_K$ vanishes on $L_K$, and by Definition \ref{def:WeinsteinModel}\eqref{itm:zero-section}, $\vp_K|_{L_K} = \Psi_K$. Thus, the diagram
	\begin{equation}
		\begin{tikzcd}
			L_K
			\arrow[d, hook]
			\arrow[drr, "\Psi_K"]
			\\
			L_{K\top}
			\arrow[r, "\vp_{K\top} \circ \Gamma F_{K}" below]
			&
			X_{\Delta_{K\top}}
			\arrow[r, "\mathrm{pr}_{X_{\Delta_K}}" below]
			&
			X_{\Delta_K}
		\end{tikzcd}
	\end{equation}
	commutes. It follows by Lemma \ref{lem:cobordismOverProduct} that, after identifying $i^*L_{K \top} \simeq L_K$, we have $i^*(\vp_{K\top} \circ \Gamma F_{K}) = \Psi_K$, so $W_{K \top}|_{\Delta_K} = W_K$. Finally,
	\begin{equation}
		\widehat{W}|_{\Delta_{K}}
		=
		(\widehat{W}|_{\Delta_{K\top}})|_{\Delta_K}
		=
		W_{K \top}|_{\Delta_K}
		=
		W_K,
	\end{equation}
	which proves the lemma.
\end{proof}

\begin{corollary}
	\label{cor:globalLag}
	$\widehat{W}$ is a Lagrangian cobordism of type $\widehat{L}$ over $\widehat{\Lambda}$ relative to $\vlg$.
\end{corollary}

\begin{proof}
	By Lemma \ref{lem:globalLag} we have a global immersion $\widehat{\Psi}:\widehat{L}\to X_{\widehat{\Lambda}}$ whose restriction to each facet collar $L_{KT}$ is $\varphi_{KT}\circ \Gamma F_K$.
	By Proposition \ref{prop:extension_of_generating_functions}\eqref{itm:boundaryVanishing2} the function $F_K$ satisfies $F_K|_{C_0^1(L_{K\top})}\equiv 0$, hence Lemma \ref{lem:exact-perturbations}
	implies that $\varphi_{K\top}\circ \Gamma F_K$ is a Lagrangian immersion representing a
	Lagrangian cobordism over $\Delta_{K\top}$ relative to $\vlg$. In particular, it
	satisfies conditions \eqref{itm:bsubmersion}-\eqref{itm:verticalCorners} of Definition \ref{def:LagrangiankCobordism}.
	Since these conditions are
	local on the source and target,
	$\widehat{\Psi}$ satisfies Definition \eqref{def:LagrangiankCobordism} globally.
\end{proof}

\subsection{The Kan complex of Lagrangian cobordisms}
\label{subsec:kan-complex}
In this Section we verify that $\cob(X,\vlg)$ is a simplicial set
and prove Theorem~\ref{thm:hornFilling}.

\begin{lemma}
	\label{lem:fixed-locus}
	Let $\Psi:L\to M$ be an immersion of manifolds with corners and
	let $\varphi:L\to L$ be a diffeomorphism satisfying
	$\Psi\circ\varphi=\Psi$. Then the fixed locus of $\varphi$ is
	a union of connected components of $L$. In particular, $\Psi$ is
	free if every such $\varphi$ fixes a point on each connected
	component of $L$.
\end{lemma}

\begin{proof}
	The fixed locus is closed. If $\varphi(q)=q$, choose an open
	neighborhood $U$ of $q$ on which $\Psi$ is injective. For every
	$p\in U\cap\varphi^{-1}(U)$, the equality
	$\Psi(\varphi(p))=\Psi(p)$ gives $\varphi(p)=p$. Thus the fixed
	locus is also open.
\end{proof}

\begin{proof}[Proof of Lemma~\ref{lem:simplicial-set-cob}]
	The simplicial identities follow from functoriality of the
	cotangent pullback formula~\eqref{eq:pullback}, as in
	Corollary~\ref{cor:functorialPullback}. The face maps give compact freely immersed Lagrangians, since faces are closed and every face of a face is again
	a face.

	Let $W=[\Psi:L\to X_{\Delta^n}]\in\cob(X,\vlg)_n$ and put
	$\sigma=\bar\sigma_{n,i}$. The surjection $\sigma$ is weakly smooth and has surjective
	differential. Its restriction to the interior of each face
	is a submersion on the interior of its image face.
	Thus $\sigma^\ast W$ is a Lagrangian cobordism by
	Lemma~\ref{lem:weakly-smooth-pullback}, of
	topological type
	\[
		P=L\times_{T^\ast\Delta^n}\sigma^\ast T^\ast\Delta^n
		\cong L\times_{\Delta^n}\Delta^{n+1},
	\]
	which is compact since it is
	closed in $L\times\Delta^{n+1}$. It remains to check freeness
	on every nonempty face.
	Let $E\subseteq\Delta^{n+1}$ be a nonempty face, put
	$F=\sigma(E)$ and write $\tau=\sigma|_E$. Then
	$(\sigma^\ast W)|_E=\tau^\ast(W|_F)$. Choose a face $D$ of $E$
	on which $\tau$ is an isomorphism onto $F$. We may take $D=E$
	unless $E$ contains both vertices $i,i+1$, in which case we
	omit vertex $i$. Writing $W|_F=[\Psi_F:L_F\to X_F]$, any
	diffeomorphism preserving $\tau^\ast\Psi_F$ preserves the
	$E$-coordinate on $L_F\times_F E$. Over $D$ it therefore
	induces a diffeomorphism preserving $\Psi_F$, and hence fixes
	every point there by freeness of $W|_F$. Each fiber of
	$L_F\times_F E\to L_F$ is a point or an interval meeting the
	part over $D$. Lemma~\ref{lem:fixed-locus} now gives freeness
	on $E$, including when $E=\Delta^{n+1}$.
\end{proof}

\begin{lemma}
	\label{lem:compact-horn-extension}
	Let $n\geq2$ and let $W:\Lambda_i[n]\to\cob(X,\vlg)$ be a
	Lagrangian over the $(n,i)$-horn. There exists a compact
	Lagrangian $n$-cobordism $\widetilde W\subset X_{\Delta^n}$
	relative to $\vlg$ such that
	\[
		\widetilde W|_{\Delta_K}=W_K
		\qquad\text{for every }K\in\cc F_{\Lambda_i[n]}.
	\]
\end{lemma}

\begin{figure}[htbp]
	\centering
	\begin{tikzpicture}[scale=1.0, line join=round, line cap=round]
		\coordinate (A) at (0,2.887);
		\coordinate (B) at (-2.5,-1.443);
		\coordinate (C) at (2.5,-1.443);

		\def\t{0.18}

		\coordinate (AtoC) at ($(A)!\t!(C)$);
		\coordinate (BtoC) at ($(B)!\t!(C)$);
		\coordinate (BtoA) at ($(B)!\t!(A)$);
		\coordinate (CtoA) at ($(C)!\t!(A)$);
		\coordinate (CtoB) at ($(C)!\t!(B)$);
		\coordinate (AtoB) at ($(A)!\t!(B)$);
		\coordinate (BtoCent) at ($(B)!0.4!(0,0)$);
		\coordinate (CtoAshift) at ($(CtoA)+(-0.4,-0.1)$);
		\coordinate (AtoCshift) at ($(AtoC)+(-0.4,0)$);


		\fill[blue!25, opacity=0.35] (A)--(B)--(BtoC)--(AtoC)--cycle;  
		\fill[green!30, opacity=0.35] (B)--(C)--(CtoA)--(BtoA)--cycle; 

		\draw[thick] (A)--(B)--(C);
		\draw (C)--(CtoA);
		\draw (A)--(AtoC);
		\node[above] at (A) {$\Delta_2$};
		\node[below left] at (B) {$\Delta_0$};
		\node[below right] at (C) {$\Delta_1$};

		\node[left] at ($(A)!0.5!(B)$) {$\Delta_{02}$};
		\node[below] at ($(B)!0.5!(C)$) {$\Delta_{01}$};

		\draw[thick, red, smooth] plot coordinates {(C) (CtoAshift) (BtoCent) (AtoCshift) (A)};

		\fill[pattern=dots, pattern color=red, opacity=0.5]
		(B) -- (C)
		plot[smooth] coordinates {(C) (CtoAshift) (BtoCent) (AtoCshift) (A)}
		-- (B) -- cycle;

	\end{tikzpicture}
	\caption{The neighborhood $U$ of $\Lambda_0^2 \subset \widehat \Lambda $ (dotted in red).}
	\label{fig: hornExtension}
\end{figure}
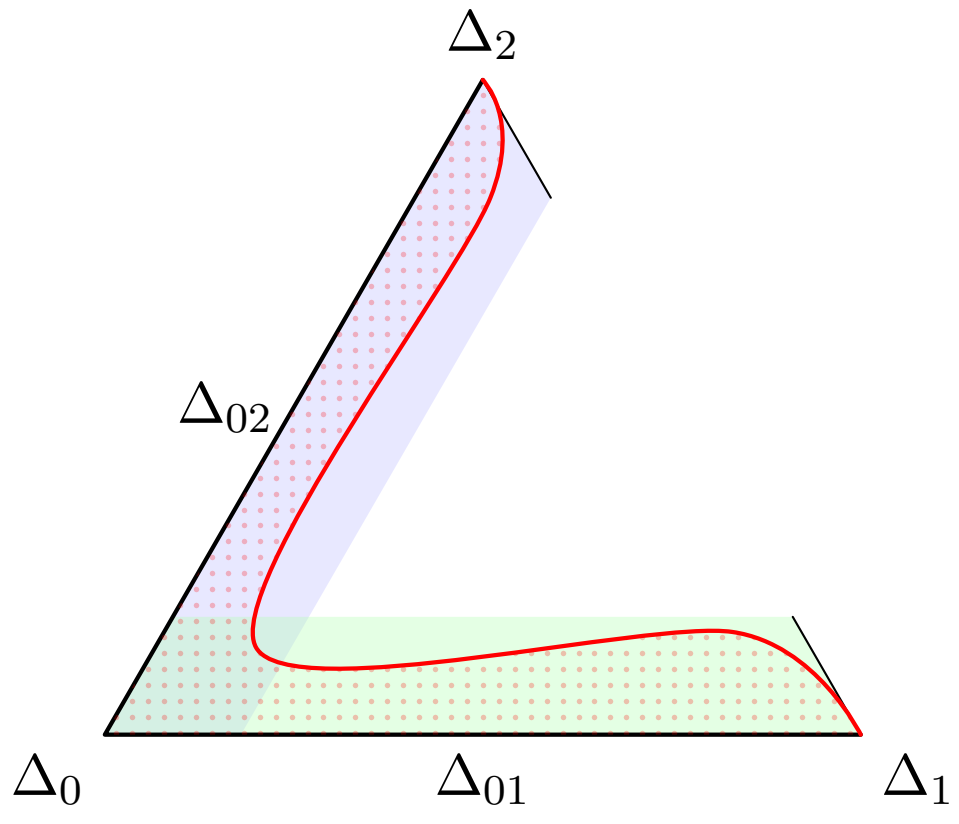

\begin{proof}
	Let
	\begin{equation}
		\Lambda_i^n := \left\{ t\in \Delta^n \mid t_j =0 \text{  for at least one } j \neq i \right\}
	\end{equation}
	be the geometric $(n,i)$-horn in the geometric $n$-simplex. Let $\widehat W$ as in Corollary \ref{cor:globalLag}.
    Choose an open neighborhood $U$ of $\Lambda_i^n$ with $\overline U\subset\widehat\Lambda$ as in Figure~\ref{fig: hornExtension}. Obtain $\overline U$ by sliding the missing facet towards the $i$-th vertex along the segments joining it to that vertex. Choose a smooth displacement that vanishes near the boundary of the missing facet. Reparametrizing these segments gives a diffeomorphism $\phi: \Delta^n \to \overline{U}$ relative to $\Lambda_i^n$.
    Denote the new facet of $\overline U$ by $H$. By a small ambient perturbation fixing $\partial\Delta^n$, we may assume that $H\cap\mathrm{int}(\Delta^n)$ is transverse to the restriction of $\widehat\Psi_{\widehat\Lambda}$ to each stratum $S^k(\widehat L)$, $k\geq 0$.

	The inclusion $\iota:\overline U\hookrightarrow\widehat\Lambda$ is
	weakly smooth in the sense of Definition~\ref{def:smoothMaps}, and its differential is invertible. Near the horn it is locally
	the identity. Wherever $H$ meets the interior of the missing facet
	$\{t_i=0\}$, the two facets have the same tangent space. Thus
	\eqref{eq:transversality-tangent-spaces1} and
	\eqref{eq:transversality-tangent-spaces2} hold, and
	$\iota^\ast\widehat W$ is a Lagrangian cobordism by
	Lemma~\ref{lem:weakly-smooth-pullback}.
	Then
	\begin{equation}
		\widetilde{W} = \phi^\ast \iota^\ast \widehat{W}
	\end{equation}
	is a Lagrangian cobordism over $\Delta^n$, relative to $\vlg$. It is compact since $\widehat\Psi_{\widehat\Lambda}^{-1}(\overline U)$
	is closed in a finite union of compact facet subcollars.
	Since $\phi$ preserves $\Lambda_i^n$, we have $\widetilde{W}|_{\Delta_K} = \widehat{W}|_{\Delta_K} = W_K$ by Lemma \ref{lem:globalLag}.
\end{proof}

\begin{lemma}
	\label{lem:free-horn-extension}
	The cobordism $\widetilde W$ in
	Lemma~\ref{lem:compact-horn-extension} can be chosen so that
	$\widetilde W|_F$ is free for every nonempty face
	$F\subseteq\Delta^n$, including $\Delta^n$ itself.
\end{lemma}

\begin{proof}
	Write $\widetilde W=[\widetilde\Psi:\widetilde L\to X_{\Delta^n}]$
	and retain only the connected components of $\widetilde L$ whose
	projections meet the horn. Any diffeomorphism preserving
	$\widetilde\Psi$ induces a diffeomorphism preserving the pullback
	immersion on each prescribed facet. Since these immersions are
	free, it fixes every point over the horn. Each retained component
	meets this set, so Lemma~\ref{lem:fixed-locus} gives freeness on
	$\Delta^n$. This argument also applies after the modification
	below, which is supported away from the horn.

	Put $J=\top\setminus\{i\}$. Every nonempty proper face other than
	$\Delta_J$ belongs to the horn and is free by assumption. Let
	$\Psi_J:L_J\to X_{\Delta_J}$ be the pullback immersion on
	$\Delta_J$. Since every proper face of $\Delta_J$ belongs to the
	horn, the same argument shows that any diffeomorphism preserving
	$\Psi_J$ fixes pointwise each component of $L_J$ whose projection
	meets $\partial\Delta_J$. We give each remaining component a
	singly covered point.

	Choose an interior point on each remaining component, with
	distinct images. At each chosen point $q$, take a Weinstein
	chart $\vp_J$ for an embedded restriction $\Psi_J|_U$, with the
	chart images pairwise disjoint. Choose $(x,\xi)$ arbitrarily
	close to $(q,0)$ so that $\vp_J(x,\xi)$ lies off all the finitely
	many local branches of $\Psi_J(L_J)$, and take
	$h\in C_c^\infty(U)$ with $dh_x=\xi$.

	Pull $\widetilde W$ back along $\coll_{J\top}$ to the collar
	$\Delta_{J\top}=\Delta_J\times[0,\varepsilon)$. In the product
	chart $\vp_{J\top}=\vp_J\times\mathrm{id}$, apply the fiberwise
	translation $\cc{T}(\chi(t)h)$ to the chosen branch only, where
	$\chi\in C_c^\infty([0,\varepsilon))$ equals $1$ near $0$.
	This modification preserves the collar coordinate $t$ and can
	be taken arbitrarily small and supported away from the vertical
	boundary and the horn, so it gives a Lagrangian cobordism relative
	to $\vlg$. On $\Delta_J$ the modified branch is
	$\vp_J\circ\Gamma h$. It covers $\vp_J(x,\xi)$ exactly once,
	and every other branch misses this point. Thus every
	diffeomorphism preserving the modified face immersion has a fixed
	point on each remaining component. Lemma~\ref{lem:fixed-locus}
	gives freeness on $\Delta_J$.
\end{proof}

\begin{proof}[Proof of Theorem~\ref{thm:hornFilling}]
	For $n=1$, take the degeneracy $s_{0,0}$ of the prescribed
	vertex. For $n\geq2$, Lemmas~\ref{lem:compact-horn-extension}
	and~\ref{lem:free-horn-extension} give an element
	$\widetilde W\in\cob(X,\vlg)_n$ with the prescribed
	restrictions. In both cases Lemma~\ref{lem:horn-filling}
	gives the desired extension.
\end{proof}

\subsection{The groupoid \texorpdfstring{$\cob(X,\vlg)$}{Cob(X,V)}}
\label{subsec:groupoid}

In this subsection we describe the homotopy category $h \cob(X,\vlg)$ explicitly. Objects are elements of $\cob(X,\vlg)_0 $, canonically identified with freely immersed compact Lagrangians in $ X $ with boundary in $\vlg$. We henceforth write $d_i$ and $s_i$ for the maps $d_{n,i}$ and $s_{n,i}$ defined in Lemma~\ref{lem:simplicial-set-cob}.

\begin{definition}
	Given $ W_{01} \in \cob(X,\vlg)_1 $ with
	\begin{equation}
		d_0(W_{01})=W_1,
		\quad
		d_1(W_{01})= W_0,
	\end{equation}
	we write $ W_{01}:W_0 \to W_1 $.
	For $ W_0 \in \cob(X,\vlg)_0 $, the \emph{identity cobordism} is
	\begin{equation}
		\label{eq:identityCob}
		\mathrm{id}_{W_0}: W_0 \to W_0,
		\quad
		\mathrm{id}_{W_0} := s_0(W_0) \in \cob(X,\vlg)_1,
	\end{equation}
	i.e. the pullback of $ W_0 $ along the constant map $ \bar\sigma_{0,0}: \Delta^1 \to \Delta^0 $, $ (s_0,s_1)\mapsto s_0+s_1 $.
\end{definition}

\begin{definition}
	Let $ W_0,W_1,W_2 \in \cob(X,\vlg)_0 $ and let
	\begin{equation}
		W_{01}: W_0 \to W_1,
		\quad
		W_{12}: W_1 \to W_2,
		\quad
		W_{02}: W_0 \to W_2,
	\end{equation}
	in $ \cob(X,\vlg)_1 $.  For $ \sigma \in \cob(X,\vlg)_2 $ we write
	\begin{equation}
		\sigma : W_{12} \circ W_{01} \Longrightarrow W_{02}
	\end{equation}
	if $ d_0(\sigma)=W_{12} $, $ d_1(\sigma)=W_{02} $ and $ d_2(\sigma)=W_{01} $:
	
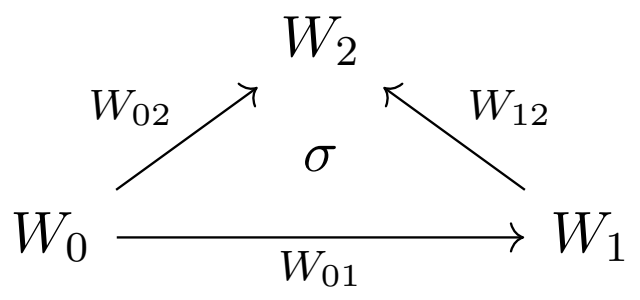
\begin{equation}
		\begin{tikzcd}
			&
			W_2
			&
			\\
			W_0
			\arrow[rr, "W_{01}"']
			\arrow[ru, "W_{02}"]
			& {}
			\arrow[u, "\sigma", phantom]
			& W_1 \arrow[lu, "W_{12}"']
		\end{tikzcd}
	\end{equation}
\end{definition}

\begin{lemma}
	\label{lem:degenerate2simplex}
	Let $ W_{01}: W_0 \to W_1 $. Then
	\begin{equation}
		s_{0}(W_{01}) : W_{01} \circ \mathrm{id}_{W_0} \Longrightarrow W_{01},
		\qquad
		s_{1}(W_{01}) : \mathrm{id}_{W_1} \circ W_{01}  \Longrightarrow W_{01}.
	\end{equation}
\end{lemma}
\begin{proof}
	Follows from the simplicial identities.
\end{proof}

\begin{definition}
	Let $ W_{01},W_{01}': W_0 \to W_1 $ in $ \cob(X,\vlg)_1 $. We say that $ W_{01}, W_{01}' $ are \emph{2-homotopic} if there exists a 2-cobordism $ \sigma \in \cob(X,\vlg)_2 $ such that $ \sigma: W_{01} \circ \mathrm{id}_{W_0} \Longrightarrow W_{01}' $:
	\begin{equation}
		\begin{tikzcd}
			&
			W_1
			&
			\\
			W_0
			\arrow[rr, "\mathrm{id}_{W_0}"']
			\arrow[ru, "W_{01}'"]
			&
			{} 	\arrow[u, "\sigma", phantom]
			&
			W_0
			\arrow[lu, "W_{01}"']
		\end{tikzcd}
	\end{equation}
	Write $ s_0 = s_0(W_{01}') $, $ s_1 = s_1(W_{01}') $ and consider the Lagrangian over $ {\Lambda_0[3]} $
	
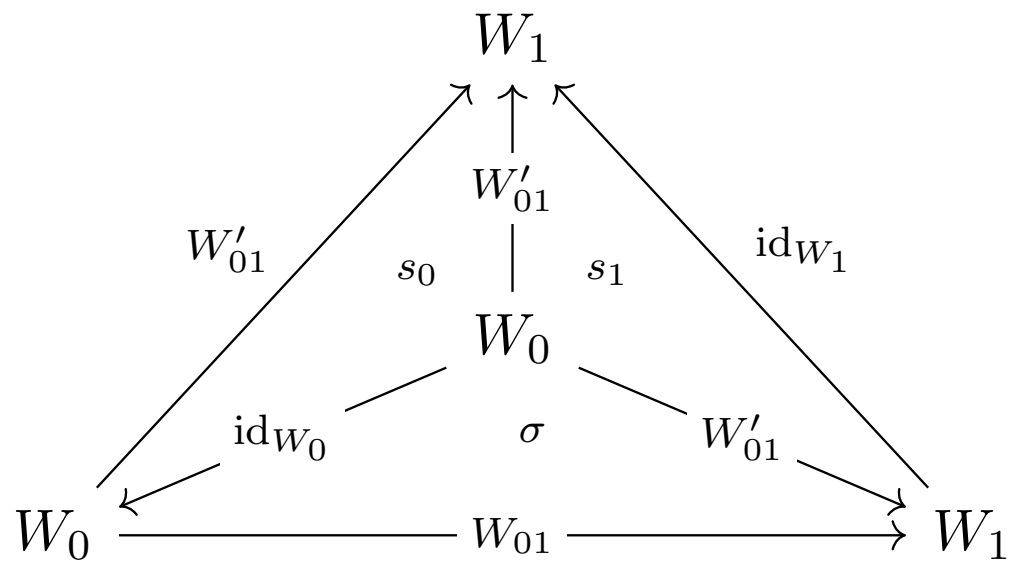
\begin{equation}
		\begin{tikzcd}
			& {}
			\arrow[rrrddd, "\scriptstyle s_1", phantom]
			&
			W_1
			& {} &
			\\
			&&&&
			\\
			&&
			W_0
			\arrow[uu, "W_{01}'" description]
			\arrow[rrd, "W_{01}'" description]
			\arrow[lld, "\mathrm{id}_{W_{0}}" description]
			&&
			{}
			\\
			W_0
			\arrow[rruuu, "W_{01}'"]
			\arrow[rrruuu, "\scriptstyle s_0", phantom]
			\arrow[rrrru, "\scriptstyle \sigma", phantom]
			\arrow[rrrr, "W_{01}"' description]
			&&&&
			W_1
			\arrow[lluuu, "\mathrm{id}_{W_{1}}"']
		\end{tikzcd}
	\end{equation}
	filling this 3-horn by Theorem \ref{thm:hornFilling},
	we  obtain a 2-simplex $ \sigma': \mathrm{id}_{W_1} \circ W_{01} \Longrightarrow W_{01}'$, which gives an equivalent criterion for 2-homotopy of $ W_{01},W_{01}' $.
\end{definition}

\begin{lemma}
	$ 2 $-homotopy is an equivalence relation.
\end{lemma}
\begin{proof}
	Let $ W_{01} \in \cob(X,\vlg)_1 $. Then $ s_0(W_{01}) : W_{01} \circ \mathrm{id}_{W_0}  \Longrightarrow W_{01} $ proves reflexivity. Let
	\begin{equation}
		\sigma:  W_{01} \circ \mathrm{id}_{W_0} \Longrightarrow W_{01}'.
	\end{equation}
	Set $ s=  s_0(\mathrm{id}_{W_0})$, $ s' = s_0(W_{01}') $ and consider the Lagrangian over $ {\Lambda_0[3]} $
	
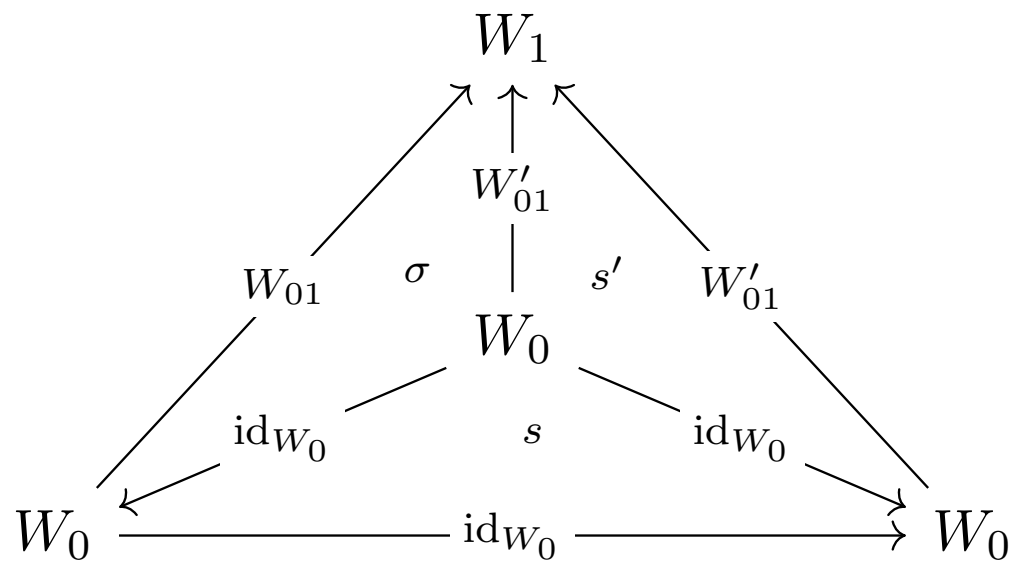
\begin{equation}
		\begin{tikzcd}
			& {}
			\arrow[rrrddd, "\scriptstyle s'", phantom]
			&
			W_1
			& {} &
			\\
			&&&&
			\\
			&&
			W_0
			\arrow[uu, "W_{01}'" description]
			\arrow[rrd, "\mathrm{id}_{W_0}" description]
			\arrow[lld, "\mathrm{id}_{W_0}" description]
			&&
			{}
			\\
			W_0
			\arrow[rruuu, "W_{01}" description]
			\arrow[rrruuu, "\scriptstyle \sigma", phantom]
			\arrow[rrrru, "\scriptstyle s", phantom]
			\arrow[rrrr, "\mathrm{id}_{W_0}"' description]
			&&&&
			W_0
			\arrow[lluuu, "W_{01}'" description]
		\end{tikzcd}
	\end{equation}
	Applying Theorem \ref{thm:hornFilling} we obtain $ \sigma' : W_{01}' \circ \mathrm{id}_{W_0} \Longrightarrow W_{01} $, proving symmetry. Let
	\begin{equation}
		\sigma: W_{01} \circ \mathrm{id}_{W_0} \Longrightarrow W_{01}',
		\qquad
		\tau:  \mathrm{id}_{W_1} \circ W_{01}'  \Longrightarrow W_{01}'',
	\end{equation}
	Write $ s = s_0(W_{01}'') $. Consider the Lagrangian over $ {\Lambda_0[3]} $

	
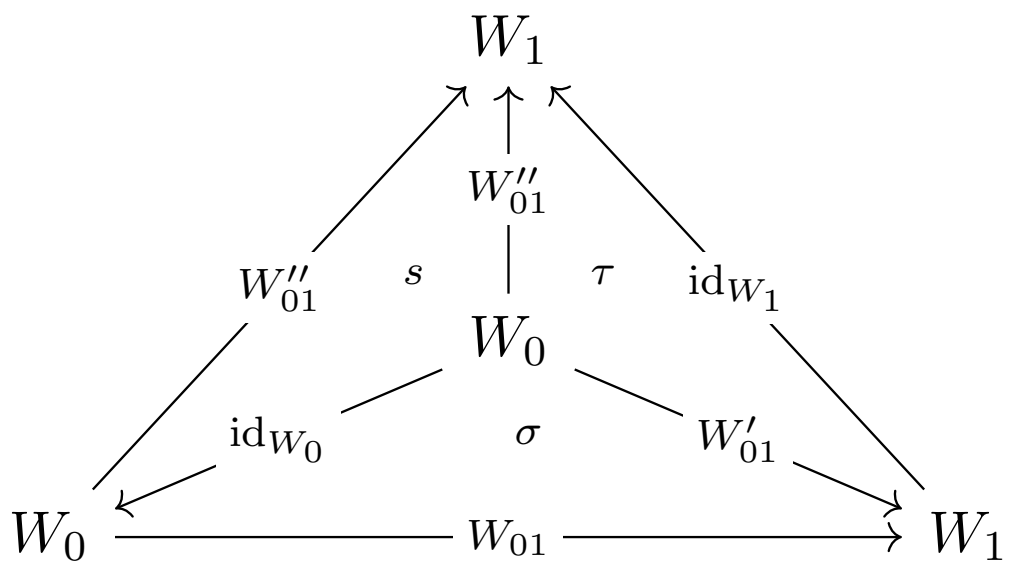
\begin{equation}
		\begin{tikzcd}
			& {}
			\arrow[rrrddd, "\scriptstyle \tau", phantom]
			&
			W_1
			& {} &
			\\
			&&&&
			\\
			&&
			W_0
			\arrow[uu, "W_{01}''" description]
			\arrow[rrd, "W_{01}'" description]
			\arrow[lld, "\mathrm{id}_{W_0}" description]
			&&
			{}
			\\
			W_0
			\arrow[rruuu, "W_{01}''" description]
			\arrow[rrruuu, "\scriptstyle s", phantom]
			\arrow[rrrru, "\scriptstyle \sigma", phantom]
			\arrow[rrrr, "W_{01}"' description]
			&&&&
			W_1
			\arrow[lluuu, "\mathrm{id}_{W_1}" description]
		\end{tikzcd}
	\end{equation}
	Applying Theorem \ref{thm:hornFilling} we obtain a 2-cobordism $  \mathrm{id}_{W_1} \circ W_{01}  \Longrightarrow W_{01}'' $, proving transitivity.
\end{proof}

\begin{lemma}
	Let
	\begin{equation}
		[W_{01}]: W_0 \to W_1,
		\quad \quad
		[W_{12}]:  W_1 \to W_2
	\end{equation}
	be 2-homotopy classes of Lagrangian 1-cobordisms $ W_{01}: W_0 \to W_1 $ and $ W_{12}: W_1 \to W_2 $, respectively.
	\begin{enumerate}
		\item There exists  $ W_{02}: W_0 \to W_2 $ and a 2-cobordism $ \sigma:  W_{12}\circ W_{01} \Longrightarrow W_{02}$.
		\item Define
		      \begin{equation}
			      \label{eq:classComposition}
			      [W_{12}] \circ [W_{01}] : = [W_{02}].
		      \end{equation}
		      Then the class $ [W_{12}] \circ [W_{01}] $ is independent of the choices of representatives $ W_{01}, W_{12}$  and $ W_{02} $.
	\end{enumerate}
\end{lemma}

\begin{proof}
	Suppose $ W_{01} \sim W_{01}' $ and let $\tau:  W_{01} \circ \mathrm{id}_{W_0} \Longrightarrow W_{01}' $. Let
	\begin{equation}
		\sigma': W_{12} \circ W_{01}' \Longrightarrow W_{02}',
		\qquad
		\sigma: W_{12} \circ W_{01} \Longrightarrow W_{02}.
	\end{equation}
	Consider the Lagrangian over $ {\Lambda_2[3]} $
	
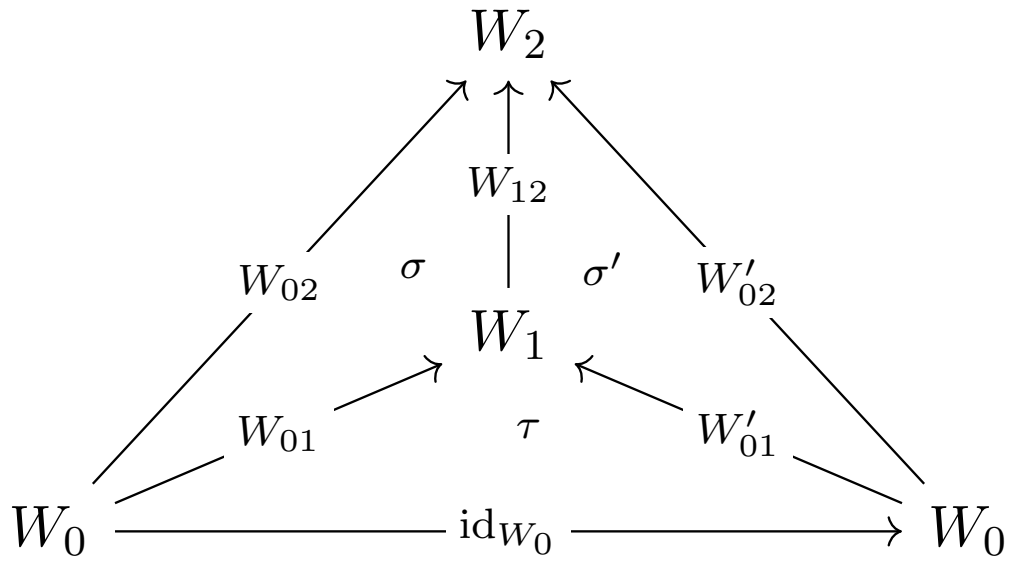
\begin{equation}
		\begin{tikzcd}
			& {}
			\arrow[rrrddd, "\scriptstyle \sigma'", phantom]
			&
			W_2
			& {} &
			\\
			&&&&
			\\
			&&
			W_1
			\arrow[uu, "W_{12}" description]
			&&
			{}
			\\
			W_0
			\arrow[rruuu, "W_{02}" description]
			\arrow[rrruuu, "\scriptstyle \sigma", phantom]
			\arrow[rrrru, "\scriptstyle \tau", phantom]
			\arrow[rrrr, "\mathrm{id}_{W_0}"' description]
			\arrow[rru, "W_{01}" description]
			&&&&
			W_0
			\arrow[llu, "W_{01}'" description]
			\arrow[lluuu, "W_{02}'" description]
		\end{tikzcd}
	\end{equation}
	Applying Theorem \ref{thm:hornFilling} we obtain equivalence $ W_{02} \sim W_{02}' $, hence the composition is independent of the choice of representative of $ [W_{01}] $. The other cases are analogous; In each case one builds the corresponding 3-horn and applies Theorem~\ref{thm:hornFilling} to obtain a 2-cobordism between the resulting outputs.
\end{proof}

\begin{prop}
	The simplicial set $ \cob(X,\vlg): \mathbf\Delta^\mathrm{op} \to \mathrm{Set} $ determines a groupoid $ h\cob(X,\vlg) $ with objects $ \cob(X,\vlg)_0 $, morphisms are 2-homotopy classes of 1-cobordisms in $ \cob(X,\vlg)_1 $, the identity morphism is represented by \eqref{eq:identityCob} and composition is given by \eqref{eq:classComposition}.
\end{prop}

\begin{proof}
	Left and right unit laws follow from Lemma \ref{lem:degenerate2simplex}. Let
	\begin{equation}
		[W_{01}]: W_0 \to W_1,
		\qquad
		[W_{12}]: W_1 \to W_2,
		\qquad
		[W_{23}]: W_2 \to W_3,
	\end{equation}
	and
	\begin{equation}
		W_{02} \in [W_{12}] \circ [W_{01}],
		\quad
		W_{13} \in  [W_{23}]\circ [W_{12}]
		\quad
		W_{03} \in \left( [W_{23}]\circ [W_{12}] \right) \circ [W_{01}].
	\end{equation}
	Consider the following Lagrangian complex over ${\Lambda_1[3]} $
	
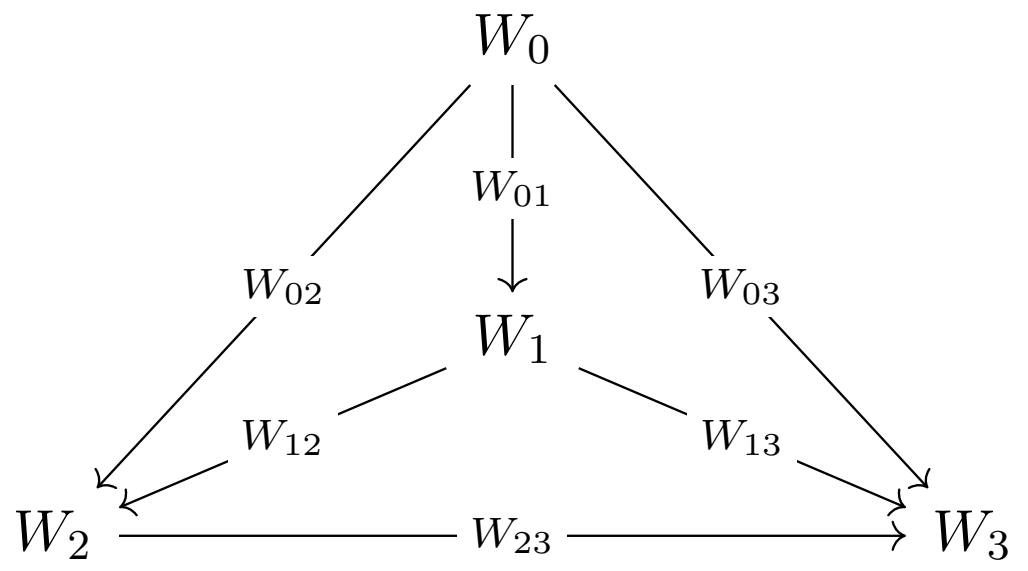
\begin{equation}
		\begin{tikzcd}
			& {}
			&
			W_0
			\arrow[llddd, "W_{02}" description]
			\arrow[rrddd, "W_{03}" description]
			\arrow[dd, "W_{01}" description]
			& {} &
			\\
			&&&&
			\\
			&&
			W_1
			\arrow[rrd, "W_{13}" description]
			\arrow[lld, "W_{12}" description]
			&&
			{}
			\\
			W_2
			\arrow[rrrr, "W_{23}"' description]
			&&&&
			W_3
		\end{tikzcd}
	\end{equation}
	Applying Theorem \ref{thm:hornFilling} yields a 2-cobordism between the two composites, hence associativity:
	\begin{equation}
		\left( [W_{23}]\circ [W_{12}] \right) \circ [W_{01}] = [W_{23}] \circ \left( [W_{12}] \circ [W_{01}] \right).
	\end{equation}
	Finally, any morphism $[W_{01}]: W_0 \to W_1$ admits an inverse: using Theorem \ref{thm:hornFilling}, choose $\sigma \in \cob(X,\vlg)_2$ filling the Lagrangian over ${\Lambda_0[2]}$
	
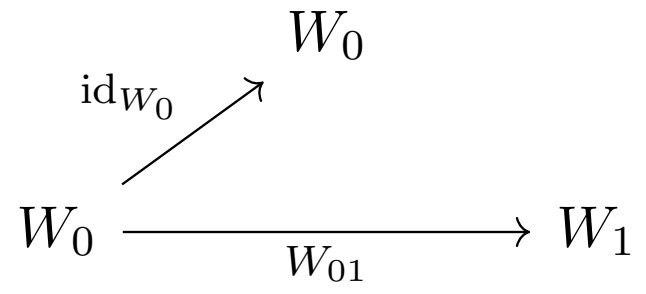
\begin{equation}
		\begin{tikzcd}
			&
			W_0
			&
			\\
			W_0
			\arrow[rr, "W_{01}"']
			\arrow[ru, "\mathrm{id}_{W_0}"]
			& {}
			& W_1
		\end{tikzcd}
	\end{equation}
	then $[d_0(\sigma)] $ is a left inverse of $[W_{01}]$. Similarly,
	using Theorem \ref{thm:hornFilling}, choose $\sigma \in \cob(X,\vlg)_2$ filling the Lagrangian over ${\Lambda_2[2]}$
    
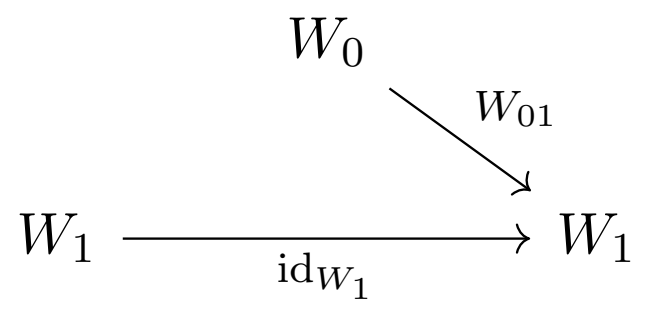
\begin{equation}
        \begin{tikzcd}
            & W_0 \arrow[dr, "W_{01}"] & \\
            W_1 \arrow[rr, "\mathrm{id}_{W_1}"'] & {} & W_1
        \end{tikzcd}
    \end{equation}
	then $[d_2(\sigma)]$ is a right inverse of $[W_{01}]$, hence $h\cob(X,\vlg)$ is a groupoid.
\end{proof}

\appendix

\section{Manifolds with corners}
\label{app: ManifoldsWithCorners}

\subsection{Basic definitions}

We recall the basic definitions of manifolds with corners following~\cite{Joyce12}.
An \emph{$ n $-dimensional manifold with corners} is a paracompact Hausdorff topological space locally modeled on open subsets of $  [0,\infty)^k \times \bb{R}^{n-k} $ for some $ 0 \leq k \leq n $.
\begin{definition}[{\cite[Definition~2.3]{Joyce12}}]
	Let $Y$ be an $n$-dimensional manifold with corners, and let  $p \in Y$. If $(U,\phi)$ is a coordinate chart with
	\begin{equation}
		\phi(p)=(u_1,\dots,u_k,t_{1},\dots,t_{n-k}) \in  [0,\infty)^k \times \bb{R}^{n-k},
	\end{equation}
	define the \emph{depth} of $p$ to be the number of coordinates among $u_1,\dots,u_k$ that are zero. This is independent of the choice of chart. The \emph{depth $k$ stratum of $Y$} is
	\begin{equation}
		S^k(Y) = \{ p \in Y \,\, | \,\,\mathrm{depth}(p)=k \},
	\end{equation}
	which is an $(n-k)$-dimensional manifold without boundary~\cite[Proposition~2.4]{Joyce12}.
\end{definition}

\begin{definition}[{\cite[Definition~2.5]{Joyce12}}]
	\label{def:boundary&corners}
	Let $Y$ be an $n$-dimensional manifold with corners. For $y\in Y$, a \textit{local boundary component} $\mathfrak{b}$ at $y$ is a choice of connected component of $S^1(Y)$ near $y$. The number of distinct local boundary components at $y$ is exactly $\text{depth}(y)$.
	Define
	\begin{equation}
		\partial Y := \{ (y,\mathfrak{b}) \,\, | \,\, y \in Y, \, \mathfrak{b} \text{ is a local boundary component for $Y$ at $y$} \}.
	\end{equation}
	Then $\partial Y$ is naturally an $(n-1)$-manifold with corners.  By iterating this construction we can define $\partial^2 Y, \dots , \partial^n Y$, where $\partial^k Y$ is an $(n-k)$-manifold with corners. By~\cite[Proposition 2.9]{Joyce12} there exists a canonical isomorphism
	\begin{equation}
		\partial^k Y
		\simeq
		\{ (y,\mathfrak{b}_1,\dots,\mathfrak{b}_k) \,\, | \,\, y \in Y, \, \mathfrak{b}_i \text{ are distinct l.b.c for $Y$ at $y$} \},
	\end{equation}
    where l.b.c stands for local boundary components.
	Let $ S_k$ be the symmetric group on  $ \{1,\dots,k\}$, acting on $\partial^k Y$ by permuting the boundary components. The \emph{$ k $-corners of $ Y $} are
	\begin{align}
		C_k(Y)
		 & := \partial^k Y/ S_k
		\\
		 & \simeq
		\{ (y,\{\mathfrak{b}_1,\dots,\mathfrak{b}_k\}) \,\, | \,\, y \in Y, \, \mathfrak{b}_i \text{ are distinct l.b.c for $Y$ at $y$} \},
	\end{align}
	which is an $ (n-k)$-manifold with corners.  Note that $Y \simeq C_0(Y)$ and $\partial Y \simeq C_1(Y)$. We denote by
	\begin{equation}
		\label{eq:cornersNaturalInclusion}
		\iota_{k}: C_k(Y) \to Y
	\end{equation}
	the canonical immersions given by forgetting the $k$ boundary components. Note that the immersions \eqref{eq:cornersNaturalInclusion} are generally non-injective.
\end{definition}

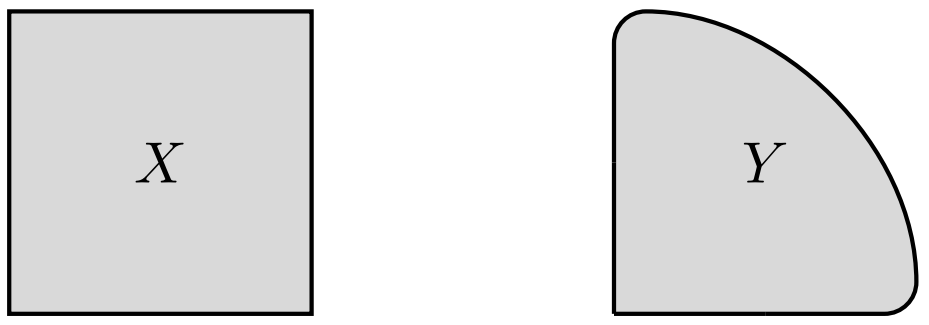
\begin{figure}[htbp]
	\centering
	\begin{tikzpicture}
		\def\R{2}

		\fill[gray!30, rounded corners=6pt]
		(4+\R,0) arc[start angle=0, end angle=90, radius=\R]
		-- (4,0) -- cycle;

		\fill[gray!30]
		(4,0) -- (4+\R/2,0)
		-- (4,\R/2) -- (4,0);

		\draw[thick, rounded corners=6pt]
		(4+\R/2,0) -- (4+\R,0) arc[start angle=0, end angle=90, radius=\R]
		-- (4,\R/2);
		\draw[thick] (4,0) -- (4+\R/2,0);
		\draw[thick] (4,0) -- (4,\R/2);

		\filldraw[fill=gray!30, thick] (0,0) rectangle (\R,\R);

		\node at (1,1) {$X$};
		\node at (5,1) {$Y$};
	\end{tikzpicture}
	\caption{For $X=$ square, $\partial X$ is a union of four closed intervals, $\partial^2 X$ is a set of $8$ points, and $C_2(X)$ is a set of four points. For $Y$, $\partial Y$ is a closed interval, $\partial^2 Y$ is a set of two points, and $C_2(Y)$ is a point.}
\end{figure}

\begin{definition}[{\cite[Definition~3.1]{Joyce12}}]
	\label{def:smoothMaps}
	Let $f:Y\to Z$ be a continuous map between manifolds with corners. The map $f$ is \emph{weakly smooth} if it is smooth in local coordinates. It is \emph{smooth} if, in addition, the following boundary condition holds.
	For any $y\in Y$, let $z=f(y)\in Z$, and let $\mathfrak b$ be a local boundary component of $Z$ at $z$. Choose a neighborhood $V$ of $z$ and a boundary defining function $b:V\to[0,\infty)$ for $\mathfrak b$ \cite[Definition~2.14]{Joyce12}. Then $f^{-1}(V)$ is open and $b\circ f:f^{-1}(V)\to[0,\infty)$ is weakly smooth. We require that either
	\begin{itemize}
		\item $b\circ f\equiv 0$ on some neighborhood of $y$, or
		\item there exists a \emph{unique} local boundary component $\tilde{\mathfrak b}$ of $Y$ at $y$ such that $b\circ f$ is a boundary defining function for $\tilde{\mathfrak b}$ at $(y,\tilde{\mathfrak b})\in\partial Y$.
	\end{itemize}
\end{definition}

\begin{definition}[{\cite[Definition~3.2]{Joyce12}}]
	\label{def:b_submersions}
	Let $f: Y \to Z$ be a smooth map of manifolds with corners.
	\begin{enumerate}
		\item $f$ is a \textit{submersion} if, for every $y \in Y$ with $y \in S^l(Y)$ and $f(y) \in S^k(Z)$, the differentials
		      \begin{equation}
			      df|_y: T_y Y \to T_{f(y)}Z,
			      \quad \quad
			      df|_y: T_y(S^l(Y)) \to T_{f(y)}(S^k(Z))
		      \end{equation}
		      are surjective.
		\item $f$ is a \textit{b-submersion} if, for every $y \in Y$ with $y \in S^l(Y)$ and $f(y) \in S^k(Z)$, the induced map on the boundary-normal quotients
		      \begin{equation}
			      (df|_y)_\ast: T_yY /T_y(S^l(Y)) \to T_{f(y)}Z/T_{f(y)}(S^k(Z))
		      \end{equation}
		      is surjective. In particular, every submersion is a b-submersion.
	\end{enumerate}
\end{definition}

\begin{lemma}
	\label{lem:b-submersive_to_submersion}
	Let $Y$ be a manifold with corners and let $\pi: Y \to [0,r)^k$ be a b-submersion. Suppose there exists $0<r_0<r $ such that $K_0:=\pi^{-1}([0,r_0]^k)$ is compact. Then, there exists $0< \varepsilon < r_0 $ such that
	\begin{equation}
		\pi|_{\pi^{-1}([0,\varepsilon)^k)}: \pi^{-1}([0,\varepsilon)^k) \to [0,\varepsilon)^k
	\end{equation}
	is a submersion.
\end{lemma}

\begin{proof}
	Let $y\in \pi^{-1}(0)$ with $y \in S^l(Y)$. Since $S^k([0,r)^k) = \{0\} $ and $\pi$ is a b-submersion, the map
	\begin{equation}
		(d\pi|_y)_\ast: T_y Y / T_y(S^l(Y)) \to T_0 ([0,r)^k)
	\end{equation}
	is surjective. In particular, $d\pi|_y : T_y Y \to T_0 ([0,r)^k)$ is surjective, and there exists a neighborhood $U_y \subset Y$ of $y$ such that $\pi|_{U_y}$ is a submersion. The compact $\pi^{-1}(0)$ admits a finite cover $U_1,\dots,U_m$ by such neighborhoods. Set $U = \bigcup_{i=1}^m U_i$.
	We claim there exists $0<\varepsilon<r_0$ such that $\pi^{-1}([0,\varepsilon)^k)\subset U$. If not, we could choose a sequence $y_n\in K_0\setminus U$ with $\pi(y_n) \to 0$. By compactness of $K_0$, a subsequence converges to some $y_0\in K_0\setminus U$, and by continuity $\pi(y_0)=0$, hence $y_0\in \pi^{-1}(0)\subset U$, a contradiction. Thus for some $\varepsilon\in(0,r_0)$ we have $\pi^{-1}([0,\varepsilon)^k)\subset U$, and since $\pi$ is a submersion on $U$, we are done.
\end{proof}

\begin{prop}[{\cite[Proposition~5.1]{Joyce12}}]
	\label{prop:submersion_local_product}
	Let $X,Y$ be manifolds with corners, $f: X \to Y$ a smooth submersion, and $x \in X$ with $f(x)=y$. Then there exist open neighborhoods $X'$, $Y'$ of $x,y$ in $X,Y$ with $f(X')=Y'$, a manifold $Z'$ with corners with $\dim X = \dim Y + \dim Z'$, and a diffeomorphism $X' \simeq Y ' \times Z'$ such that $f|_{X'}: X' \to Y'$ is identified with the projection $\pi_{Y'}: Y' \times Z' \to Y'$.
\end{prop}

\begin{definition}[{\cite[Definition~6.1]{Joyce12}}]
	\label{def:transversality}
	Let $X,Y,Z$ be manifolds with corners and $f : X \to Z$, $g : Y \to Z$ be smooth maps. We call $f,g$ \emph{transverse} if the following holds. Suppose $x \in X$, $y \in Y$ and $z \in Z$ with $f(x) = z = g(y)$, so that there are induced linear maps of tangent spaces $(df)_x : T_x X \to T_z Z$ and $(dg)_y : T_y Y \to T_z Z$. Let $x \in S^j(X)$, $y \in S^k(Y)$ and $z \in S^l(Z)$, so that as in Definition~\ref{def:b_submersions} we have $(df)_x\!\big(T_x S^j(X)\big) \subset T_z S^l(Z)$ and $(dg)_y\!\big(T_y S^k(Y)\big) \subset T_z S^l(Z)$. Then we require that
	\begin{align}
		\label{eq:transversality-tangent-spaces1}
		T_z Z      & =\; (df)_x\big(T_x X\big) \;+\; (dg)_y\big(T_y Y\big)
		\\
		\label{eq:transversality-tangent-spaces2}
		T_z S^l(Z) & =\; (df)_x\big(T_x S^j(X)\big) \;+\; (dg)_y\big(T_y S^k(Y)\big)
	\end{align}
	for all such $x,y,z$. From Definition~\ref{def:b_submersions}, if one of $f,g$ is a submersion then $f,g$ are automatically transverse.
\end{definition}

\begin{theorem}[{\cite[Theorem~6.4]{Joyce12}}]
	\label{thm:fiber_products_corners}
	Let $X,Y,Z$ be manifolds with corners, and let $f: X \to Z$, $g: Y \to Z$ be smooth transverse maps. Then the fiber product $X \times_{f,Z,g} Y$ exists in the category of manifolds with corners, and
	\[
		\dim\bigl(X \times_{f,Z,g} Y\bigr) \;=\; \dim X + \dim Y - \dim Z.
	\]
\end{theorem}

\begin{lemma}[{\cite[Lemma~4.4]{Joyce12}}]
	Let $f: Y \to Z$ be a b\mbox{-}submersion. Then the fiber product $Y \times_Z \partial Z$ is a manifold with corners, with a canonical embedding
	\begin{equation}
		\label{eq:horizontalEmbedding}
		\xi:\ Y \times_Z \partial Z \hookrightarrow \partial Y.
	\end{equation}
	We denote by $\partial_-^f Y$ the image of $\xi$, and by $\partial_+^f Y$ the complement $\partial Y \setminus \partial_-^f Y$.
\end{lemma}

\begin{definition}
	\label{def:horizontalAndVerticalBoundary}
	Let $f: Y \to Z$ be a b\mbox{-}submersion. We call $\partial_-^f Y$ the \emph{horizontal} boundary and $\partial_+^f Y$ the \emph{vertical} boundary.
\end{definition}

\begin{lemma}
	\label{lem:weakly-smooth-fiber-product}
	Let $g:Y\to Z$ be a b-submersion and let $f:X\to Z$ be weakly smooth.
	Suppose that $f$ and $g$ satisfy
	\eqref{eq:transversality-tangent-spaces1}
	and~\eqref{eq:transversality-tangent-spaces2}
	whenever $f(x)=g(y)$. Then $P=X\times_Z Y$, with the subspace
	topology, has the structure of a manifold with corners such that
	\begin{equation}
		\label{eq:fiber-product-tangent-space}
		T_{(x,y)}P
		=\{(\xi,\eta)\in T_xX\oplus T_yY
		\mid df_x(\xi)=dg_y(\eta)\}.
	\end{equation}
	The projection $P\to X$ is a b-submersion, and its local vertical
	boundary components at $(x,y)$ are in bijection with the
	local boundary components of $Y$ at $y$ that are vertical for $g$.
\end{lemma}

\begin{proof}
	Fix $(x_0,y_0)\in P$ and choose coordinates $(r,w)$ centered at
	$f(x_0)=g(y_0)$, where the $r_i$ are boundary defining functions
	and $w\in\mathbb R^m$. By~\cite[Lemma~4.4]{Joyce12}, the functions
	$r_i\circ g$ define distinct local boundary components of $Y$ at
	$y_0$. Complete $s=r\circ g$ to coordinates $(s,u)$ at $y_0$.
	Then $g(s,u)=(s,h(s,u))$, and the boundary coordinates among $u$
	define precisely the local vertical boundary components of $g$.
	Writing $f(x)=(R(x),F(x))$ and substituting $s=R(x)$ identifies
	$P$ locally, with its subspace topology, with the zero set of
	\[
		G(x,u)=h(R(x),u)-F(x).
	\]
	Since its target is Euclidean, $G$ is smooth
	by~\cite[Remark~3.3]{Joyce12}.
	Each component $r_i \circ f$ of $R$ is nonnegative and vanishes at $x_0$, so
	$dR_{x_0}$ vanishes on the tangent space to the stratum of $x_0$.
	Write $u_0=u(y_0)$. For $\xi$ tangent to the stratum of $x_0$
	and $(0,\eta)$ tangent to the stratum of $y_0$, we have
	\[
		dg_{y_0}(0,\eta)-df_{x_0}(\xi)
		=\bigl(0,d_u h_{(0,u_0)}(\eta)-dF_{x_0}(\xi)\bigr)
		=\bigl(0,dG_{(x_0,u_0)}(\xi,\eta)\bigr).
	\]
	By~\eqref{eq:transversality-tangent-spaces2}, the left-hand side
	ranges over the tangent space to the stratum of $f(x_0)$, which in
	the coordinates $(r,w)$ is $\{0\}\times\mathbb R^m$.
	Hence the restriction of $dG_{(x_0,u_0)}$ to the tangent space
	of the stratum of $(x_0,u_0)$ is surjective.

	The implicit function theorem
	therefore allows us to solve $G=0$ for $m$ real coordinates among the non-boundary coordinates in
	$(x,u)$. This gives a corner
	chart on $P$, whose coordinates are restrictions of coordinates on
	$X\times Y$. The implicit functions and the formula $s=R(x)$ recover all
	coordinates of $(x,y)$ from the chart coordinates. Evaluating the coordinate functions of another chart
	at this point therefore gives a weakly smooth transition map.
	Interchanging the charts gives the same conclusion for its inverse,
	so the charts are compatible by~\cite[Definition~2.1]{Joyce12}.

	In these charts the projection $P\to X$ pulls each boundary
	coordinate of $X$ back to a distinct boundary coordinate of $P$.
	It is therefore a b-submersion, and its local vertical boundary
	components are exactly those defined by the boundary coordinates
	among $u$. Finally, differentiating $f(x)=g(y)$ gives the inclusion
	of $T_{(x,y)}P$ in the right hand side of equation~\eqref{eq:fiber-product-tangent-space}. Both have dimension
	$\dim X+\dim Y-\dim Z$, by the chart construction and
	\eqref{eq:transversality-tangent-spaces1}, respectively.
\end{proof}

\begin{definition}[{\cite[Definition~4.5]{Joyce12}}]
	\label{def:C(f)}
	Let $f:Y\to Z$ be a b-submersion. Define
	\[
		\begin{aligned}
			C(f)\colon \coprod_{d=0}^{\dim Y} C_d(Y)
			&\longrightarrow \coprod_{h=0}^{\dim Z} C_h(Z),\\
			(y,\{\tilde{\mathfrak b}_1,\dots,\tilde{\mathfrak b}_d\})
			&\longmapsto (z,\{\mathfrak b_1,\dots,\mathfrak b_h\}),
		\end{aligned}
	\]
	where $z=f(y)$,  $\{\mathfrak b_1,\dots,\mathfrak b_h\}$ is 
	\[
        \Bigl\{\ \mathfrak b\ :\ (y,(z,\mathfrak b))\in Y\times_Z \partial Z,\ \ \xi\bigl(y,(z,\mathfrak b)\bigr)=(y,\tilde{\mathfrak b}_i)\ \text{for some }i\in\{1,\dots,d\}\ \Bigr\},
	\]
	and $\xi:Y\times_Z\partial Z\hookrightarrow \partial Y$ is the canonical embedding \eqref{eq:horizontalEmbedding}.
\end{definition}

\begin{remark}
	The map $C(f)$ can be defined for general smooth maps. The assignment $Y\mapsto \coprod_d C_d(Y)$, $f\mapsto C(f)$ has good functorial properties~\cite[Theorem~4.7]{Joyce12}. The map $C(f)$ is smooth in the generalized sense of Joyce, and if $f$ is a b\mbox{-}submersion then
	\(
	C(f)\big(C_d(Y)\big)\subset \coprod_{0\le h\le d} C_h(Z).
	\)
\end{remark}

\begin{theorem}[{\cite[Theorem~6.11]{Joyce12}}]
	Let $X,Y,Z$ be manifolds with corners and let $f:Y\to Z$, $g:X\to Z$ be smooth transverse maps with $f$ a b\mbox{-}submersion. Then there is a canonical diffeomorphism
    \begin{equation}
		\label{eq:fiberProdCorners}
		\begin{aligned}[t]
			C_d(X\times_Z Y)
			&\cong \coprod_{\substack{h,i,k\ge0\\ d=h+k-i}}
				\bigl(C(g)^{-1}(C_i(Z))\cap C_h(X)\bigr)\\
			&\qquad{}\times_{C_i(Z)}
				\bigl(C(f)^{-1}(C_i(Z))\cap C_k(Y)\bigr).
		\end{aligned}
	\end{equation}
\end{theorem}

\subsection{Horizontal and vertical corners}
\label{subsec:tensors}

\begin{definition}
	\label{def:submersion_mwc_tensors}
	Let $f:Y\to Z$ be a b\mbox{-}submersion and fix $d$ with $0\le d\le \dim Y$. For $h,v\ge 0$ with $h+v=d$, define $\tensor{C}{_h^{\,v}}(Y)$ by the points $(y,\{\tilde{\mathfrak b}_1,\dots,\tilde{\mathfrak b}_d\})\in C_d(Y)$
	for which exactly $h$ of the chosen local boundary components are
	\emph{horizontal} with respect to $f$.
	Equivalently,
	\begin{equation}
		\label{eq:inverseCornerFunction}
		\tensor{C}{_h^{\,v}}(Y)
		\;=\; C(f)^{-1}\!\bigl(C_h(Z)\bigr)\ \cap\ C_d(Y).
	\end{equation}
	Thus $f$ induces a decomposition
	\[
		C_d(Y)\;=\;\coprod_{h+v=d}\tensor{C}{_h^{\,v}}(Y).
	\]
	We denote by $\tensor{f}{_h^{\,v}}:\tensor{C}{_h^{\,v}}(Y)\to C_h(Z)$ the restriction $C(f)\big|_{\tensor{C}{_h^{\,v}}(Y)}$.
\end{definition}

\begin{definition}
	For $0 \le h \le d$ set $v:=d-h$ and introduce the alternate notation
	\begin{equation}
		\label{eq:checkTensor}
		\tensor{\Cd}{_d^{\,h}}(Y)\;:=\;\tensor{C}{_h^{\,v}}(Y),
		\qquad
		C_d(Y)\;=\;\coprod_{h=0}^{d}\tensor{\Cd}{_d^{\,h}}(Y).
	\end{equation}
	We denote by
	\[
		\tensor{\check f}{_d^{\,h}}
		\;:=\;
		C(f)\big|_{\tensor{\Cd}{_d^{\,h}}(Y)}
		\;:\;
		\tensor{\Cd}{_d^{\,h}}(Y)\longrightarrow C_h(Z).
	\]
	The ``check'' decoration highlights that in $\tensor{\Cd}{_d^{\,h}}$ the \emph{lower} index $d$ records the \emph{corner depth}, while the \emph{upper} index $h$ records the number of horizontal components. We use both notations in view of their different roles in Lemma~\ref{lem:fiberAndTensorProducts} below.
\end{definition}

\begin{example}
	Let $Z=X\times Y$ and let $\pi:Z\to X$ be the projection. Then, \emph{relative to $\pi$}, there is a canonical identification
	\[
		\tensor{C}{_h^{\,v}}(Z)\;\cong\; C_h(X)\times C_v(Y).
	\]
\end{example}

\begin{lemma}
	\label{lem:horizontalPullback}
	Let $f:Y\to Z$ be a b\mbox{-}submersion and let $\iota_h:C_h(Z)\to Z$ be the canonical immersion. Then there is a canonical diffeomorphism
	\[
		\tensor{C}{_h^{\,0}}(Y)\;\cong\; C_h(Z)\times_{\iota_h,Z,f} Y.
	\]
\end{lemma}

\begin{proof}
	The canonical immersion $\tensor{C}{_h^{\,0}}(Y)\to Y$ together with
	$\tensor{f}{_h^{\,0}}:\tensor{C}{_h^{\,0}}(Y)\to C_h(Z)$ induces a bijective smooth map
	\[
		(y,\{\tilde{\mathfrak b}_1,\dots,\tilde{\mathfrak b}_h\})\;\longmapsto\;
		\bigl((f(y),\{\mathfrak b_1,\dots,\mathfrak b_h\}),\,y\bigr)
	\]
	into $C_h(Z)\times_{\iota_h,Z,f} Y$, where the $\mathfrak b_i$ are the (unique) horizontal boundary components at $z=f(y)$ corresponding to the $\tilde{\mathfrak b}_i$ under $f$  (Definition~\ref{def:smoothMaps}).
\end{proof}

\begin{lemma}
	\label{lem:verticalPullback}
	Let $f: Y \to Z$, $g: X \to Z$ be smooth transverse maps such that $f$ is a b\mbox{-}submersion. There exists a canonical diffeomorphism
	\begin{equation}
		\tns{C}{_0^v}(X \times_{g,Z,f} Y) \simeq X \times_{g,Z,\tns{f}{_0^v}} \tns{C}{_0^v}(Y),
	\end{equation}
	where the left-hand side is taken with respect to the canonical projection $\pi_X: X\times_Z Y \to X$.
\end{lemma}

\begin{proof}
	We have
	\[
		\tns{C}{_0^v}(X \times_{Z} Y)  = C(\pi_X)^{-1}(X) \cap C_v(X\times_ZY)
	\]
	This intersection is precisely the component with $d=v$, $h=i=0$, $k=v$ in \eqref{eq:fiberProdCorners}, which gives the claimed identification.
\end{proof}

\begin{lemma}
	Let $f: Y \to Z$, $g: X \to Z$ be smooth transverse b-submersions. There exists a canonical isomorphism
	\begin{equation}
		\tns{C}{_h^v}(X \times_Z Y)
		\simeq
		\coprod_{0 \leq i \leq h }
		\tns{\Cd}{_h^i}(X) \times_{\tns{\check g}{_h^i},C_i(Z), \tns{f}{_i^v}} \tns{C}{_i^{v}}(Y),
	\end{equation}
	where the left-hand side is taken with respect to the canonical projection  $ \pi_X: X\times_Z Y \to X $.
\end{lemma}

\begin{proof}
	Set $d=h+v$. Then
	\[
		\tns{C}{_h^v}(X \times_Z Y)
		\;=\;
		C(\pi_X)^{-1}(C_h(X))\ \cap\  C_{d}(X \times_Z Y).
	\]
	By \eqref{eq:fiberProdCorners}, for triples $(h,i,k)$ with $d=h+k-i$ we have components
	\[
		\bigl(C(g)^{-1}(C_i(Z))\cap C_h(X)\bigr)\ \times_{C_i(Z)}\ \bigl(C(f)^{-1}(C_i(Z))\cap C_k(Y)\bigr).
	\]
	Fixing $h$ and $d=h+v$ yields $k=v+i$. Using the notations
	$\tensor{\Cd}{_h^{\,i}}(X)=C(g)^{-1}(C_i(Z))\cap C_h(X)$ and
	$\tensor{C}{_i^{\,v}}(Y)=C(f)^{-1}(C_i(Z))\cap C_{v+i}(Y)$
	gives the stated identification.
\end{proof}

\begin{definition}
	\label{def:associatedTensors}
	Let $f:Y\to Z$ be a b\mbox{-}submersion. We associate to $f$ the arrays
	\[
		\check C(Y)\;=\;\bigl(\,\tensor{\Cd}{_d^{\,h}}(Y)\,\bigr)_{\substack{0 \le d \le \dim Y\\ 0 \le h \le \dim Z}},
		\qquad
		C(Y)\;=\;\bigl(\,\tensor{C}{_h^{\,v}}(Y)\,\bigr)_{\substack{0 \le h \le \dim Z\\ 0 \le v \le \dim Y-h}}.
	\]
	We will informally refer to $C(Y)$ and $\check C(Y)$ as \emph{$(1,1)$-tensors of manifolds with corners}.
	If $f:Y\to Z$ and $g:X\to Z$ are smooth transverse b\mbox{-}submersions, define the contracted tensor
	\[
		\tensor{\bigl(\check C(X)\otimes_{C(Z)} C(Y)\bigr)}{_h^{\,v}}
		\;:=\;
		\coprod_{i=0}^{h}\ \tensor{\Cd}{_h^{\,i}}(X) \times_{\,C_i(Z)} \tensor{C}{_i^{\,v}}(Y).
	\]
\end{definition}

\begin{figure}[htbp]
	\centering
	\begin{tikzpicture}[scale=0.6]
		\coordinate (B1) at (-0.058,-0.97);
		\coordinate (B2) at (4.097,0.951);

		\coordinate (X1) at (-6,2);
		\coordinate (X2) at (-2,4);
		\coordinate (X3) at (-6,6);
		\coordinate (X4) at (-10,4);

		\coordinate (Y1) at (0,5);
		\coordinate (Y2) at (4,7);
		\coordinate (Y3) at (4,12);
		\coordinate (Y4) at (0,10);

		\coordinate (W1) at (-6,8);
		\coordinate (W2) at (-2,10);
		\coordinate (W3) at (-6,12);
		\coordinate (W4) at (-10,10);

		\coordinate (W1') at (-6,13);
		\coordinate (W2') at (-2,15);
		\coordinate (W3') at (-6,17);
		\coordinate (W4') at (-10,15);

		\filldraw[fill=gray!20, thick] (X1) -- (X2) -- (X3) -- (X4) -- cycle;
		\draw[thick, fill=gray!20] (Y1) -- (Y2) -- (Y3) -- (Y4) -- cycle;
		\draw[very thick] (B1) -- (B2);

		\fill[blue!20] (W3) -- (W4) -- (W4') -- (W3');
		\fill[blue!20] (W2) -- (W3) -- (W3') -- (W2');
		\filldraw[fill=red!20, draw=orange, fill opacity=1, very thick] (W1) -- (W2) -- (W3) -- (W4) -- cycle;
		\filldraw[fill=red!20, draw=orange, fill opacity=0.7, very thick] (W1') -- (W2') -- (W3') -- (W4') -- cycle;
		\fill[blue!20, fill opacity=0.35] (W1) -- (W1') -- (W4') -- (W4);
		\fill[blue!20, fill opacity=0.35] (W1) -- (W2) -- (W2') -- (W1');

		\draw[very thick, olive] (W1) -- (W1');
		\draw[very thick,olive] (W2) -- (W2');
		\draw[very thick,olive] (W3) -- (W3');
		\draw[very thick,olive] (W4) -- (W4');

		\fill[purple] (W1) circle (3pt);
		\fill[purple] (W2) circle (3pt);
		\fill[purple] (W4) circle (3pt);

		\fill[purple] (W1') circle (3pt);
		\fill[purple] (W2') circle (3pt);
		\fill[purple] (W3') circle (3pt);
		\fill[purple] (W4') circle (3pt);

		\node[left, font=\tiny] at ($(X1)!0.5!(X4)$) {$\tensor{\check C}{_1^1}(X)\;$};
		\node[above, font=\tiny] at ($(X3)!0.5!(X4)$) {$\tensor{\check C}{_1^0}(X)\;\;$};
		\node[below, font=\tiny] at (X1) {$\tensor{\check C}{_2^1}(X)$};
		\node[font=\tiny] at ($(X1)!0.5!(X3)$) {$\tensor{\check C}{_0^0}(X)$};

		\node[right, font=\tiny] at ($(Y2)!0.5!(Y3)$) {$\tensor{ C}{_1^0}(Y)$};
		\node[above, font=\tiny] at ($(Y3)!0.5!(Y4)$) {$\tensor{C}{_0^1}(Y)\;\;$};
		\node[right, font=\tiny] at (Y2) {$\tensor{C}{_1^1}(Y)$};
		\node[font=\tiny] at ($(Y1)!0.5!(Y3)$) {$\tensor{C}{_0^0}(Y)$};

		\node[below,font=\tiny] at ($(B1)!0.5!(B2)$) {$\;\;\tensor{C}{_0}(Z)$};
		\node[left,font=\tiny] at (B1) {$\;\tensor{C}{_1}(Z)$};

		\node[right, font=\tiny, text=blue] at ($(W4)!0.5!(W4')$) {$\tensor{C}{_1^0}(W)\;\;$};

		\node[left, font=\tiny, text=red] at ($(W1)!0.8!(W2)$) {$\tensor{C}{_0^1}(W)\;\;$};

		\node[right, font=\tiny, text=olive] at ($(W2)!0.5!(W2')$) {$\tensor{C}{_2^0}(W)$};
		\node[left, font=\tiny, text=orange] at ($(W3')!0.4!(W4')$) {$\tensor{C}{_1^1}(W)\;\;$};

		\node[right, font=\tiny, text=purple] at (W2') {$\tensor{C}{_2^1}(W)\;\;$};

		\draw[-latex] (2,4) -- (2,2) node[midway, right] {$f$};
		\draw[-latex] (-2,2) -- (-0.613,1.266) node[midway, above] {$g$};
	\end{tikzpicture}
	\caption{Example of a fiber product with tensor notation.}
	\label{fig: fiberTensorProduct}
\end{figure}
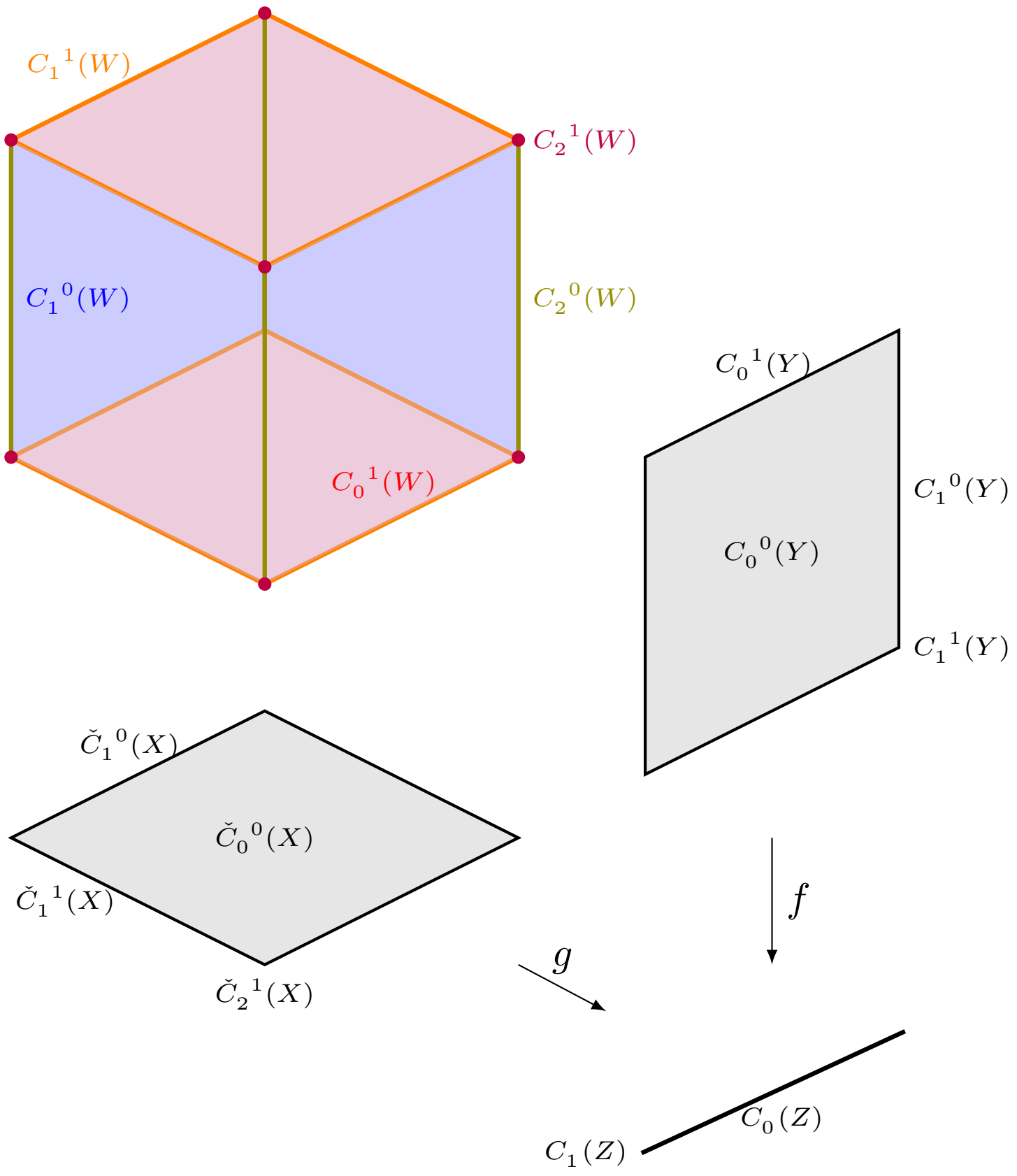

\begin{corollary}
	\label{lem:fiberAndTensorProducts}
	Let $f:Y\to Z$ and $g:X\to Z$ be smooth transverse b\mbox{-}submersions. Then
	\[
		C\bigl(X\times_Z Y\bigr)\;=\;\check C(X) \otimes_{C(Z)} C(Y),
	\]
	where $C(X\times_Z Y)$ is the $(1,1)$-tensor defined with respect to the b\mbox{-}submersion $\pi_X:X\times_Z Y\to X$.
\end{corollary}

\begin{example}
	Consider $X=Y=[0,1]^2$ and $Z=[0,1]$, with $f,g:[0,1]^2 \to [0,1]$ the projection $(x_1,x_2)\mapsto x_1$. Denote $W = X\times_Z Y$. Then:
	\begin{align}
		\tns{C}{_0^1}(W)
		 & = \tns{\Cd}{_0^0}(X) \times_{C_0(Z)} \tns{C}{_0^1}(Y) ,
		\\
		\tns{C}{_1^0}(W)
		 & = \left(\tns{\Cd}{_1^0}(X) \times_{C_0(Z)} \tns{C}{_0^0}(Y) \right)  \sqcup \left(\tns{\Cd}{_1^1}(X) \times_{C_1(Z)} \tns{C}{_1^0}(Y) \right),
		\\
		\tns{C}{_1^1}(W)
		 & = \left( \tns{\Cd}{_1^0}(X) \times_{C_0(Z)} \tns{C}{_0^1}(Y) \right)  \sqcup \left( \tns{\Cd}{_1^1}(X) \times_{C_1(Z)} \tns{C}{_1^1}(Y) \right),
		\\
		\tns{C}{_2^0}(W)
		 & = \tns{\Cd}{_2^1}(X) \times_{C_1(Z)} \tns{C}{_1^0}(Y) ,
		\\
		\tns{C}{_2^1}(W)
		 & = \tns{\Cd}{_2^1}(X) \times_{C_1(Z)} \tns{C}{_1^1}(Y).
	\end{align}
\end{example}

\bibliographystyle{amsalpha}
\bibliography{refs}

\end{document}